# Von Neumann Stability Analysis of DG-like and PNPM-like Schemes for PDEs that have Globally Curl-Preserving Evolution of Vector Fields

By

Dinshaw S. Balsara[1] (dbalsara@nd.edu) and Roger Käppeli[2] (roger.kaeppeli@sam.math.ethz.ch)

[1]Physics and ACMS Departments, University of Notre Dame

[2]Seminar for Applied Mathematics (SAM), Department of Mathematics, ETH Zurich, CH-8092, Zurich, Switzerland

**Abstract**

This paper examines a class of involution-constrained PDEs where some part of the PDE system evolves a vector field whose curl remains zero or grows in proportion to specified source terms. Such PDEs are referred to as curl-free or curl-preserving respectively. They arise very frequently in equations for hyperelasticity and compressible multiphase flow, in certain formulations of general relativity and in the numerical solution of Schrödinger's equation. Experience has shown that if nothing special is done to account for the curl-preserving vector field, it can blow up in a finite amount of simulation time. In this paper we catalogue a class of DG-like schemes for such PDEs. To retain the globally curl-free or curl-preserving constraints, the components of the vector field, as well as their higher moments, have to be collocated at the edges of the mesh. They are updated by using potentials that are collocated at the vertices of the mesh. The resulting schemes: 1) do not blow up even after very long integration times, 2) do not need any special cleaning treatment, 3) can operate with large explicit timesteps, 4) do not require the solution of an elliptic system and 5) can be extended to higher orders using DG-like methods. The methods rely on a special curl-preserving reconstruction and they also rely on multidimensional upwinding. The Galerkin projection, so crucial to the design of a DG method, is now carried out in the edges of the mesh and yields a weak form update that uses potentials that are obtained at the vertices of the mesh with the help of a multidimensional Riemann solver. A von Neumann stability analysis of the curl-preserving methods is carried out and the limiting CFL numbers of this entire family of methods is catalogued in this work. The stability analysis confirms that with increasing order of accuracy, our novel curl-free methods have superlative phase accuracy while substantially abstractreducing dissipation. We also show that PNPM-like methods, which only evolve the lower moments while reconstructing the higher moments, retain much of the excellent wave propagation characteristics of the DG-like methods while offering a much larger CFL number and lower computational complexity. The quadratic energy preservation of these methods is also shown to be excellent, especially at higher orders. The methods have also been shown to be curl-preserving over long integration times.

**I) Introduction**

Novel PDEs of importance to science and engineering problems are routinely discovered. Many of those PDE systems have involution constraints whose deeper study leads to mimetic numerical solution strategies that retain greater fidelity with the physics of those systems. While involution constraints can come in many forms, the last couple of years have seen the emergence of novel classes of PDEs that place constraints on the evolution of the curl of a vector field. The structure of the evolutionary equation for the vector field of interest is such that it either keeps the curl of the vector field zero, or evolves it in proportion to certain source terms. The former PDEs are referred to as curl-free whereas the latter class of PDEs are referred to as curl-preserving. Many of the hyperbolic systems resulting from the Godunov-Peshkov-Romenski (GPR) formulation for hyperelasticity and compressible multiphase flow with and without surface tension have such curl-preserving update equations (Godunov and Romenski [32], Romenski [40], Romenski *et al.* [41], Peshkov and Romenski [37], [38], Dumbser *et al.* [27], [28], [31], Schmidmayer *et al.* [42]). The equations of General Relativity when cast in the FO-CCZ4 formulation also have such a structure (Alic *et al.* [1], [2], 2012, Brown *et al.* [18], Dumbser *et al.* [29], Dumbser, *et al.* [30]). Similarly, it has recently become possible to recast Schrödinger's equation in first order hyperbolic form, and the time-evolution of this very important equation also has curl-preserving constraints (Dhaouadi *et al.* [25], Busto *et al.* [19]).

A motivating PDE system would help the reader a lot here. Let us take a simple example involving a fluid with thermal conduction in the GPR formulation. Let us denote the density by $\rho$, the fluid velocity by $\mathbf{v}$, the fluid pressure by "$P$", the fluid temperature by "$T$", the internal thermal energy density by "$e$", the total energy density by $E \equiv e + \rho \mathbf{v}^2/2$, the thermal impulse by a vector $\mathbf{J}$, the heat flux by a vector $\mathbf{q}$ and the thermal stress by the second rank tensor $\boldsymbol{\sigma}$. The equations for a fluid with thermal conduction can be written as

$$\frac{\partial \rho}{\partial t} + \nabla \cdot (\rho \mathbf{v}) = 0 \tag{1.1a}$$

$$\frac{\partial (\rho \mathbf{v})}{\partial t} + \nabla \cdot (\rho \mathbf{v} \otimes \mathbf{v} + P\mathbf{I} + \boldsymbol{\sigma}) = 0 \tag{1.1b}$$

$$\frac{\partial E}{\partial t} + \nabla \cdot \left( (E+P)\mathbf{v} + \mathbf{v} \cdot \boldsymbol{\sigma} + \mathbf{q} \right) = 0 \tag{1.1c}$$

$$\frac{\partial \mathbf{J}}{\partial t} + \nabla (\mathbf{J} \cdot \mathbf{v} + T) - \mathbf{v} \times (\nabla \times \mathbf{J}) = -\frac{\rho T}{\tau} \mathbf{J} \tag{1.1d}$$

The identity matrix is denoted by $\mathbf{I}$ in the above equations. To complete our description of the above system, we also mention the constitutive relation for the thermal stress tensor $\sigma_{ij} = \rho c_h^2 J_i J_j$ and the other constitutive relation for the thermal conduction vector $q_i = \rho T c_h^2 J_i$. Here $c_h$ denotes the hyperbolic speed of heat waves, i.e. the second sound. The first three of the four equations in eqn. (1.1) above reveal themselves to be the equations for mass, momentum and energy conservation for a fluid, with additional contributions from the thermal conduction vector, $\mathbf{q}$, and the thermal stress tensor, $\boldsymbol{\sigma}$. The fourth equation in eqn. (1.1) is a novel contribution from the GPR formulation, see (Romenski [40]).

Now let us focus on the last equation in eqn. (1.1). Let us consider the limit where the relaxation time is very large, so that the source term is irrelevant. Since the vector field $\mathbf{J}$ starts off curl-free, it is easy to see that it remains curl-free by considering the remaining two parts of that equation. The first part of the update equation, given by $\nabla (\mathbf{J} \cdot \mathbf{v} + T)$, is just the gradient of a scalar. Since the curl of a gradient is zero, the first term will not contribute to the curl if none is present initially. The second part of the update equation, given by $\mathbf{v} \times (\nabla \times \mathbf{J})$, will also be zero if the vector $\mathbf{J}$ is initially curl-free. We see, therefore, that the vector field $\mathbf{J}$ stays curl-free if it is initially curl-free in the limit of very large relaxation time. Of course, when the relaxation time cannot be ignored, the curl of the vector field does indeed evolve in response to the presence of the stiff source term $-\rho T\, \mathbf{J}/\tau$. It is important to realize that if the fourth equation in eqn. (1.1) does not have a consistent discretization then the curl of the vector field $\mathbf{J}$ will only be specified by the accuracy of the numerical method. As a result, even for regions of the flow that should have no thermal conduction, there will indeed continue to be some small amount of thermal conduction. This affects the fidelity of the method and its results. For this particular PDE system, Balsara *et al.* [16] showed via direct comparison that a curl-preserving formulation produces a desirable result whereas a direct zone-centered collocation of variables does not. The curl-preserving formulation did not require the inclusion of any additional equations and was able to operate with a robust CFL.

To consider another system, Dumbser, *et al*. [30] analyzed the FO-CCZ4 formulation of general relativity and showed that it has a large number of equations that look like eqn. (1.1d). The full PDE system is too large to detail here. They showed that zone-centered discretizations simply blow up! They found that the only way to mitigate the build-up of circulation in a zone-centered context consisted of adding one extra generalized Lagrange multiplier (GLM) system that evolved an additional vector field and another generalized Lagrange multiplier system that evolved an additional scalar field for *each and every* equation that looks like eqn. (1.1d). This more than doubles the computational cost; besides the mitigation is imperfect! Furthermore, the mitigation requires that the GLM fields propagate at speeds that are two or three times larger than the other signal speeds in the problem! As a result, the timestep is reduced by a factor of two or three. Through the description of the above two PDE systems, we see the importance of a consistent, curl-preserving discretization and evolution strategy.

A detailed pointwise description of the implementation of a finite volume scheme that is curl-preserving is also provided in Section III.5 of Balsara *et al*. [16]. The implementation consists of using a zone-centered reconstruction of the fluid variables in eqn. (1.1a,b,c) and a curl-preserving reconstruction of the vector field in eqn. (1.1d). This gives us higher order spatial accuracy. This is followed by the application of one-dimensional and two-dimensional Riemann problems at the faces and vertices of the mesh; which gives us upwinding. Well-known strong stability-preserving Runge-Kutta time stepping then provides higher order temporal accuracy.

We see, therefore, that many very useful PDE systems have a curl-preserving involution. The simplest example of such a PDE can be written as

$$\frac{\partial \mathbf{J}}{\partial t} + \nabla \left( \mathbf{J} \cdot \mathbf{v} + \varphi(\rho) \right) - \mathbf{v} \times \left( \nabla \times \mathbf{J} \right) = \mathbf{S}(\mathbf{J}, \rho) \quad (1.2)$$

Here $\mathbf{J}$ is the vector field of interest and $\mathbf{v}$ is some specified velocity field. By "$\rho$" we generically denote other variables that might be part of a larger set of equations to which eqn. (1.2) belongs. As a result, the potential $\phi \equiv \mathbf{J} \cdot \mathbf{v} + \varphi(\rho)$ can also depend on the other variables and the source term $\mathbf{S}(\mathbf{J}, \rho)$ depends on "$\mathbf{J}$" as well as "$\rho$". By setting the source term "$\mathbf{S}$" to zero and taking the curl of the above equation, it is easy to see that if we initially have $(\nabla \times \mathbf{J}) = 0$ then the structure of the above equation retains $(\nabla \times \mathbf{J}) = 0$ for all time; in other words, the evolution of the vector field $\mathbf{J}$ is curl-free. If the above equation has a non-zero source term, the evolution of the vector field $\mathbf{J}$ would be curl-preserving. In principle, any curl-preserving scheme should be able to reach the curl-free limit as the importance of the source term goes to zero. For the purposes of this study, it is adequate to study curl-free evolution of vector fields because the inclusion or exclusion of the source term does not change the discussion that is to follow. In other words, the inclusion of source terms requires a separate study of how stiff terms are included in the time-evolution of any hyperbolic PDE; and that is not the point of focus for this paper. Eqn. (1.2) shows up quite frequently as a part of many of the hyperbolic systems mentioned before. When eqn. (1.2) appears in such larger systems, it is usually the source of many numerical difficulties.

Eqn. (1.2), and the larger systems that include it, cannot be evolved with a traditional higher order Godunov methodology, even when we have accounted for non-conservative products. Indeed, Dumbser *et al.* [30] have shown that if a classical higher order Godunov scheme is applied to eqn. (1.2), the solution blows up very rapidly. The blow-up manifests itself in an explosive increase in the curl of the vector field when the simulation is run over long periods of time; please see Fig. 5 of Dumbser *et al.* [30] which shows the explosive blow-up of the curl of the vector field. The blow-up occurs even when the source terms are zero, indicating that the source terms are not the cause of the difficulty. The problem is not specific to the PDE system considered in that paper because it was shown in Boscheri *et al.* [17] that a very simple model system that is based on eqn. (1.2) will also blow up if treated with a classical higher order Godunov scheme. A generalized Lagrange multiplier-based (GLM) cleaning procedure has been developed in Dumbser *et al.* [30] for suppressing the fictitious numerical build up of the circulation, but it requires greatly increasing the signal speed that is used in the cleaning equations and also adds many more vector fields than are originally necessary. Indeed, Dumbser *et al.* [30] in their study of a general relativistic system had to use signal speeds in their GLM-style approach that were so large that they were substantially larger than the speed of light! It is unusual to introduce superluminal signal speeds in a theory that is based on recognizing the speed of light as the maximal signal speed. As a result, GLM cleaning-based schemes force a very strong reduction in timestep, which is inconsistent with our notions of how the timestep should evolve in a higher order Godunov scheme. Boscheri *et al.* [17] were able to obtain curl constraint-preservation but only at the expense of solving an elliptic PDE system at every timestep. Furthermore, the method was restricted to second order of accuracy. This too is inconsistent with our notion that a higher order Godunov scheme should be time-explicit, easy to solve and extensible to all higher orders.

In Balsara *et al*. [16] we first realized that the constrained evolution mandated by eqn. (1.2) requires a special treatment. Two innovations were introduced in that paper. First, a curl-constraint preserving reconstruction was proposed which requires us to start with a vector field whose components are collocated at the edges of the mesh, and are indeed aligned with the edge directions. Second, it was realized that a multidimensional Riemann solver that is invoked at the vertices of the mesh can indeed give us stabilization through multidimensional upwinding. When coupled with strong stability preserving Runge-Kutta (SSP-RK) timestepping (Shu and Osher [45], [46], Shu [47], Spiteri and Ruuth [43], [44], Gottlieb *et al*. [33]), we obtained a robust finite volume-based numerical scheme. As a result, both ingredients in the design of a higher order Godunov scheme – i.e. the reconstruction as well as the Riemann solver – had to be fundamentally rethought when dealing with eqn. (1.2). These two ingredients proved to be highly beneficial because Balsara *et al*. [16] were able to obtain stable finite volume-like schemes for systems that used eqn. (1.2) which: **1)** did not blow up even after very long integration times, **2)** did not need any GLM-style cleaning with its deleterious side-effect of needing very high signal speeds, **3)** could operate with large explicit timesteps, **4)** did not require the solution of an elliptic system and **5)** could be extended to higher orders by using WENO-like methods. The WENO-like methods draw on ideas from weighted essentially non-oscillatory schemes (Jiang and Shu [35], Balsara and Shu [3], Balsara, Garain and Shu [11]). As a result, Balsara *et al*. [16] were able to integrate non-linear hybridization into their novel curl-preserving reconstruction, thus making the curl-preserving reconstruction suitable for use with higher order Godunov schemes.

In Balsara *et al*. [16] curl-preserving WENO-like methods for evolving systems that used eqn. (1.2) were developed. We called such methods WENO-like because the reconstruction used many insights from WENO schemes, while being substantially different from traditional, finite volume-based WENO schemes. Recall that a conservation law has a flux form that ensures that the integrated conserved variables in any zone (or connected set of zones) evolve in response to the fluxes at the boundary of that volume. This telescoping property for the fluxes gives the discrete version of the conservation law a globally conservative property. In an entirely analogous fashion, the curl-free schemes presented in Balsara *et al*. [16] are such that the discrete circulation evaluated over the edges of any face (or collection of faces) depends only on the potentials at the vertices of that facial area. In that sense, the schemes developed in Balsara *et al*. [16] are globally curl-preserving because they have a telescoping property on the potentials. In that paper, we were also able to present a von Neumann stability analysis of WENO-like globally curl-preserving schemes, showing that with increasing order of accuracy the schemes became progressively less dissipative and their dispersion error was also reduced. But it is well-known (Reed and Hill [39], Cockburn and Shu [20], [22], [24], Cockburn *et al*. [21], [23], Liu *et al*. [36], Zhang and Shu [48]) that finite volume DG schemes have superior wave propagation properties relative to finite volume WENO schemes of comparable order. In their study of DG-like schemes for magnetohydrodynamics and Maxwell's equations that have one or more constrained, divergence-preserving vector fields, Balsara and Käppeli [12], [14] found a similar trend; and that trend was numerically confirmed in Hazra *et al*. [34] and Balsara *et al*. [15]. We may, therefore, expect that curl-preserving DG-like schemes for eqn. (1.2) should have wave propagation characteristics that are substantially better than curl-preserving WENO-like schemes for eqn. (1.2) at comparable orders of accuracy. This paper, therefore, has the following three goals:-

**1)** The *first goal* of this paper is to lay out the conceptual foundations for DG-like schemes that preserve the global curl constraint. Recall that a classical DG scheme starts with zone-centered mean values for the flow and endows it with higher moments. These higher moments are then evolved in time by a classical DG scheme. The time evolution of the higher moments is carried out in a fashion that is consistent with the governing equations. In an exactly analogous fashion, the globally curl-free DG-like schemes start with edge-centered components. These components are then endowed with higher moments whose time-evolution is carried out consistent with the governing equations.

**2)** The *second goal* of this paper is to carry out a von Neumann stability analysis of the newly-obtained curl-free DG-like schemes. We use this stability analysis to find the maximum CFL number that is available at all orders up to fourth order. We also use the stability analysis to show that our new class of DG-like schemes have superior wave propagation characteristics. The curl of a vector field only manifests itself in two or three dimensions. As a result, our von Neumann stability analysis is also two-dimensional. Because the stability analysis is multidimensional by necessity, it is not technically feasible to extend it to very high orders. The computer algebra systems that we use for carrying out this stability analysis cannot be pushed beyond fourth order.

**3)** Dumbser *et al*. [26] proposed PNPM schemes where all the modes that are up to $N^{th}$ degree were evolved while all the higher modes up to $M^{th}$ degree are reconstructed. The PNPM schemes had the great advantage that they displayed wave propagation properties that were almost as good as classical DG schemes of comparable order while permitting substantially larger timesteps. Balsara and Käppeli [12], [14] found a similar trend in their study of divergence-preserving PNPM-like schemes. Therefore, the *third goal* of this paper is to analyze PNPM-like schemes that are globally curl-preserving. We intend to show that such PNPM-like schemes are competitive with their DG-like counterparts at comparable order of accuracy. We also intend to show that the maximum CFL number of PNPM-like schemes is much larger than that of DG-like schemes.

In this paper we analyze eqn. (1.2) with $\phi(\rho)=0$ and $\mathbf{S}(\mathbf{J},\rho)=0$. A constant velocity "**v**" is specified. The plan of the paper is as follows. Section II introduces a DG-like formulation for a curl-preserving model equation. Section III provides the essential ideas behind curl-free and curl-preserving reconstruction while pointing the reader to the further literature. Section IV presents the von Neumann stability analysis. Section V presents results from the von Neumann analysis of globally curl-free DG-like schemes. Section VI presents some numerical results to show that the proposed schemes meet their order of accuracy. Section VII presents conclusions.

**II) DG-like Formulation for the Curl-Preserving Model Equation**

Let us consider eqn. (1.2) with the intention of bringing out its physics. This will lead us to a curl-preserving DG-like formulation for that equation. Let us first consider eqn. (1.2) in its curl-free form (i.e. with $\nabla \times \mathbf{J} = 0$) since the curl has to be kept mathematically zero in that limit. Since we are considering a mimimalist system, we can also take $\varphi(\rho)=0$ because we are ignoring any further PDEs in the system. We write the vector components in two-dimensions as $\mathbf{J}=(J^x, J^y)$ and assume a constant velocity $\mathbf{v}=(v_x, v_y)$, so that the dot product becomes $\mathbf{J}\cdot\mathbf{v} = v^x J^x + v^y J^y$. The equations that we have to solve can be written as

$$\frac{\partial J^x}{\partial t} + \frac{\partial \left(v^x J^x + v^y J^y\right)}{\partial x} = 0 \quad ; \quad \frac{\partial J^y}{\partial t} + \frac{\partial \left(v^x J^x + v^y J^y\right)}{\partial y} = 0 \quad \text{with the constraint} \quad \left(\frac{\partial J^y}{\partial x} - \frac{\partial J^x}{\partial y}\right) = 0$$

(2.1)

Using the constraint in each of the two equations, they can be written as

$$\frac{\partial J^x}{\partial t} + v^x \frac{\partial J^x}{\partial x} + v^y \frac{\partial J^x}{\partial y} = 0 \quad ; \quad \frac{\partial J^y}{\partial t} + v^x \frac{\partial J^y}{\partial x} + v^y \frac{\partial J^y}{\partial y} = 0 \tag{2.2}$$

In other words, eqn. (2.2) tells us that we should be able to advect two components of a vector field. *However, and this is very important*, we should be able to carry out this advection of the vector field in a fashion that preserves the curl-free aspect of the vector field over each zone. We can only carry out such an advection if the component $J^x$ is collocated at the x-edges of the mesh and the component $J^y$ is collocated at the y-edges of the mesh. Furthermore, the potential, $v^x J^x + v^y J^y$, should be collocated at the vertices of the mesh, as shown in Fig. 1. Furthermore, eqn. (2.2) shows us that the potential should be obtained via multidimensional upwinding at the vertices of the mesh. For this simple example, multidimensional upwinding can be enforced visually depending on the direction of the velocity field. For example, Fig. 1 shows the multidimensionally upwinded potentials in the case where both components of the velocity are positive. On a mesh with zone sizes $\Delta x$ and $\Delta y$ in the x- and y-directions, the curl-free update equations at first order (with positive velocity components) can be written as

$$\frac{\partial J^x_{i,j+1/2}}{\partial t} + \left[\frac{\left(v^x J^x_{i,j+1/2} + v^y J^y_{i+1/2,j}\right) - \left(v^x J^x_{i-1,j+1/2} + v^y J^y_{i-1/2,j}\right)}{\Delta x}\right] = 0 \quad ;$$

$$\frac{\partial J^y_{i+1/2,j}}{\partial t} + \left[\frac{\left(v^x J^x_{i,j+1/2} + v^y J^y_{i+1/2,j}\right) - \left(v^x J^x_{i,j-1/2} + v^y J^y_{i+1/2,j-1}\right)}{\Delta y}\right] = 0$$

(2.3)

Eqn. (2.3) shows us that the collocation described in Fig. 1 is crucial for maintaining a curl-free update and such a collocation will have to be built into our curl-free DG-like scheme. The use of identical potentials at each vertex of the mesh allows us to claim that the discrete version of the curl-free constraint:-

$$\frac{J^y_{i+1/2,j} - J^y_{i-1/2,j}}{\Delta x} - \frac{J^x_{i,j+1/2} - J^x_{i,j-1/2}}{\Delta y} = 0 \tag{2.4}$$

is preserved forever, and it is preserved for each and every zone of the mesh. In other words, the update equations are globally curl-free. Eqn. (2.3) also shows us the importance of multidimensional upwinding because it provides a unique potential at each vertex of the mesh. To prove that eqn. (2.4) follows from eqn. (2.3) at all orders of accuracy, please write out equations like eqn. (2.3) for $J^y_{i+1/2,j}$, $J^y_{i-1/2,j}$, $J^x_{i,j+1/2}$ and $J^x_{i,j-1/2}$ in Fig. 1. Please write the equations out in terms of the potentials $\phi_{i+1/2,j+1/2}$, $\phi_{i-1/2,j+1/2}$, $\phi_{i+1/2,j-1/2}$ and $\phi_{i-1/2,j-1/2}$ shown at the vertices in Fig. 1. Notice that these potentials can be made as accurate as we desire by the use of a higher order reconstruction of the vector field. The pairwise cancellation will show that eqn. (2.4) is satisfied at all orders – in other words, we have a mimetic scheme. More generally, when eqn. (1.2) is incorporated into a larger PDE system, the multidimensional Riemann solver (Balsara, [4], [5], [7], [9], Balsara, Dumbser and Abgrall [6], Balsara and Dumbser [8], Balsara *et al*. [10], Balsara and Nkonga [13]) provides us with multidimensional upwinding that is consistent with the waves that are propagating in all directions.

The process of designing a higher order DG-like scheme consists of endowing the primal variables with higher order moments and then evolving those moments consistent with the governing equations. Consider the zone (*i,j*) in Fig. 1. We center the coordinates at the zone center so that the two-dimensional zone has extent $[-\Delta x/2, \Delta x/2] \times [-\Delta y/2, \Delta y/2]$. At the right y-edge and the top x-edge of the zone we assert the moments up to fourth order as

$$J^y(y,t) = J^y_0(t) + J^y_y(t)\left(\frac{y}{\Delta y}\right) + J^y_{yy}(t)\left(\left(\frac{y}{\Delta y}\right)^2 - \frac{1}{12}\right) + J^y_{yyy}(t)\left(\left(\frac{y}{\Delta y}\right)^3 - \frac{3}{20}\left(\frac{y}{\Delta y}\right)\right) \ ;$$

$$J^x(x,t) = J^x_0(t) + J^x_x(t)\left(\frac{x}{\Delta x}\right) + J^x_{xx}(t)\left(\left(\frac{x}{\Delta x}\right)^2 - \frac{1}{12}\right) + J^x_{xxx}(t)\left(\left(\frac{x}{\Delta x}\right)^3 - \frac{3}{20}\left(\frac{x}{\Delta x}\right)\right)$$

(2.5)

In a DG-like scheme, the modes in eqn. (2.5) become time-evolutionary. If only the linear part of eqn. (2.5) is retained, we get a second order DG-like scheme; if the quadratic part of eqn. (2.5) is retained, we get a third order DG-like scheme; if all the terms of eqn. (2.5) are retained, we get a fourth order DG-like scheme. Eqn. (2.5) makes it easy to see the trial functions that are asserted in each of the edges of the mesh. We will use test functions that are identical to the trial functions. Let us first write the Galerkin projection in the abstract and then specialize it to eqn. (2.5). It is useful to realize that a traditional finite volume DG method is derived from a Gauss' law-based vector identity

$$\nabla \cdot (\psi \mathbf{F}) = \psi \nabla \cdot \mathbf{F} + \mathbf{F} \cdot \nabla \psi$$

From the one-dimensional gradients involved in eqn. (2.1) we realize that the DG-like methods that we seek depend on the product rule for derivatives applied one-dimensionally. In other words, we rely on the identity

$$\frac{\partial (\psi \phi)}{\partial x} = \psi \frac{\partial \phi}{\partial x} + \phi \frac{\partial \psi}{\partial x}$$

In the above equation, "$\psi$" is a test function that lives in the edges of the mesh. Applying the above identity to the y-component of eqn. (1.2) we can make the following Galerkin projection in the y-edge of the mesh as follows:-

$$\frac{\partial}{\partial t} \left( \int_{y=-\Delta y/2}^{\Delta y/2} \psi(y) J^y(y,t) dy \right) + \psi(y=\Delta y/2) \phi(y=\Delta y/2) - \psi(y=-\Delta y/2) \phi(y=-\Delta y/2)$$
$$- \left( \int_{y=-\Delta y/2}^{\Delta y/2} \psi'(y) \phi(y) dy \right) + \left( \int_{y=-\Delta y/2}^{\Delta y/2} \psi(y) v_x (\nabla \times \mathbf{J})_z dy \right) = \left( \int_{y=-\Delta y/2}^{\Delta y/2} \psi(y) S_y (\mathbf{J}, \rho) dy \right) \quad (2.6)$$

Here $\psi(y)$ is a test function that lives in the y-edge of Fig. 1. Similarly, applying the above identity to the x-component of eqn. (1.2) we get

$$\frac{\partial}{\partial t} \left( \int_{x=-\Delta x/2}^{\Delta x/2} \psi(x) J^x(x,t) dx \right) + \psi(x=\Delta x/2) \phi(x=\Delta x/2) - \psi(x=-\Delta x/2) \phi(x=-\Delta x/2)$$
$$- \left( \int_{x=-\Delta x/2}^{\Delta x/2} \psi'(x) \phi(x) dx \right) - \left( \int_{x=-\Delta x/2}^{\Delta x/2} \psi(x) v_y (\nabla \times \mathbf{J})_z dx \right) = \left( \int_{x=-\Delta x/2}^{\Delta x/2} \psi(x) S_x (\mathbf{J}, \rho) dx \right) \quad (2.7)$$

There is not much going on in eqns. (2.6) and (2.7) other than an integration by parts along with a one-dimensional Galerkin projection. However, in the next paragraph we interpret the above two equations in order to bring out the physics of the situation.

Eqn. (2.6) allows us to write the update equations for the evolutionary modes of $J^y(y,t)$ as

$$\frac{dJ_0^y(t)}{dt} + \frac{1}{\Delta y}\left[\phi^{**}(y=\Delta y/2) - \phi^{**}(y=-\Delta y/2)\right] + \langle v_x(\nabla\times\mathbf{J})_z\rangle = \langle S_y(\mathbf{J}^*,\rho^*)\rangle \qquad (2.8a)$$

$$\frac{1}{12}\frac{dJ_y^y(t)}{dt} + \frac{1}{2\Delta y}\left[\phi^{**}(y=\Delta y/2) + \phi^{**}(y=-\Delta y/2)\right] - \frac{1}{\Delta y}\langle\phi^*(y)\rangle$$
$$+ \left\langle\left(\frac{y}{\Delta y}\right)v_x(\nabla\times\mathbf{J})_z\right\rangle = \left\langle\left(\frac{y}{\Delta y}\right)S_y(\mathbf{J}^*,\rho^*)\right\rangle \qquad (2.8b)$$

$$\frac{1}{180}\frac{dJ_{yy}^y(t)}{dt} + \frac{1}{6\Delta y}\left[\phi^{**}(y=\Delta y/2) - \phi^{**}(y=-\Delta y/2)\right] - \frac{2}{\Delta y}\left\langle\left(\frac{y}{\Delta y}\right)\phi^*(y)\right\rangle$$
$$+ \left\langle\left(\left(\frac{y}{\Delta y}\right)^2 - \frac{1}{12}\right)v_x(\nabla\times\mathbf{J})_z\right\rangle = \left\langle\left(\left(\frac{y}{\Delta y}\right)^2 - \frac{1}{12}\right)S_y(\mathbf{J}^*,\rho^*)\right\rangle$$

$$(2.8c)$$

$$\frac{1}{2800}\frac{dJ_{yyy}^y(t)}{dt} + \frac{1}{20\Delta y}\left[\phi^{**}(y=\Delta y/2) + \phi^{**}(y=-\Delta y/2)\right] - \frac{3}{\Delta y}\left\langle\left(\left(\frac{y}{\Delta y}\right)^2 - \frac{1}{20}\right)\phi^*(y)\right\rangle$$
$$+ \left\langle\left(\left(\frac{y}{\Delta y}\right)^3 - \frac{3}{20}\left(\frac{y}{\Delta y}\right)\right)v_x(\nabla\times\mathbf{J})_z\right\rangle = \left\langle\left(\left(\frac{y}{\Delta y}\right)^3 - \frac{3}{20}\left(\frac{y}{\Delta y}\right)\right)S_y(\mathbf{J}^*,\rho^*)\right\rangle \qquad (2.8d)$$

The angled brackets, $\langle\,\rangle$, represent line integrated averages of sufficiently high order within a y-edge (and later, similarly, for the x-edge). The potentials $\phi^{**}$ with the double star superscripts denote the potentials that are obtained at the vertices of the mesh using a multidimensional Riemann solver. The potentials $\phi^*(y)$ with the single star superscripts denote potentials that are obtained by application of one-dimensional Riemann solvers at the y-edge of the zone being considered. These one-dimensional Riemann solvers may be invoked at multiple quadrature points in the y-edge so that the terms $\langle\phi^*(y)\rangle$ and $\left\langle\left(\frac{y}{\Delta y}\right)\phi^*(y)\right\rangle$ are accurately evaluated. The source terms have a similar interpretation so that the $\mathbf{J}^*$ and $\rho^*$ variables in $S_y(\mathbf{J}^*,\rho^*)$ are obtained from the one-dimensional Riemann solvers. We see from eqn. (2.8a) that the update of the mean value, $J_0^y(t)$, will be curl-preserving and will be able to approach curl-free evolution in the limit where the source term tends to zero. Eqns. (2.8b), (2.8c) and (2.8d) show the same type of body terms that arise in a classical DG scheme due to the integration by parts with a test function; with the key difference that they are now applied to the edges of the mesh. The terms involving $(\nabla \times \mathbf{J})_z$ in eqn. (2.8) also have a special interpretation. These terms are exactly zero when the evolution is curl-free. When the evolution is only curl-preserving, these terms will be proportional to the discrete circulation around the zone, but only if a curl-preserving reconstruction from Balsara *et al.* [16] is used. From each of the two sides of a two-dimensional mesh, we can obtain terms that provide $(\nabla \times \mathbf{J})_z$. The resulting $\langle v_x (\nabla \times \mathbf{J})_z \rangle$ in eqn. (2.8a) is therefore an arithmetic average of the curl evaluated from either side of that edge. We make analogous interpretations for the terms with $(\nabla \times \mathbf{J})_z$ in eqns. (2.8b), (2.8c) and (2.8d).

Eqn. (2.7) allows us to write the update equations for the evolutionary modes of $J^x(x,t)$ as

$$\frac{dJ_0^x(t)}{dt} + \frac{1}{\Delta x}\left[\phi^{**}(x=\Delta x/2) - \phi^{**}(x=-\Delta x/2)\right] - \langle v_y (\nabla \times \mathbf{J})_z \rangle = \langle S_x(\mathbf{J}^*, \rho^*) \rangle \quad (2.9\text{a})$$

$$\frac{1}{12}\frac{dJ_x^x(t)}{dt} + \frac{1}{2\Delta x}\left[\phi^{**}(x=\Delta x/2) + \phi^{**}(x=-\Delta x/2)\right] - \frac{1}{\Delta x}\langle \phi^*(x) \rangle$$
$$- \left\langle \left(\frac{x}{\Delta x}\right) v_y (\nabla \times \mathbf{J})_z \right\rangle = \left\langle \left(\frac{x}{\Delta x}\right) S_x(\mathbf{J}^*, \rho^*) \right\rangle \quad (2.9\text{b})$$

$$\frac{1}{180}\frac{dJ_{xx}^x(t)}{dt} + \frac{1}{6\Delta x}\left[\phi^{**}(x=\Delta x/2) - \phi^{**}(x=-\Delta x/2)\right] - \frac{2}{\Delta x}\left\langle \left(\frac{x}{\Delta x}\right)\phi^*(x) \right\rangle$$
$$- \left\langle \left(\left(\frac{x}{\Delta x}\right)^2 - \frac{1}{12}\right) v_y (\nabla \times \mathbf{J})_z \right\rangle = \left\langle \left(\left(\frac{x}{\Delta x}\right)^2 - \frac{1}{12}\right) S_x(\mathbf{J}^*, \rho^*) \right\rangle \quad (2.9\text{c})$$

$$\frac{1}{2800}\frac{dJ_{xxx}^x(t)}{dt} + \frac{1}{20\Delta x}\left[\phi^{**}(x=\Delta x/2) + \phi^{**}(x=-\Delta x/2)\right] - \frac{3}{\Delta x}\left\langle \left(\left(\frac{x}{\Delta x}\right)^2 - \frac{1}{20}\right)\phi^*(x) \right\rangle$$
$$- \left\langle \left(\left(\frac{x}{\Delta x}\right)^3 - \frac{3}{20}\left(\frac{x}{\Delta x}\right)\right) v_y (\nabla \times \mathbf{J})_z \right\rangle = \left\langle \left(\left(\frac{x}{\Delta x}\right)^3 - \frac{3}{20}\left(\frac{x}{\Delta x}\right)\right) S_x(\mathbf{J}^*, \rho^*) \right\rangle \quad (2.9\text{d})$$

The interpretation of the terms in eqn. (2.9) mirrors that of eqn. (2.8). In eqn. (2.9), the angled brackets, $\langle \ \rangle$, represent line integrated averages of sufficiently high order within an x-edge.

In this work, we are interested in analyzing curl-free evolution, with the result that all terms with $(\nabla \times \mathbf{J})_z$, $S_x(\mathbf{J}^*, \rho^*)$ and $S_y(\mathbf{J}^*, \rho^*)$ can be set to zero in eqns. (2.8) and (2.9). Without the support of a larger PDE system, it is not possible to specify the source terms. It is, nevertheless, important for the sake of completeness of our discussion that the reader should understand what a curl-preserving reconstruction is. We illustrate that for the simplest of cases in the next Section.

**III) Curl-Preserving Reconstruction**

Eqns. (2.8) and (2.9) show that we need a reconstruction strategy within a zone that matches the vector components and their higher moments in the edges of the mesh. This is needed because we want the scheme to be globally curl-preserving. Consequently, each edge, as seen by its abutting zones, will have the same component of the vector field as well as its higher moments. Eqn. (1.2), as well as eqns. (2.8) and (2.9) show that when the discrete circulation evaluated around a zone is small, it should make proportionately small contributions to the edges via the $(\nabla \times \mathbf{J})_z$-dependent terms. In other words, the curl of the reconstructed vector field should match the discrete circulation (evaluated around each zone) as well as its higher moments. Fig. 2 shows us an example of how this works in two dimensions and at second order. In Balsara *et al*. [16] we have presented two and three-dimensional versions of such a reconstruction strategy at several orders. Here we just show some details for the second order case so that the reader may appreciate the core ideas as they are presented in one self-contained place.

Consider the vector field that is shown in Fig. 2. We consider the zone to be a unit square spanning $[-1/2, 1/2] \times [-1/2, 1/2]$. The modes in the left and right y-edges are given by

$$J_y^1 + (\Delta_y J_y^1) y \quad \text{and} \quad J_y^2 + (\Delta_y J_y^2) y \tag{3.1}$$

The modes in the bottom and top x-edges are given by

$$J_x^1 + (\Delta_x J_x^1) x \quad \text{and} \quad J_x^2 + (\Delta_x J_x^2) x \tag{3.2}$$

We can write the equations of a vector field which matches the values of the components and their linear variation within each edge as follows:-

$$J^x(x,y) = \left[J_x^1 + (\Delta_x J_x^1)x\right]\left(\frac{1}{2} - y\right) + \left[J_x^2 + (\Delta_x J_x^2)x\right]\left(\frac{1}{2} + y\right) + a_{yy}(1 - 4y^2) \ ;$$
$$J^y(x,y) = \left[J_y^1 + (\Delta_y J_y^1)y\right]\left(\frac{1}{2} - x\right) + \left[J_y^2 + (\Delta_y J_y^2)y\right]\left(\frac{1}{2} + x\right) + b_{xx}(1 - 4x^2)$$
(3.3)

Note that the coefficients $a_{yy}$ and $b_{xx}$ are needed for curl constraint-satisfaction. The discrete circulation within the zone of interest is given by $\left[J_x^1 - J_x^2 + J_y^2 - J_y^1\right]$, with the result that we can obtain the higher modes of that variable up to linear variation and write it as

$$R^z(x,y) = \left[J_x^1 - J_x^2 + J_y^2 - J_y^1\right] + (\Delta_x R^z)x + (\Delta_y R^z)y \tag{3.4}$$

Notice that if the discrete circulation $\left[J_x^1 - J_x^2 + J_y^2 - J_y^1\right]$ is zero, its variation will also be zero, so that we will get $(\Delta_x R^z)$ and $(\Delta_y R^z)$ as zero values. If the discrete circulation $\left[J_x^1 - J_x^2 + J_y^2 - J_y^1\right]$ is small, its variation as represented by $(\Delta_x R^z)$ and $(\Delta_y R^z)$ will also be proportionately small. We want to fix the coefficients $a_{yy}$ and $b_{xx}$ so that the curl of the vector field in eqn. (3.3) exactly matches eqn. (3.4). This is obtained by setting

$$b_{xx} = \frac{1}{8}\left[-(\Delta_x R^z) + (\Delta_x J_x^1) - (\Delta_x J_x^2)\right] \quad ; \quad a_{yy} = \frac{1}{8}\left[(\Delta_y R^z) + (\Delta_y J_y^1) - (\Delta_y J_y^2)\right]$$
(3.5)

This shows that the reconstructed vector field in eqn. (3.3) has been reconstructed in curl-preserving fashion. Therefore, the growth of the curl in eqns. (2.8) and (2.9) is perfectly well-controlled and consistent with the PDE in eqn. (1.2). Furthermore, when the discrete circulation is exactly zero, the $(\nabla \times \mathbf{J})_z$-dependent terms in eqns. (2.8) and (2.9) contribute absolutely nothing. In other words, the curl-free limit is exactly retrieved by our choice of reconstruction and discretization.

We can use eqn. (3.3), taken along with eqn. (3.5), to write the reconstructed curl-preserving vector field in terms of an orthonormal set of modes that span the zone itself. Projecting eqn. (3.3) into an orthonormal basis made of tensor product Legendre polynomials, we get

$$J^x(x,y) = \left[\left(J_x^1 + J_x^2\right)/2 + \left(\left(\Delta_y R^z\right) + \left(\Delta_y J_y^1\right) - \left(\Delta_y J_y^2\right)\right)/12\right] + \left[\left(\left(\Delta_x J_x^1\right) + \left(\Delta_x J_x^2\right)\right)/2\right]x + \left[-J_x^1 + J_x^2\right]y$$
$$+ \left[-\left(\left(\Delta_y J_y^1\right) - \left(\Delta_y J_y^2\right) + \left(\Delta_y R^z\right)\right)/2\right]\left(y^2 - 1/12\right) + \left[-\left(\Delta_x J_x^1\right) + \left(\Delta_x J_x^2\right)\right]xy$$
$$J^y(x,y) = \left[\left(J_y^1 + J_y^2\right)/2 + \left(-\left(\Delta_x R^z\right) + \left(\Delta_x J_x^1\right) - \left(\Delta_x J_x^2\right)\right)/12\right] + \left[-J_y^1 + J_y^2\right]x + \left[\left(\left(\Delta_y J_y^1\right) + \left(\Delta_y J_y^2\right)\right)/2\right]y$$
$$+ \left[-\left(\left(\Delta_x J_x^1\right) - \left(\Delta_x J_x^2\right) - \left(\Delta_x R^z\right)\right)/2\right]\left(x^2 - 1/12\right) + \left[-\left(\Delta_y J_y^1\right) + \left(\Delta_y J_y^2\right)\right]xy$$

(3.6)

Notice that eqn. (3.6) shows us that there is a transcription from the modes that we use in a curl-preserving DG-like scheme to the modes that we would use for a finite volume-based DG scheme. We see that the modes of a curl-preserving DG-like scheme just combine differently, consistent with the constraints, to give us the modes in a traditional, finite volume-based DG scheme. This ensures that the order property is always retained. Note though that curl constraint-preservation results in some higher order modes in eqn. (3.6) that would not be present in a classical second order finite volume-based DG scheme. Therefore, the reverse transcription, i.e. going from the modes of a finite volume-based DG scheme of a certain order to the modes of a curl-preserving DG-like scheme of the same order, does not hold. At second order, one cannot make much from this transcription because all the coefficients in eqn. (3.6) are fully determined. However, as shown in Section II.4 of Balsara *et al*. [16], at fourth order and beyond, some of the modes of the higher order curl-preserving reconstruction have to be obtained volumetrically while others are obtained from the edges of the mesh.

**IV) von Neumann Stability Analysis of Curl-Free DG Schemes – Second Order Example**

The von Neumann stability analysis of a DG scheme can be carried out in two different styles. The first is to convert the DG equations into a finite-difference-like form (Liu *et al*. [36], Zhang and Shu [48], Balsara and Käppeli [12]). The second approach is to identify the minimal number of modes, endow them with harmonic variation and then to directly carry out the stability analysis on the primal variables of the DG scheme (Balsara and Käppeli [14]). The latter approach works very well because it quickly allows us to identify the smallest number of variables that should be retained in a constraint-preserving DG scheme.

Here we describe the basic ingredients that go into carrying out a von Neumann stability analysis for a second order accurate curl-free DG-like scheme. Please focus on Fig. 3. In the right y-edge of the central zone we identify the modes $J_0^{y+}(t)$ and $J_y^{y+}(t)$ as the mean y-component of the vector field and its linear variation in the y-direction. Similarly, in the top x-edge of the central zone we identify the modes $J_0^{x+}(t)$ and $J_x^{x+}(t)$ as the mean x-component of the vector field and its linear variation in the x-direction. In the spirit of a DG scheme, the modes are endowed with time-dependence. In the spirit of a harmonic variation, we assume rectangular zones of size $\Delta x$ and $\Delta y$ so that the Fourier modes vary as $e^{-i(k_x x + k_y y)}$ where the wave vector is given by $(k_x, k_y)$. In fact, we simplify even further by assuming square zones in most parts of this paper. Because we use Fourier modes in a von Neumann stability analysis, the modes in the left y-edge are related to the modes in the right y-edge. Similarly, the modes in the bottom x-edge are related to the modes in the top x-edge. The relationship goes as follows

$$J_0^{y-}(t) = J_0^{y+}(t) e^{-ik_x \Delta x} \quad ; \quad J_y^{y-}(t) = J_y^{y+}(t) e^{-ik_x \Delta x} \quad ;$$
$$J_0^{x-}(t) = J_0^{x+}(t) e^{-ik_y \Delta y} \quad ; \quad J_x^{x-}(t) = J_x^{x+}(t) e^{-ik_y \Delta y} \tag{4.1}$$

Fig. 3 shows even further inter-relationships between the modes that reside in the edges once the Fourier modal variation is assumed. At first blush it would seem that each zone in Fig. 3 has four independent pieces of information given by $J_0^{y+}(t)$, $J_y^{y+}(t)$, $J_0^{x+}(t)$ and $J_x^{x+}(t)$. However, because of the curl-free constraint, there are only three independent pieces of information. This becomes apparent when we use eqn. (4.1) to write the discrete curl-free condition in zone $(i,j)$ of Fig. 3 as follows

$$\frac{J_0^{y+}(t) - J_0^{y+}(t) e^{-ik_x \Delta x}}{\Delta x} - \frac{J_0^{x+}(t) - J_0^{x+}(t) e^{-ik_y \Delta y}}{\Delta y} = 0 \quad \Leftrightarrow \quad J_0^{x+}(t) = J_0^{y+}(t) \frac{\Delta y}{\Delta x} \frac{1 - e^{-ik_x \Delta x}}{1 - e^{-ik_y \Delta y}} \tag{4.2}$$

As a result of eqn. (4.2), the only independent variables in the von Neumann stability analysis of a second order accurate, curl-free DG-like scheme are $J_0^{y+}(t)$, $J_y^{y+}(t)$ and $J_x^{x+}(t)$. This simplifies the analysis quite considerably.

Fig. 3 shows how the facial modes at all the faces of the mesh are inter-related because of the Fourier modes and their spatial variation. As a result, the curl-free reconstruction of the vector field in each zone of Fig. 3 can be symbolically carried out using a computer algebra system. Eqns. (2.8a), (2.8b) and (2.9b) (in their curl-free forms) can then again be symbolically expressed using the same computer algebra system. The result is that we obtain expressions for the time rate of change of $J_0^{y+}(t)$, $J_y^{y+}(t)$ and $J_x^{x+}(t)$ that can be written in terms of $J_0^{y+}(t)$, $J_y^{y+}(t)$ and $J_x^{x+}(t)$. In other words, we have reduced the problem of evaluating a single stage in the multistage RK-timestepping to the problem of obtaining a linear system of ODEs that look as follows:-

$$\frac{\partial}{\partial t}\begin{pmatrix} J_0^{y+}(t) \\ J_y^{y+}(t) \\ J_x^{x+}(t) \end{pmatrix} = \begin{pmatrix} A_{11} & A_{12} & A_{13} \\ A_{21} & A_{22} & A_{23} \\ A_{31} & A_{32} & A_{33} \end{pmatrix} \begin{pmatrix} J_0^{y+}(t) \\ J_y^{y+}(t) \\ J_x^{x+}(t) \end{pmatrix} \tag{4.3}$$

The nine coefficients in the matrix shown in eqn. (4.3) depend only on the wave numbers $k_x$ and $k_y$, the velocities $v_x$ and $v_y$, and the zone sizes $\Delta x$ and $\Delta y$. They are explicitly given in Appendix A. We can also formally define the vector of unknowns as $\mathbf{V}(t) = \left(J_0^{y+}(t),\ J_y^{y+}(t),\ J_x^{x+}(t)\right)^T$. As a result, eqn. (4.3) can be formally written as $\partial_t \mathbf{V}(t) = \mathbf{A}\,\mathbf{V}(t)$, where "**A**" is the 3×3 matrix shown in eqn. (4.3).

We then discretize eqn. (4.3) in time with an explicit *m*-stage Runge-Kutta scheme having a timestep $\Delta t$ of the form

$$\mathbf{V}^{(0)} = \mathbf{V}(t^n)$$
$$\mathbf{V}^{(i)} = \sum_{k=0}^{i-1}\left(\alpha_{i,k}\mathbf{I} + \Delta t\,\beta_{i,k}\mathbf{A}\right)\mathbf{V}^{(k)} \quad \text{for } i = 1,\ldots,m \tag{4.4}$$
$$\mathbf{V}(t^{n+1}) = \mathbf{V}^{(m)}$$

Here "**I**" is the identity matrix. The expressions for the coefficients $\alpha_{i,k}$ and $\beta_{i,k}$ can be found in Gottlieb *et al*. [33] and also Spiteri and Ruuth [43], [44]. Given the linearity of our DG scheme, we can write the time update as

$$\mathbf{V}(t^{n+1}) = \mathbf{G}\,\mathbf{V}(t^n) \tag{4.5}$$

Here "**G**" is known as the amplification matrix of the scheme. It depends on the coefficients of the Runge-Kutta scheme, on the timestep $\Delta t$ and the matrix "**A**" from eqn. (4.3). For the second order SSP-RK scheme we can write the amplification matrix as

$$\mathbf{G} = \mathbf{I} + \Delta t \mathbf{A} + \frac{\Delta t^2}{2}\mathbf{A}^2 \tag{4.6}$$

Likewise, for the third order SSP-RK scheme we can write the amplification matrix as

$$\mathbf{G} = \mathbf{I} + \Delta t \mathbf{A} + \frac{\Delta t^2}{2}\mathbf{A}^2 + \frac{\Delta t^3}{3!}\mathbf{A}^3 \tag{4.7}$$

In the next Section we will use this amplification matrix to devise our von Neumann stability analysis. This completes our description of the mathematics associated with the von Neumann stability analysis at second order. Higher orders can be done similarly.

**V) Results from the von Neumann Stability Analysis of Globally Curl-free DG-like Schemes**

By taking a close look at eqn. (2.1) we realize that the von Neumann stability analysis will depend on the angle that the velocity vector, $(v_x, v_y)$ makes with respect to the x-axis of the mesh. Furthermore, eqns. (4.1) and (4.2) show us that the von Neumann stability analysis will also depend on the angle that the wave vector $(k_x, k_y)$ makes with respect to the velocity vector $(v_x, v_y)$. For this reason, the stability analysis depends on multiple parameters. Besides, owing to the fact that the curl only manifests itself in two or more dimensions, it has to be multidimensional. For all of these reasons, we have only been able to carry out a von Neumann stability analysis for curl-free WENO-like, PNPM-like and DG-like schemes up to fourth order of accuracy. However, we realize that such a stability analysis that is done in two dimensions and for a full scheme can give us a wealth of information; and that information is catalogued in the ensuing two Sub-sections.

The first insight that we would like to extract from such a stability analysis is the maximal CFL number for which the scheme is stable. For eqn. (2.1) the Fourier modes are indeed propagating with the velocity vector; therefore, the velocity vector sets the signal speed. For each choice of spatial accuracy, we can choose a temporal accuracy for our SSP-RK scheme that is comparable or greater than the spatial accuracy. The upshot is that for each choice of spatial and temporal accuracy, we can identify a maximal CFL number. Please realize that this involves sweeping through all velocities in two-dimensions and for each choice of velocity we have to sweep over all the wavenumbers that are permitted on the mesh. Stable CFL numbers are identified as the ones for which all possible wave vectors return an amplification matrix all of whose eigenvalues have an absolute value that is less than or equal to unity. Such a study of the maximal CFL number is documented in Sub-section V.1.

DG-like schemes can be very accurate, even when they are compared to their WENO-like counterparts. But we need to visually appreciate that. For that reason, we choose velocity vectors that make angles of 0º, 15º, 30º and 45º to the mesh. For each of those velocity vectors, we sweep through all possible angles that the wave vector can make with respect to the velocity. This allows us to visualize the dissipation and dispersion errors of the schemes that we analyze. This information is shown in Sub-section V.2.

### V.a) Maximal CFL Numbers from Stability Analysis

We identify the CFL number in each of the two directions by $C_x = v_x \Delta t / \Delta x$ and $C_y = v_y \Delta t / \Delta y$. For each choice of the CFL number in either of the two directions, we sweep over all possible wave numbers $(k_x \Delta x, k_y \Delta y) \in [-\pi/2, \pi/2] \times [-\pi/2, \pi/2]$. Fig. 4 shows a colorized plot of the eigenvector of the amplification matrix with the largest absolute value for second, third and fourth order curl-free DG-like schemes. The white polygons in Fig 4 identify the domain of stability for which the absolute value described above is less than or equal to unity. The white circles in Fig. 4 are the largest circles that can be inscribed in the polygons. The radii of those circles give us the largest effective CFL that we should use for each of those DG-like schemes. Fig. 4 is intended to give us a flavor of the process that goes into finding the largest effective CFL number that we can find from our von Neumann stability analysis.

Fig. 4 corresponds to a situation where the order of temporal accuracy of our SSP-RK time stepping schemes indeed matched the spatial accuracy of the DG-like discretization. But we can use several possible SSP-RK schemes with each DG-like discretization, as long as the temporal accuracy is at least as large as the spatial accuracy. Table I shows the largest effective CFL number for a range of curl-free DG-like spatial discretizations and a range of temporal accuracies. In all cases, SSP-RK schemes were used for the temporal update. We see that for each scheme, an increasing temporal accuracy results in a larger effective CFL, a result what conforms with the findings of Zhang and Shu [48] and Liu *et al*. [36].

**Table I shows the largest effective CFL number for a range of curl-free DG-like spatial discretizations and a range of temporal accuracies. In all cases, SSP-RK schemes were used for the temporal update.**

|            | P=0    | P=1    | P=2    | P=3    |
|------------|--------|--------|--------|--------|
| RK1        | 0.7071 | _____  | _____  | _____  |
| SSP-RK2    | 0.7071 | 0.3162 | _____  | _____  |
| SSP-RK3    | 0.8884 | 0.3906 | 0.2069 | _____  |
| SSP-RK(5,4)| 1.5495 | 0.6367 | 0.3401 | 0.2143 |

Because we have erected the machinery of the von Neumann stability analysis, we can also use it to analyze the largest effective CFL number when other spatial discretizations are used. For example, when we retain only the time-evolution of the zeroth mode in eqns. (2.8a) and (2.9a), we get a family of curl-free WENO-like schemes. For those schemes, all the higher modes, up to the desired order of accuracy, have to be reconstructed. Table II shows the largest effective CFL number for a range of curl-free WENO-like spatial discretizations and a range of temporal accuracies. In all cases, SSP-RK schemes were used for the temporal update. As expected, we see that the curl-free WENO-like schemes have CFL numbers that are much larger than the curl-free DG-like schemes.

**Table II shows the largest effective CFL number for a range of curl-free WENO-like spatial discretizations and a range of temporal accuracies. In all cases, SSP-RK schemes were used for the temporal update.**

|     | P0P0   | P0P1  | P0P2  | P0P3  |
|-----|--------|-------|-------|-------|
| RK1 | 0.7071 | _____ | _____ | _____ |
| SSP-RK2 | 0.7071 | 0.7071 | ______ | ______ |
| SSP-RK3 | 0.8884 | 0.8318 | 1.1507 | ______ |
| SSP-RK(5,4) | 1.5495 | 1.2252 | 1.4859 | 1.3040 |

Table III shows the largest effective CFL number for a range of curl-free P1PM-like spatial discretizations and a range of temporal accuracies. In addition to retaining the time evolution of the modes in eqns. (2.8a) and (2.9a), such schemes also evolve the first moments from eqns. (2.8b) and (2.9b). As before, SSP-RK schemes were used for the temporal update. Comparing the CFL numbers from Table III to those from Tables I and II, we see that curl-free P1PM-like schemes give us CFL numbers that are somewhat smaller than those of their WENO counterparts but substantially larger than their DG counterparts. We will further see in the next Sub-section that P1PM-like schemes retain their first moments and that gives them dissipation and dispersion properties that are closer to their DG counterparts. The physical reason for that is because the linear mode retains most of the variation in the zone; consequently, much of the accuracy is retained. We, therefore, understand why curl-free P1PM-like schemes retain an important utilitarian position in the full range of schemes studied here.

**Table III shows the largest effective CFL number for a range of curl-free P1PM-like spatial discretizations and a range of temporal accuracies. In all cases, SSP-RK schemes were used for the temporal update.**

| | P1P2 | P1P3 |
|---|---|---|
| SSP-RK3 | 0.3903 | ______ |
| SSP-RK(5,4) | 0.6260 | 0.6799 |

This completes our study of the CFL number of curl-free DG-like schemes and their cousins.

**V.b) Dissipation and Dispersion Properties of DG and PNPM Schemes**

We now wish to study the dissipation and dispersion properties of the curl-free WENO-like, PNPM-like and DG-like schemes. We will study these properties for second, third and fourth order, so that we have a clear understanding of the improving wave propagation properties of these schemes with increasing order. By the same token, we will also be able to inter-compare between the WENO-like, PNPM-like and DG-like schemes. We expect that retaining more moments, and evolving them consistent with the governing PDE, should give us schemes with improved wave propagation characteristics. For all the data shown in this Sub-section, the temporal accuracy was made to match the spatial accuracy. While the von Neumann stability analysis for curl-free WENO-like schemes was already documented in Balsara *et al*. [16], we present it again here so that one can inter-compare with the PNPM-like and DG-like schemes. The von Neumann stability analysis of the curl-free PNPM-like and DG-like schemes is being presented for the very first time here.

In each instance we choose velocity vectors that make angles of 0º, 15º, 30º and 45º to the mesh and use a CFL number that is 0.9 times the maximum shown in either Table I or II or III. We then let the wave number $(k_x, k_y)$ sweep through all angles, from $-180º$ to $+180º$ relative to the velocity vector $(v_x, v_y)$. For each of those angles, we plot out $1-|\text{amplification factor}|$ for the scheme. If this number is non-negative and close to zero, it indicates that the scheme has low dissipation. The phase of the amplification factor gives us a measure of the propagation speed of the waves. For each angle between the velocity vector and the wave number, we also plot out the error in the phase speed. As the wavelength increases relative to the zone size, we expect the schemes to propagate waves with increasing accuracy. As a result, we consider wavelengths that are $5\,\Delta x$, $10\,\Delta x$ and $15\,\Delta x$. In the next eight figures that follow, wavelengths of $5\,\Delta x$ are always shown with a blue color, wavelengths of $10\,\Delta x$ are always shown with a green color and wavelengths of $15\,\Delta x$ are always shown with a red color.

Figs. 5, 6 and 7 show the wave propagation characteristics of curl-free WENO-like schemes at second, third and fourth order respectively. We see that as we go from second to fourth order, the dissipation (as measured by $1-|\text{amplification factor}|$) improves by an order of magnitude for each of the three wavelengths considered here! Similarly, as we go from second to fourth order, the phase error is also reduced by an order of magnitude. Table IV shows the minimum of the absolute value of the amplification factor for all possible velocity directions and all angles between the velocity and wave number vectors for curl-free WENO-like schemes when we have waves with wavelength 5 $\Delta x$, 10 $\Delta x$ and 15 $\Delta x$. In that same table, we also show the maximum phase error for the similar situation and for the same wavelengths. In other words, Table IV was extracted from Figs. 5, 6 and 7 and allows us to quantify the most significant aspect of those figures. Table IV, therefore, allows us to make an important practical decision. Say we want to carry out a simulation with WENO-like schemes we want to meet a target set of dissipation and dispersion properties, Table IV shows us what our options are. We may indeed choose a lower order scheme and use a lot of zones to cover the characteristic wavelength in the simulation. But we see that we can also choose a higher order scheme and use fewer zones to cover the characteristic wavelength in the simulation.

**Table IV shows the minimum of the absolute value of the amplification factor for all possible velocity directions and all angles between the velocity and wave number vectors for curl-free WENO-like schemes when we have waves with wavelength 5 $\Delta x$, 10 $\Delta x$ and 15 $\Delta x$. We also show the maximum phase error for the same wavelengths.**

| Min of |Amplification factor| | $\lambda = 5\ \Delta x$ | $\lambda = 10\ \Delta x$ | $\lambda = 15\ \Delta x$ |
|---|---|---|---|
| WENO-O2 | 0.8672298 | 0.9908930 | 0.9981729 |
| WENO-O3 | 0.7455074 | 0.9787628 | 0.9955671 |
| WENO-O4 | 0.9105516 | 0.9980383 | 0.9998192 |
| Max of Phase Error | $\lambda = 5\ \Delta x$ | $\lambda = 10\ \Delta x$ | $\lambda = 15\ \Delta x$ |
| WENO-O2 | 1.6211953E-01 | 5.5976172E-02 | 2.6459753E-02 |
| WENO-O3 | 6.8417271E-02 | 5.4542411E-03 | 1.1453074E-03 |
| WENO-O4 | 2.5814369E-02 | 1.0046737E-03 | 2.0775987E-04 |

There is no second order P1P1 scheme, because such a scheme would be identical to a second order DG scheme. However, our study of CFL numbers has shown us that curl-free third order P1P2-like and fourth order P1P3-like schemes still retain a very robust CFL number. We, therefore, want to know whether such schemes have superior wave propagation characteristics relative to the WENO-like schemes that we studied in the previous paragraph. Figs. 8 and 9 show the wave propagation characteristics of curl-free P1P2-like and P1P3-like schemes at third and fourth order respectively. We should, therefore, compare Fig. 8 to Fig. 6 because they both pertain to third order schemes. Similarly, we should compare Fig. 9 to Fig. 7 because they both pertain to fourth order schemes. The results are quite interesting. We see that Fig. 8 and Fig. 6 show comparable quality of wave propagation indicating that at third order the advantages are minimal. It is possible that this lack of significant improvement has to do with the fact that SSP-RK3 timestepping has excessive stabilization. Now let us turn to comparing Fig. 9 and Fig. 7. At fourth order we do see that the P1P3-like scheme outperforms the WENO-O4 scheme by almost an order of magnitude. It shows the value of designing PNPM schemes as half way houses between WENO and DG schemes. Table V shows the minimum of the absolute value of the amplification factor for all possible velocity directions and all angles between the velocity and wave number vectors for curl-free PNPM-like schemes when we have waves with wavelength 5 $\Delta x$, 10 $\Delta x$ and 15 $\Delta x$. In that same table, we also show the maximum phase error for the similar situation and for the same wavelengths. In other words, Table V was extracted from Figs. 8 and 9 and allows us to quantify the most significant aspect of those figures. Again, Table V can help with practical decision-making. It shows us, for example, that fourth order P1P3-like schemes do give us a substantial improvement over fourth order WENO-like scheme while incurring only a modest increase in computational complexity.

**Table V shows the minimum of the absolute value of the amplification factor for all possible velocity directions and all angles between the velocity and wave number vectors for curl-free PNPM-like schemes when we have waves with wavelength 5 $\Delta x$, 10 $\Delta x$ and 15 $\Delta x$. We also show the maximum phase error for the same wavelengths.**

| Min of |Amplification factor| | $\lambda = 5\ \Delta x$ | $\lambda = 10\ \Delta x$ | $\lambda = 15\ \Delta x$ |
|---|---|---|---|
| P1P2 | 0.9869830 | 0.9990722 | 0.9998118 |
| P1P3 | 0.9943549 | 0.9998913 | 0.9999898 |
| Max of Phase Error | $\lambda = 5\ \Delta x$ | $\lambda = 10\ \Delta x$ | $\lambda = 15\ \Delta x$ |
| P1P2 | 5.2001351E-03 | 3.1972379E-04 | 6.4931856E-05 |
| P1P3 | 1.1220642E-03 | 1.4773508E-04 | 3.2952515E-05 |

While they have the smallest CFL numbers, DG-like schemes hold out the promise of almost spectral-like accuracy with increasing order of accuracy; for finite volume-based approaches this is now viewed as accepted fact. We are now in a position to test that contention as it pertains to curl-free DG-like schemes. Figs. 10, 11 and 12 show the wave propagation characteristics of curl-free DG-like schemes at second, third and fourth order respectively. With increasing order, the curl-free DG-like schemes do show significant improvement when inter-compared amongst themselves. Let us, therefore, compare across algorithms since we have all the data concatenated in one place. Fig. 10 should be compared to Fig. 5. Fig. 11 should be compared to Figs. 6 and 8. Fig. 12 should be compared to Figs. 7 and 9. We see that the wave propagation characteristics of the second order DG-like scheme are entirely competitive with the wave propagation characteristics of the fourth order WENO-like scheme! The fourth order DG-like scheme is also somewhat superior to the fourth order P1P3-like scheme; but please recall that this comes with a substantial increase in programming complexity and a decrease in the CFL number. Table VI shows the minimum of the absolute value of the amplification factor for all possible velocity directions and all angles between the velocity and wave number vectors for curl-free DG-like schemes when we have waves with wavelength 5 $\Delta x$, 10 $\Delta x$ and 15 $\Delta x$. In that same table, we also show the maximum phase error for the similar situation and for the same wavelengths. In other words, Table VI was extracted from Figs. 10, 11 and 12 and allows us to quantify the most significant aspect of those figures. We see that the dissipation and dispersion characteristics of the curl-free DG-like schemes that we have designed are indeed excellent. Comparing Tables V and VI we also see that the curl-free PNPM-like schemes are not far behind. Therefore, we see that the DG-like schemes are the go-to scheme when superlative performance is the only driving consideration. However, if one wants lower computational complexity and more robust timesteps, the PNPM-like schemes also present themselves as very attractive choices.

**Table VI shows the minimum of the absolute value of the amplification factor for all possible velocity directions and all angles between the velocity and wave number vectors for curl-free DG-like schemes when we have waves with wavelength 5 $\Delta x$, 10 $\Delta x$ and 15 $\Delta x$. We also show the maximum phase error for the same wavelengths.**

| Min of \|Amplification factor\| | $\lambda = 5\ \Delta x$ | $\lambda = 10\ \Delta x$ | $\lambda = 15\ \Delta x$ |
|---|---|---|---|
| DG P=1 | 0.9889383 | 0.9991534 | 0.9998251 |
| DG P=2 | 0.9937189 | 0.9995565 | 0.9999105 |
| DG P=3 | 0.9994633 | 0.9999897 | 0.9999991 |
| Max of Phase Error | $\lambda = 5\,\Delta x$ | $\lambda = 10\,\Delta x$ | $\lambda = 15\,\Delta x$ |
| DG P=1 | 3.0344813E-02 | 6.4200877E-03 | 2.7378616E-03 |
| DG P=2 | 7.6077271E-03 | 5.1942472E-04 | 1.0415238E-04 |
| DG P=3 | 3.2546521E-03 | 2.5127499E-04 | 5.1804468E-05 |

While this Section has been focused on curl-free methods, we point out that curl-preserving methods only require a few additional terms in eqns. (2.8) and (2.9) compared to curl-free schemes. Without an underlying fluid-dynamical PDE system that supplies additional terms like density, velocity and temperature, it is not possible to obtain the source terms and other types of terms which would make eqn. (1.2) curl-preserving instead of curl-free. Furthermore, computer algebra systems are just not adept enough to support a more extensive stability analysis for larger PDE systems. For these reasons, the analysis presented here is focused on curl-free methods; but the insights developed here will extend to all curl-preserving methods.

## VI) Numerical Results

In this Section, we present numerical experiments confirming that the developed curl-free DG schemes reach their expected design accuracies. Results for all the DG-like (PNPN) and PNPM (P0PN for WENO-like and P1PN for HWENO-like) up to fourth order are reported for two smooth test problems. All the tests are run with 95% of the maximal CFL number (Tables I-III) of the respective scheme.

### VI.a) Plane wave test problem

The first test problem is the propagation of a plane wave in a Cartesian domain of size $[-1/2,+1/2]^2$ with periodic boundaries. The plane wave is advected with velocity $v_x = v_y = 1$ diagonally through the domain. The curl-free field **J** is initialized from a potential

$$\phi(x,y) = \cos(k_x x + k_y y),$$

where we set $k_x = k_y = 2\pi$. The x- and y-field components of the curl-free field are then given by

$$J_x = \frac{\partial \phi}{\partial x} \text{ and } J_y = \frac{\partial \phi}{\partial y}.$$

The setup is run up to time $t_f = 1$, by which time the plane wave has propagated once through the domain, and the accuracy of the schemes is evaluated. Moreover, we also show the preservation of the quadratic field energy $(J_x^2 + J_y^2)/2$ highlighting the dissipation characteristics of the schemes.

Table VII shows the $L_1$ and $L_\infty$ errors at the final time for the PNPN (N=1,2,3) DG-like schemes for resolutions from 8×8 up to 64×64. Table VIII shows the convergence analysis of the P0PN (N=1,2,3) WENO-like schemes. Table IX shows the convergence analysis of the P1PN (N=1,2,3) HWENO-like schemes. The tables also catalogue the final quadratic energy as a fraction of the initial quadratic energy. We observe that all schemes reach their design accuracy. We also take note of the improved quadratic field energy preservation with increasing order of accuracy. Note that we did nothing special in the scheme to ensure that quadratic field energy is conserved; as a result, the rather good preservation of quadratic energy is entirely a consequence of the accuracy of the method. This is especially true for the DG-like schemes which preserve quadratic energy very well especially as the resolution is increased. The WENO-like schemes show slightly inferior energy preservation characteristics. However, the latter allow much larger time steps due to their larger allowed CFL numbers. The HWENO-like schemes show nearly the same quadratic energy preservation properties as the DG-like schemes and, furthermore, allow larger time steps similar to the WENO-like schemes. Consequently we see that the HWENO-like schemes may be viewed as an efficient compromise between the extreme accuracy of the DG-like schemes and the much larger time steps of the WENO-like schemes.

**Table VII Accuracy analysis (Plane Wave Test) of the PNPN (N=1,2,3) DG-like schemes. The total quadratic energy on the mesh as a fraction of its initial value is also shown.**

| P1P1 | L$_1$ Error | L$_1$ Accuracy | L$_\infty$ Error | L$_\infty$ Accuracy | Total Quadratic Energy |
|---|---|---|---|---|---|
| 8×8 | 1.054E+00 | - | 1.710E+00 | - | 0.767072144659713 |
| 16×16 | 1.959E-01 | 2.43 | 3.041E-01 | 2.49 | 0.963508927621496 |
| 32×32 | 3.642E-02 | 2.43 | 5.699E-02 | 2.42 | 0.995169598723023 |
| 64×64 | 7.897E-03 | 2.21 | 1.240E-02 | 2.20 | 0.999386133084480 |
| P2P2 | L$_1$ Error | L$_1$ Accuracy | L$_\infty$ Error | L$_\infty$ Accuracy | Total Quadratic Energy |
| 8×8 | 8.529E-01 | - | 1.335E+00 | - | 0.798900332986684 |
| 16×16 | 1.229E-01 | 2.79 | 1.931E-01 | 2.79 | 0.969506454484418 |
| 32×32 | 1.584E-02 | 2.96 | 2.488E-02 | 2.96 | 0.996044633506995 |
| 64×64 | 1.993E-03 | 2.99 | 3.130E-03 | 2.99 | 0.999501861659339 |
| P3P3 | L$_1$ Error | L$_1$ Accuracy | L$_\infty$ Error | L$_\infty$ Accuracy | Total Quadratic Energy |
| 8×8 | 1.150E-01 | - | 1.711E-01 | - | 0.982433477747556 |
| 16×16 | 8.007E-03 | 3.84 | 1.235E-02 | 3.79 | 0.999160570223597 |
| 32×32 | 5.131E-04 | 3.96 | 8.082E-04 | 3.93 | 0.999961937114448 |
| 64×64 | 3.256E-05 | 3.98 | 5.115E-05 | 3.98 | 0.999998147852102 |

**Table VIII Accuracy analysis (Plane Wave Test) of the P0PN (N=1,2,3) WENO-like schemes. The total quadratic energy on the mesh as a fraction of its initial value is also shown.**

| P0P1 | L$_1$ Error | L$_1$ Accuracy | L$_\infty$ Error | L$_\infty$ Accuracy | Total Quadratic Energy |
|---|---|---|---|---|---|
| 8×8 | 4.993E+00 | - | 6.970E+00 | - | 0.147331805631007 |
| 16×16 | 1.687E+00 | 1.57 | 3.081E+00 | 1.18 | 0.672786308056742 |
| 32×32 | 7.354E-01 | 1.20 | 1.342E+00 | 1.20 | 0.961889800888593 |
| 64×64 | 1.939E-01 | 1.92 | 5.060E-01 | 1.41 | 0.996184224345619 |
| P0P2 | L$_1$ Error | L$_1$ Accuracy | L$_\infty$ Error | L$_\infty$ Accuracy | Total Quadratic Energy |
| 8×8 | 2.377E+00 | - | 3.458E+00 | - | 0.493454736716243 |
| 16×16 | 3.817E-01 | 2.64 | 5.868E-01 | 2.56 | 0.906990382419879 |
| 32×32 | 5.000E-02 | 2.93 | 7.805E-02 | 2.91 | 0.987543393334568 |
| 64×64 | 6.291E-03 | 2.99 | 9.866E-03 | 2.98 | 0.998428021784668 |
| P0P3 | $L_1$ Error | $L_1$ Accuracy | $L_\infty$ Error | $L_\infty$ Accuracy | Total Quadratic Energy |
| 8×8 | 5.523E-01 | - | 9.562E-01 | - | 0.863259629337563 |
| 16×16 | 1.244E-02 | 5.47 | 3.046E-02 | 4.97 | 0.996553938792429 |
| 32×32 | 3.951E-04 | 4.98 | 8.497E-04 | 5.16 | 0.999902248490279 |
| 64×64 | 1.387E-05 | 4.83 | 2.503E-05 | 5.09 | 0.999997008920715 |

**Table IX Accuracy analysis (Plane Wave Test) of the P1PN (N=2,3) HWENO-like schemes. The total quadratic energy on the mesh as a fraction of its initial value is also shown.**

| P1P2 | $L_1$ Error | $L_1$ Accuracy | $L_\infty$ Error | $L_\infty$ Accuracy | Total Quadratic Energy |
|---|---|---|---|---|---|
| 8×8 | 8.478E-01 | - | 1.254E+00 | - | 0.800179565838325 |
| 16×16 | 1.244E-01 | 2.77 | 1.918E-01 | 2.71 | 0.969170286097292 |
| 32×32 | 1.628E-02 | 2.93 | 2.546E-02 | 2.91 | 0.995933787552492 |
| 64×64 | 2.065E-03 | 2.98 | 3.239E-03 | 2.97 | 0.999483905512296 |
| P1P3 | $L_1$ Error | $L_1$ Accuracy | $L_\infty$ Error | $L_\infty$ Accuracy | Total Quadratic Energy |
| 8×8 | 2.331E-01 | - | 3.870E-01 | - | 0.942799109295127 |
| 16×16 | 1.341E-02 | 4.12 | 2.090E-02 | 4.21 | 0.997772917919100 |
| 32×32 | 7.612E-04 | 4.14 | 1.183E-03 | 4.14 | 0.999927383586396 |
| 64×64 | 4.600E-05 | 4.05 | 7.202E-05 | 4.04 | 0.999997714221129 |

**VI.b) Vortex test problem**

The second test problem consists of the advection of a localized curl-free vortex similar to the magnetic vortex for the induction equation. The Cartesian domain extents are $[-10,+10]^2$ with periodic boundary condition. The vortex is initialized from a potential

$$\phi(x,y) = e^{\frac{1}{2}(1-r^2)},$$

where $r = \sqrt{x^2+y^2}$. This results in a field given by

$$\mathbf{J} = \nabla\phi = -e^{\frac{1}{2}(1-r^2)}[x,y]^T.$$

The advection velocity is set to $v_x = v_y = 1$. The problem is simulated for a time $t_f = 20$, by which point the vortex was advected once through the square domain in the diagonal direction till it returns to its initial position. At final time, we measure the accuracy in the $L_1$ and $L_\infty$ errors norms. Moreover, we also show the preservation of the quadratic field energy $(J_x^2 + J_y^2)/2$ highlighting the dissipation characteristics of the schemes.

Table X shows the $L_1$ and $L_\infty$ errors at the final time for the PNPN (N=1,2,3) RKDG-like schemes for resolutions from 16×16 up to 256×256. Table XI shows the convergence analysis of the P0PN (N=1,2,3) WENO-like schemes. Table XII shows the convergence analysis of the P1PN (N=1,2,3) HWENO-like schemes. The tables also catalogue the final quadratic energy as a fraction of the initial quadratic energy. We observe that all schemes reach their design accuracy on the chosen mesh resolutions, even for this highly spatially localized vortex. Note that much of the field variation is confined around a circle with unit radius, which corresponds to one-tenth of the computational domain. We find that the presented schemes concurrently have quadratic energy preservation with, nevertheless, excellent accuracy.

As before, we did nothing special in the scheme to ensure that quadratic field energy is conserved; consequently, the rather good preservation of quadratic energy is entirely a consequence of the accuracy of the method. For the quadratic field energy, we observe similar trends to the plane wave test problem. To further highlight the point, Fig. 13 shows the quadratic field energy preservation characteristics from Table X-XII graphically. Each panel shows all the available schemes up to fourth order of accuracy. Fig. 13a underlines that the curl-free DG-like schemes show excellent energy preserving properties with increasing order of accuracy. In Fig. 13b, we also observe again that the energy preserving properties of the curl-free WENO-like schemes are somewhat inferior in the pre-asymptotic regime. But the improving trend with increasing order of accuracy is also clearly visible. In Fig. 13c, we see that the curl-free P1PN-like schemes share almost the same preservation properties as the curl-free DG-like schemes. The latter fact, and their substantially larger allowed CFL numbers (and hence time steps), highlight again that the P1PN-like schemes are very efficient curl constraint-preserving methods that share very desirable qualities from both full DG-like and WENO-like schemes.

Lastly, we also catalogue that the schemes designed here are curl-preserving and can indeed reach the curl-free limit even when they are integrated for long simulation times. Dumbser *et al*. [30] have shown that if a classical higher order Godunov scheme is applied to eqn. (1.2), the numerical instability manifests itself as an explosive increase in the discrete circulation when the simulation is run over long periods of time; see Fig. 5 of that paper. Therefore, we would like to demonstrate that the discrete circulation is held down to machine precision when the simulation is run for a long period of time when we use the methods designed here. We would indeed like to go one step further and plot out the time series of the maximum pointwise error in the curl of **J**. In other words, we know that the vector field is a polynomial, see eqns. (3.3) or (3.6) for instance, so that we can evaluate the pointwise curl at any point within a zone (because the polynomials are differentiable). We then choose 2×2, 3×3 or 4×4 uniformly spaced points that are internal to each zone at second, third and fourth orders respectively. We then evaluate the maximum of the absolute value of the curl at each and every of those internal points for all of the zones on the mesh and we plot this maximum value as a function of time. Let us ask why this demonstration matters? We see from eqns. (2.8) and (2.9) that in a curl-preserving scheme we will also need terms like $\langle v_x (\nabla \times \mathbf{J})_z \rangle$ and $\langle v_y (\nabla \times \mathbf{J})_z \rangle$, and other terms like it, at the edges of the mesh. These have to evaluated from either side of the edge that is being considered. Therefore, when we approach the curl-free limit, a curl-preserving scheme should naturally obtain curl-free evolution. This demonstration that the maximum pointwise error in the curl of **J** remains close to machine zero over long simulation times guarantees that such a limit is met.

In order to show that the curl remains close to machine zero at all points on the mesh even during long time integration, we have run the vortex problem on a 64×64 zone mesh to a final time of 200. This time corresponds to the vortex making ten passages through the periodic computational domain. Fig. 14 shows the maximum pointwise error of the curl of **J** as a function of time for a 64×64 zone run of the vortex problem. Fig. 14a shows the evolution of the maximum pointwise curl as a function of time for the 2$^{nd}$, 3$^{rd}$ and 4$^{th}$ order curl-free DG-like schemes. Fig. 14b shows the evolution of the maximum pointwise curl as a function of time for the 3$^{rd}$ and 4$^{th}$ order curl-free P1PN-like schemes. Fig. 14c shows the evolution of the maximum pointwise curl as a function of time for the 2$^{nd}$, 3$^{rd}$ and 4$^{th}$ order curl-free WENO-like schemes. The figure shows that all our curl-preserving schemes can preserve the curl constraint up to machine accuracy in simulations that are run for long integration times.

**Table X Accuracy analysis (Vortex Test) of the PNPN (N=1,2,3) DG-like schemes. The total quadratic energy on the mesh as a fraction of its initial value is also shown.**

| P1P1 | $L_1$ Error | $L_1$ Accuracy | $L_\infty$ Error | $L_\infty$ Accuracy | Total Quadratic Energy |
|---|---|---|---|---|---|
| 16×16 | 3.960E-02 | - | 1.296E+00 | - | 0.244138062854683 |
| 32×32 | 1.937E-02 | 1.03 | 9.775E-01 | 0.41 | 0.584068951760809 |
| 64×64 | 4.780E-03 | 2.02 | 3.237E-01 | 1.59 | 0.887813286147527 |
| 128×128 | 8.569E-04 | 2.48 | 6.715E-02 | 2.27 | 0.982478078399363 |
| 256×256 | 1.678E-04 | 2.35 | 1.243E-02 | 2.43 | 0.997712794953464 |
| P2P2 | $L_1$ Error | $L_1$ Accuracy | $L_\infty$ Error | $L_\infty$ Accuracy | Total Quadratic Energy |
| 16×16 | 3.813E-02 | - | 1.017E+00 | - | 0.449340807458768 |
| 32×32 | 1.535E-02 | 1.31 | 6.368E-01 | 0.67 | 0.743031482000765 |
| 64×64 | 3.249E-03 | 2.24 | 1.748E-01 | 1.87 | 0.938198337740548 |
| 128×128 | 4.658E-04 | 2.80 | 2.755E-02 | 2.67 | 0.990782650850455 |
| 256×256 | 5.981E-05 | 2.96 | 3.601E-03 | 2.94 | 0.998809736688618 |
| P3P3 | $L_1$ Error | $L_1$ Accuracy | $L_\infty$ Error | $L_\infty$ Accuracy | Total Quadratic Energy |
| 16×16 | 1.672E-02 | - | 5.812E-01 | - | 0.837642154450034 |
| 32×32 | 2.898E-03 | 2.53 | 1.281E-01 | 2.18 | 0.980766190163135 |
| 64×64 | 2.427E-04 | 3.58 | 1.200E-02 | 3.42 | 0.999051579122300 |
| 128×128 | 1.607E-05 | 3.92 | 8.047E-04 | 3.90 | 0.999964537681624 |
| 256×256 | 1.019E-06 | 3.98 | 5.100E-05 | 3.98 | 0.999998674687398 |

**Table XI Accuracy analysis (Vortex Test) of the P0PN (N=1,2,3) WENO-like schemes. The total quadratic energy on the mesh as a fraction of its initial value is also shown.**

| P0P1 | $L_1$ Error | $L_1$ Accuracy | $L_\infty$ Error | $L_\infty$ Accuracy | Total Quadratic Energy |
|---|---|---|---|---|---|
| 16×16 | 4.396E-02 | - | 1.408E+00 | - | 0.016052344223764 |
| 32×32 | 3.955E-02 | 0.15 | 1.771E+00 | -0.33 | 0.060627418657843 |
| 64×64 | 2.399E-02 | 0.72 | 1.357E+00 | 0.38 | 0.274789236132377 |
| 128×128 | 7.655E-03 | 1.65 | 5.670E-01 | 1.26 | 0.773714109741705 |
| 256×256 | 1.988E-03 | 1.94 | 1.543E-01 | 1.88 | 0.979590328058657 |
| P0P2 | $L_1$ Error | $L_1$ Accuracy | $L_\infty$ Error | $L_\infty$ Accuracy | Total Quadratic Energy |
| 16×16 | 4.036E-02 | - | 1.418E+00 | - | 0.014812755301855 |
| 32×32 | 3.251E-02 | 0.31 | 1.671E+00 | -0.24 | 0.113246493972690 |
| 64×64 | 1.018E-02 | 1.67 | 6.669E-01 | 1.32 | 0.692818242056277 |
| 128×128 | 1.995E-03 | 2.35 | 1.484E-01 | 2.17 | 0.942819902537954 |
| 256×256 | 2.689E-04 | 2.89 | 2.108E-02 | 2.82 | 0.992210030651288 |
| P0P3 | $L_1$ Error | $L_1$ Accuracy | $L_\infty$ Error | $L_\infty$ Accuracy | Total Quadratic Energy |
| 16×16 | 3.903E-02 | - | 1.409E+00 | - | 0.042648512931034 |
| 32×32 | 1.931E-02 | 1.02 | 1.087E+00 | 0.37 | 0.483958370120984 |
| 64×64 | 1.940E-03 | 3.31 | 1.397E-01 | 2.96 | 0.953798557685217 |
| 128×128 | 8.435E-05 | 4.52 | 5.971E-03 | 4.55 | 0.998754807116245 |
| 256×256 | 4.283E-06 | 4.30 | 2.492E-04 | 4.58 | 0.999964258479139 |

**Table XII Accuracy analysis (Vortex Test) of the P1PN (N=2,3) HWENO-like schemes. The total quadratic energy on the mesh as a fraction of its initial value is also shown.**

| P1P2 | $L_1$ Error | $L_1$ Accuracy | $L_\infty$ Error | $L_\infty$ Accuracy | Total Quadratic Energy |
|---|---|---|---|---|---|
| 16×16 | 3.629E-02 | - | 1.172E+00 | - | 0.288607282457506 |
| 32×32 | 1.563E-02 | 1.21 | 7.317E-01 | 0.68 | 0.675080995909535 |
| 64×64 | 3.281E-03 | 2.25 | 1.797E-01 | 2.03 | 0.932780213745329 |
| 128×128 | 4.730E-04 | 2.79 | 2.774E-02 | 2.70 | 0.990441076783871 |
| 256×256 | 6.101E-05 | 2.95 | 3.655E-03 | 2.92 | 0.998771733291638 |
| P1P3 | $L_1$ Error | $L_1$ Accuracy | $L_\infty$ Error | $L_\infty$ Accuracy | Total Quadratic Energy |
| 16×16 | 2.919E-02 | - | 1.055E+00 | - | 0.346983109358981 |
| 32×32 | 7.340E-03 | 1.99 | 3.727E-01 | 1.50 | 0.872263471054451 |
| 64×64 | 5.285E-04 | 3.80 | 3.912E-02 | 3.25 | 0.992401678126276 |
| 128×128 | 2.700E-05 | 4.29 | 1.909E-03 | 4.36 | 0.999731098030406 |
| 256×256 | 1.514E-06 | 4.16 | 9.190E-05 | 4.38 | 0.999991386714680 |

**VII) Conclusions**

Novel classes of PDEs have recently emerged and the physics of those PDEs requires keeping strict control of the curl of one or more vector fields. The PDEs are all hyperbolic systems of great interest to science and engineering. Many of the hyperbolic systems resulting from the Godunov-Peshkov-Romenski (GPR) formulation for hyperelasticity and compressible multiphase flow with and without surface tension have such curl-preserving update equations (Godunov and Romenski [32], Romenski [40], Romenski *et al*. [41], Peshkov and Romenski [37], [38], Dumbser *et al*. [27], [28], [31], Schmidmayer *et al*. [42]). The equations of General Relativity when cast in the FO-CCZ4 formulation also have such a structure (Alic *et al*. [1], [2], 2012, Brown *et al*. [18], Dumbser *et al*. [29], Dumbser, *et al*. [30]). Similarly, it has recently become possible to recast Schrödinger's equation in first order hyperbolic form, and the time-evolution of this very important equation also has curl-preserving constraints (Dhaouadi *et al*. [25], Busto *et al*. [19]). Experience has shown that if nothing special is done to account for the curl-preserving vector field, it can blow up in a finite amount of simulation time (Dumbser *et al*. [30]).

Prior work has shown that classical zone-centered Godunov methods can be adapted to such systems only if a GLM-type cleaning approach is included to suppress the build up of circulation on the mesh (Dumbser *et al*. [30]). The two-fold problem with this approach is that: 1) we often get a very large system of Lagrange multipliers that are not part of the original PDE system and 2) the signal speed with which the Lagrange multipliers have to be advected often exceeds the physical signal speed in the problem by a substantial margin. Another alternative is to solve an elliptical system at every timestep (Boscheri *et al*. [17]), which makes each timestep very expensive. In Balsara *et al*. [16] we first presented curl-preserving WENO-like methods that overcame both of the above-mentioned limitations. The methods were based on inventing a novel globally curl-preserving reconstruction strategy that reconstructs the vector field over the zone's volume by using the components of the vector field that were collocated at the edges of the mesh. Non-linear hybridization, via WENO methods, was seamlessly built into the curl-preserving reconstruction strategy. These edge-centered components were updated by using multidimensionally upwinded potentials that were collocated at the vertices of the mesh. Multidimensional Riemann solvers, designed by the first author, provided the requisite multidimensional upwinding. The resulting highly stable finite volume-like schemes for curl-preserving systems had the following desirable properties: **1)** They did not blow up even after very long integration times. **2)** They did not need GLM-style cleaning with very high signal speeds. **3)** They could operate with large explicit timesteps. **4)** They did not require the solution of an elliptic system. and **5)** They could be extended to higher orders while incorporating non-linear hybridization by using WENO-like methods. It is, therefore, very desirable to invent DG-like and PNPM-like variants of these WENO-like schemes so that they can inherit the same desirable features – such a task is fulfilled on in this paper.

Since we know that DG schemes and PNPM schemes provide more accurate alternatives to WENO schemes, it becomes interesting to design curl-free and curl-preserving variants of the such schemes. In this paper, we present for the very first time, globally curl-preserving DG-like and PNPM-like schemes that share the beneficial traits of the globally curl-preserving WENO-like schemes designed in Balsara *et al*. [16]. The higher moments of the vector components that live in the edges of the mesh are, therefore, endowed with time-evolution that is consistent with the governing equations. This is accomplished by making a Galerkin projection within each edge that results in a weak form of update equation for the higher order edge-centered moments. The update utilizes the multidimensionally upwinded potentials at the vertices of the mesh. Such update equations have been documented in Section II for the model eqn. (1.2) and some nuances of the curl-preserving reconstruction, and how it relates to traditional DG schemes, is brought out in Section III.

It is well-known that zone-centered DG schemes have wave propagation characteristics that are superior to zone-centered WENO schemes. Such a superior behavior can be revealed by carrying out a von-Neumann stability analysis of either scheme and inter-comparing the results. It is, therefore, very interesting to carry out a von Neumann stability analysis of our newly-developed globally curl-free and curl-preserving DG-like and PNPM-like schemes. In Section IV we present details of our von Neumann stability analysis. To bring out the curl-free aspect of the evolution, the analysis must absolutely be done in two or more dimensions. Therefore, our analysis is two-dimensional by its very design. By pushing the capabilities of computer algebra systems to the limits, we have been able to extend this two-dimensional curl-preserving von Neumann stability analysis up to fourth order of accuracy.

Section V shows the results of this von Neumann stability analysis. We present such an analysis for globally curl-free WENO-like, PNPM-like and DG-like schemes to facilitate inter-comparison. In Sub-section V.a the limiting CFL numbers for all these schemes are derived and documented in Tables I, II and III. We find, unsurprisingly, that WENO-like schemes offer the largest CFL numbers while DG-like schemes restrict us to substantially smaller CFL numbers. The PNPM-like schemes give us quite large CFL numbers at a much-reduced computational complexity. In Sub-section V.b we document the dissipation and dispersion properties of the same three schemes. We do this for waves with wavelengths that are 5, 10 and 15 times the zone size. For each family of schemes, the dissipation and dispersion properties do indeed improve with increasing order of accuracy, as expected. This is shown in the figures associated with Sub-section V.b and also in Tables IV, V and VI. We find that WENO-like schemes have dissipation and dispersion properties that are noticeably inferior to the DG-like schemes at the same order. However, we find that PNPM-like schemes have dissipation and dispersion properties that approach that of DG-like schemes at the same order while offering substantially larger CFL numbers.

Section VI presents numerical results, where we show that our methods meet their design accuracies. We also show that with increasing order of accuracy the methods become very good at preserving quadratic energy. This is a welcome result, because the methods were not intentionally designed to preserve quadratic energy; yet they seem to do a good job of doing so. When the evolution of the PDE is curl-free, our methods also hold down the discrete circulation to machine accuracy over long integration times. The importance of this fact in the design of curl-preserving schemes is also discussed.

This paper has laid the essential foundation for several novel globally curl constraint-preserving methods and catalogued their very many desirable properties. The next step would be to apply them to full PDE systems where their potential gains can be realized.

**Acknowledgements**

DSB acknowledges support via NSF grants NSF-19-04774, NSF-AST-2009776 and NASA-2020-1241.

**Ethical Statement**

**i. Compliance with Ethical Standards** : This manuscript complies with all ethical standards for scientific publishing.
**ii. (in case of Funding) Funding** : The funding has been acknowledged. DSB acknowledges support via NSF grants NSF-19-04774, NSF-AST-2009776 and NASA-2020-1241
**iii. Conflict of Interest** : The authors have no conflict of interest.
**iv. Ethical approval** : N/A
**v. Informed consent** : N/A

**Appendix A**

Here we explicitly provide the 9 coefficients of the 3×3 matrix "**A**" in eqn. (4.3). This should enable the reader to cross-check her or his mathematics. The first row is given by

$$A_{1,1} = \frac{\left(\Delta x \cos(\Delta y\, k_y) - \Delta x\right)|v_y| - i\Delta x \sin(\Delta y\, k_y) v_y + \left(\Delta y \cos(\Delta x\, k_x) - \Delta y\right)|v_x| - i\Delta y \sin(\Delta x\, k_x) v_x}{\Delta x\, \Delta y}$$

$$A_{1,2} = -\frac{i\sin(\Delta x\, k_x)|v_x| + \left(1 - \cos(\Delta x\, k_x)\right) v_x}{2\Delta x}$$

$$A_{1,3} = -\frac{\left(i\sin\left(\frac{\Delta y\, k_y + \Delta x\, k_x}{2}\right) - i\sin\left(\frac{\Delta y\, k_y - \Delta x\, k_x}{2}\right)\right)|v_y| + \left(\cos\left(\frac{\Delta y\, k_y - \Delta x\, k_x}{2}\right) - \cos\left(\frac{\Delta y\, k_y + \Delta x\, k_x}{2}\right)\right) v_y}{2\Delta x}$$

.

The second row is given by

$$A_{2,1} = \frac{6i\sin(\Delta x\, k_x)|v_x| + \left(6 - 6\cos(\Delta x\, k_x)\right) v_x}{\Delta x}$$

$$A_{2,2} = \frac{\left(\cos(\Delta y\, k_y) - 1\right)|v_y| - i\sin(\Delta y\, k_y) v_y + \left(-3\cos(\Delta x\, k_x) - 3\right)|v_x| + 3i\sin(\Delta x\, k_x) v_x}{\Delta x}$$

$$A_{2,3} = 0.$$

The third row is given by

$$A_{3,1} = \frac{3\Delta x\, e^{3i\Delta y k_y + \frac{i\Delta x k_x}{2}}|v_y|}{\Delta y^2 e^{\frac{3i\Delta y k_y}{2}+i\Delta x k_x} - \Delta y^2 e^{\frac{3i\Delta y k_y}{2}}} - \frac{3\Delta x\, e^{2i\Delta y k_y + \frac{i\Delta x k_x}{2}}|v_y|}{\Delta y^2 e^{\frac{3i\Delta y k_y}{2}+i\Delta x k_x} - \Delta y^2 e^{\frac{3i\Delta y k_y}{2}}}$$

$$-\frac{3\Delta x\, e^{i\Delta y k_y + \frac{i\Delta x k_x}{2}}|v_y|}{\Delta y^2 e^{\frac{3i\Delta y k_y}{2}+i\Delta x k_x} - \Delta y^2 e^{\frac{3i\Delta y k_y}{2}}} + \frac{3\Delta x\, e^{\frac{i\Delta x k_x}{2}}|v_y|}{\Delta y^2 e^{\frac{3i\Delta y k_y}{2}+i\Delta x k_x} - \Delta y^2 e^{\frac{3i\Delta y k_y}{2}}}$$

$$-\frac{3\Delta x\, e^{3i\Delta y k_y + \frac{i\Delta x k_x}{2}}v_y}{\Delta y^2 e^{\frac{3i\Delta y k_y}{2}+i\Delta x k_x} - \Delta y^2 e^{\frac{3i\Delta y k_y}{2}}} + \frac{9\Delta x\, e^{2i\Delta y k_y + \frac{i\Delta x k_x}{2}}v_y}{\Delta y^2 e^{\frac{3i\Delta y k_y}{2}+i\Delta x k_x} - \Delta y^2 e^{\frac{3i\Delta y k_y}{2}}}$$

$$-\frac{9\Delta x\, e^{i\Delta y k_y + \frac{i\Delta x k_x}{2}}v_y}{\Delta y^2 e^{\frac{3i\Delta y k_y}{2}+i\Delta x k_x} - \Delta y^2 e^{\frac{3i\Delta y k_y}{2}}} + \frac{3\Delta x\, e^{\frac{i\Delta x k_x}{2}}v_y}{\Delta y^2 e^{\frac{3i\Delta y k_y}{2}+i\Delta x k_x} - \Delta y^2 e^{\frac{3i\Delta y k_y}{2}}}$$

$$A_{3,2} = 0$$

$$A_{3,3} = -\frac{\left(3\cos(\Delta y k_y)+3\right)|v_y| - 3i\sin(\Delta y k_y)v_y + \left(1-\cos(\Delta x k_x)\right)|v_x| + i\sin(\Delta x k_x)v_x}{\Delta y}$$

This completes our description of the 9 coefficients of the 3×3 matrix "**A**" in eqn. (4.3).

**Figure Captions**

*Fig. 1 shows the components of the curl-free vector field around the four zones centered around (i,j), (i-1,j), (i-1,j-1) and (i,j-1). A first order curl-free reconstruction is used. The multidimensionally upwinded potentials at the vertices of the zone (i,j) are also shown for the situation where both components of the velocity are positive.*

*Fig. 2 shows the collocation of vector components along the edges of a two-dimensional control volume. As evaluated over the edges of the square element, the discrete circulation is fully specified. The mean value of the vector components and their linear variation are shown along each edge, in keeping with a second order accurate reconstruction scheme. The reconstruction problem for a curl-preserving reconstruction consists of obtaining a polynomial-based vector field that matches the specified mean circulation in the zone while simultaneously matching the edge values within each zone.*

*Fig. 3 shows how the facially collocated Fourier modes associated with the curl-free vector field relate to one another across the different faces of the mesh. These Fourier modes, and their analogues at all the other faces in the figure, are used for carrying out the von Neumann stability analysis.*

*Fig. 4 shows the domain of stability for a) a second order in space and time curl-free DG-like scheme that uses SSK-RK2 timestepping, b) a third order in space and time curl-free DG-like scheme that uses SSK-RK3 timestepping, and c) a fourth order in space and time curl-free DG-like scheme that uses SSP-RK(5,4) timestepping. The CFL numbers in the x- and y-directions are denoted by $C_x$ and $C_y$ and the color coding shows the absolute value of the largest eigenvalue of the amplification matrix. The white polygon identifies the full domain of stability and the white circle identifies the largest circle that can be fit within the domain of stability. The radius of the white circle, therefore, gives us the maximal CFL number.*

*Fig. 5 shows the wave propagation characteristics for curl-preserving second order WENO-like schemes. Figs. 5a to 5d show one minus the absolute value of the amplification factor when the velocity vector makes angles of 0º , 15º , 30º and 45º relative to the x-direction of the 2D mesh. Figs. 5e to 5h show the phase error, again for the same angles. The 2D wave vector can make any angle relative to the 2D direction of velocity propagation, therefore, the amplitude and phase information are shown w.r.t. the angle made between the velocity direction and the direction of the wave vector. In each plot, the blue curve refers to waves that span 5 cells per wavelength; the green curve refers to waves that span 10 cells per wavelength; the red curve refers to waves that span 15 waves per wavelength.*

*Fig. 6 shows the wave propagation characteristics for curl-preserving third order WENO-like schemes. Figs. 6a to 6d show one minus the absolute value of the amplification factor when the velocity vector makes angles of 0º , 15º , 30º and 45º relative to the x-direction of the 2D mesh. Figs. 6e to 6h show the phase error, again for the same angles. The 2D wave vector can make any angle relative to the 2D direction of velocity propagation, therefore, the amplitude and phase information are shown w.r.t. the angle made between the velocity direction and the direction of the wave vector. In each plot, the blue curve refers to waves that span 5 cells per wavelength; the green curve refers to waves that span 10 cells per wavelength; the red curve refers to waves that span 15 waves per wavelength.*

*Fig. 7 shows the wave propagation characteristics for curl-preserving fourth order WENO-like schemes. Figs. 7a to 7d show one minus the absolute value of the amplification factor when the velocity vector makes angles of 0º , 15º , 30º and 45º relative to the x-direction of the 2D mesh. Figs. 7e to 7h show the phase error, again for the same angles. The 2D wave vector can make any angle relative to the 2D direction of velocity propagation, therefore, the amplitude and phase information are shown w.r.t. the angle made between the velocity direction and the direction of the wave vector. In each plot, the blue curve refers to waves that span 5 cells per wavelength; the green curve refers to waves that span 10 cells per wavelength; the red curve refers to waves that span 15 waves per wavelength.*

*Fig. 8 shows the wave propagation characteristics for curl-preserving third order P1P2-like schemes. Figs. 8a to 8d show one minus the absolute value of the amplification factor when the velocity vector makes angles of 0º , 15º , 30º and 45º relative to the x-direction of the 2D mesh. Figs. 8e to 8h show the phase error, again for the same angles. The 2D wave vector can make any angle relative to the 2D direction of velocity propagation, therefore, the amplitude and phase information are shown w.r.t. the angle made between the velocity direction and the direction of the wave vector. In each plot, the blue curve refers to waves that span 5 cells per wavelength; the green curve refers to waves that span 10 cells per wavelength; the red curve refers to waves that span 15 waves per wavelength.*

*Fig. 9 shows the wave propagation characteristics for curl-preserving fourth order P1P3-like schemes. Figs. 9a to 9d show one minus the absolute value of the amplification factor when the velocity vector makes angles of 0º , 15º , 30º and 45º relative to the x-direction of the 2D mesh. Figs. 9e to 9h show the phase error, again for the same angles. The 2D wave vector can make any angle relative to the 2D direction of velocity propagation, therefore, the amplitude and phase information are shown w.r.t. the angle made between the velocity direction and the direction of the wave vector. In each plot, the blue curve refers to waves that span 5 cells per wavelength; the green curve refers to waves that span 10 cells per wavelength; the red curve refers to waves that span 15 waves per wavelength.*

*Fig. 10 shows the wave propagation characteristics for curl-preserving second order DG-like schemes. Figs. 10a to 10d show one minus the absolute value of the amplification factor when the velocity vector makes angles of 0º , 15º , 30º and 45º relative to the x-direction of the 2D mesh. Figs. 10e to 10h show the phase error, again for the same angles. The 2D wave vector can make any angle relative to the 2D direction of velocity propagation, therefore, the amplitude and phase information are shown w.r.t. the angle made between the velocity direction and the direction of the wave vector. In each plot, the blue curve refers to waves that span 5 cells per wavelength; the green curve refers to waves that span 10 cells per wavelength; the red curve refers to waves that span 15 waves per wavelength.*

*Fig. 11 shows the wave propagation characteristics for curl-preserving third order DG-like schemes. Figs. 11a to 11d show one minus the absolute value of the amplification factor when the velocity vector makes angles of 0º , 15º , 30º and 45º relative to the x-direction of the 2D mesh. Figs. 11e to 11h show the phase error, again for the same angles. The 2D wave vector can make any angle relative to the 2D direction of velocity propagation, therefore, the amplitude and phase information are shown w.r.t. the angle made between the velocity direction and the direction of the wave vector. In each plot, the blue curve refers to waves that span 5 cells per wavelength; the green curve refers to waves that span 10 cells per wavelength; the red curve refers to waves that span 15 waves per wavelength.*

*Fig. 12 shows the wave propagation characteristics for curl-preserving fourth order DG-like schemes. Figs. 12a to 12d show one minus the absolute value of the amplification factor when the velocity vector makes angles of 0º , 15º , 30º and 45º relative to the x-direction of the 2D mesh. Figs. 12e to 12h show the phase error, again for the same angles. The 2D wave vector can make any angle relative to the 2D direction of velocity propagation, therefore, the amplitude and phase information are shown w.r.t. the angle made between the velocity direction and the direction of the wave vector. In each plot, the blue curve refers to waves that span 5 cells per wavelength; the green curve refers to waves that span 10 cells per wavelength; the red curve refers to waves that span 15 waves per wavelength.*

*Fig. 13 shows the quadratic field energy from the vortex problem that is preserved on the mesh at the final time in the simulation as a function of mesh size. Panel a) displays the curl-free DG-like schemes, panel b) displays the curl-free P1PN-like schemes and panel c) displays the curl-free WENO-like schemes.*

*Fig. 14 shows the maximum pointwise error of the curl of **J** as a function of time for a 64×64 zone run of the vortex problem. Fig. 14a shows the evolution of the maximum pointwise curl as a function of time for the $2^{nd}$ , $3^{rd}$ and $4^{th}$ order curl-free DG-like schemes. Fig. 14b shows the evolution of the maximum pointwise curl as a function of time for the $3^{rd}$ and $4^{th}$ order curl-free*

*P1PN-like schemes. Fig. 14c shows the evolution of the maximum pointwise curl as a function of time for the $2^{nd}$, $3^{rd}$ and $4^{th}$ order curl-free WENO-like schemes. The figure shows that all our curl-preserving schemes can preserve the curl constraint up to machine accuracy.*

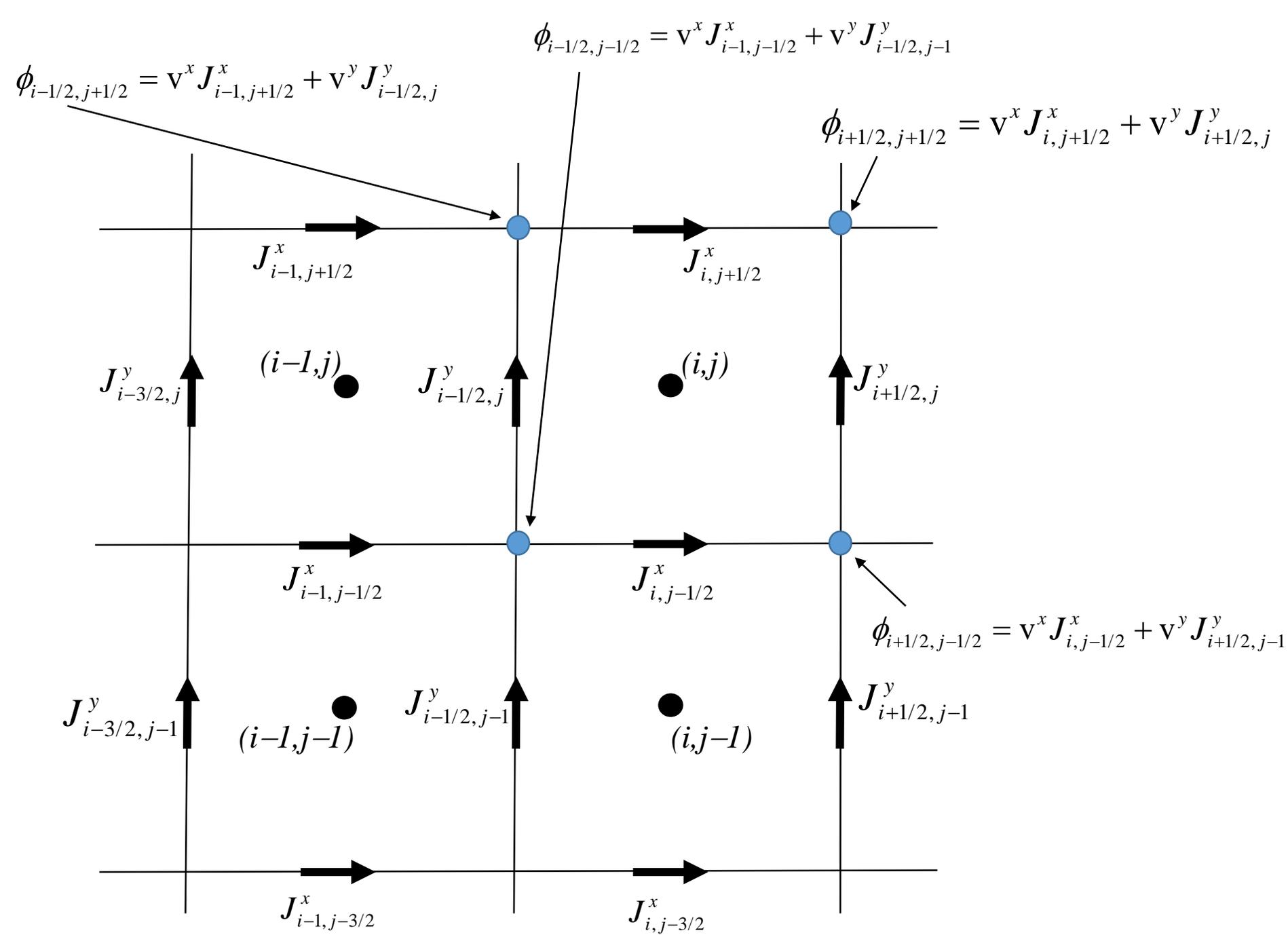

*Fig. 1 shows the components of the curl-free vector field around the four zones centered around (i,j), (i-1,j), (i-1,j-1) and (i,j-1). A first order curl-free reconstruction is used in this figure; though the use of a higher order reconstruction of the vector field is also possible. The multidimensionally upwinded potentials at the vertices of the zone (i,j) are also shown for the situation where both components of the velocity are positive. Again, while the potentials are shown explicitly for the first order case in the figure, a higher order reconstruction will yield more accurate potentials at the vertices of the mesh.*

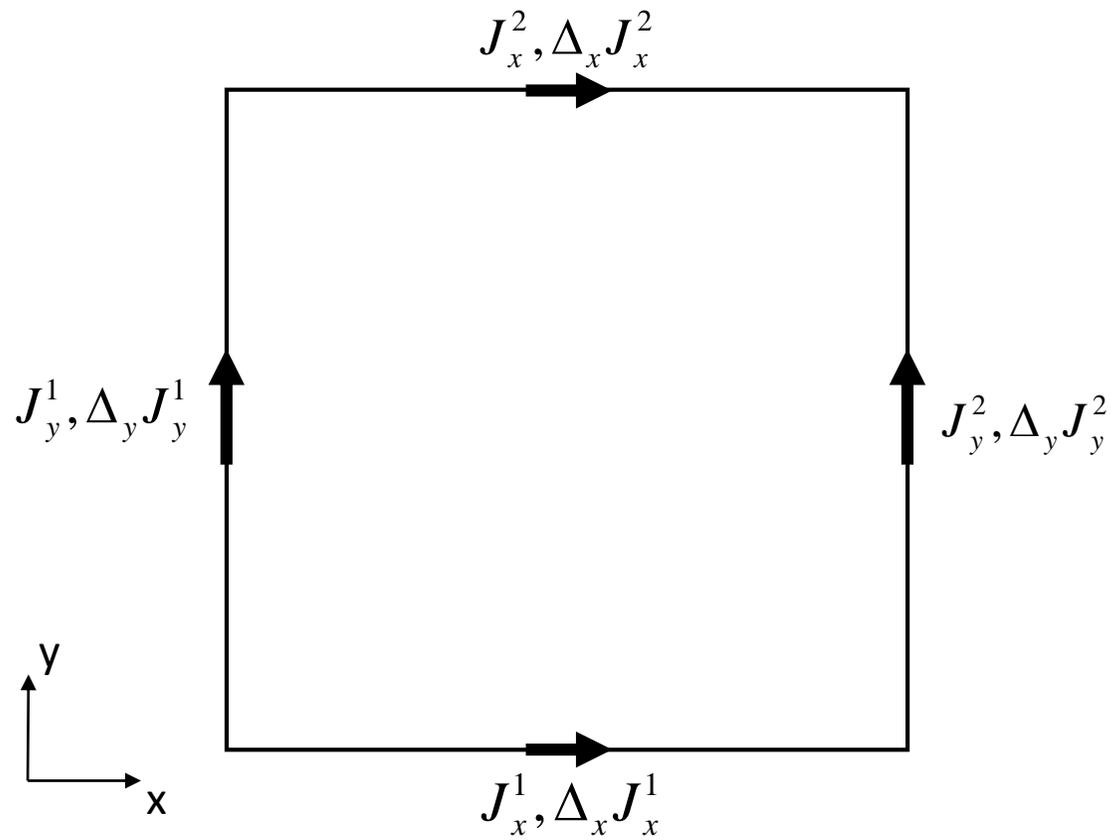

Fig. 2 shows the collocation of vector components along the edges of a two-dimensional control volume. As evaluated over the edges of the square element, the discrete circulation is fully specified. The mean value of the vector components and their linear variation are shown along each edge, in keeping with a second order accurate reconstruction scheme. The reconstruction problem for a curl-preserving reconstruction consists of obtaining a polynomial-based vector field that matches the specified mean circulation in the zone while simultaneously matching the edge values within each zone.

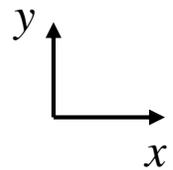

| | $J_0^{x+}(t)e^{+ik_y\Delta y}, J_x^{x+}(t)e^{+ik_y\Delta y}$ | |
|---|---|---|
| (−1,1) | (0,1) | (1,1) |
| | $J_0^{x+}(t), J_x^{x+}(t)$ | |
| $J_0^{y+}(t)e^{-2ik_x\Delta x}$ $J_y^{y+}(t)e^{-2ik_x\Delta x}$ (−1,0) | $J_0^{y+}(t)e^{-ik_x\Delta x}$ $J_y^{y+}(t)e^{-ik_x\Delta x}$ (0,0) | $J_0^{y+}(t)$ $J_y^{y+}(t)$ (1,0)    $J_0^{y+}(t)e^{+ik_x\Delta x}$ $J_y^{y+}(t)e^{+ik_x\Delta x}$ |
| | $J_0^{x+}(t)e^{-ik_y\Delta y}, J_x^{x+}(t)e^{-ik_y\Delta y}$ | |
| (−1,−1) | (0,−1) | (1,−1) |
| | $J_0^{x+}(t)e^{-2ik_y\Delta y}, J_x^{x+}(t)e^{-2ik_y\Delta y}$ | |

*Fig. 3 shows how the facially collocated Fourier modes associated with the curl-free vector field relate to one another across the different faces of the mesh. These Fourier modes, and their analogues at all the other faces in the figure, are used for carrying out the von Neumann stability analysis.*

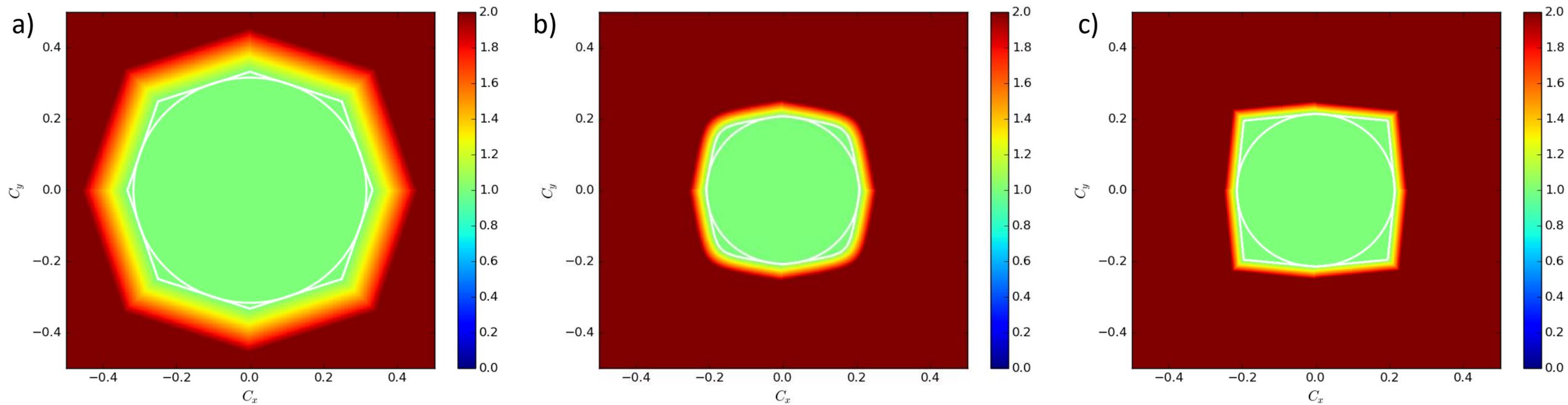

*Fig. 4 shows the domain of stability for a) a second order in space and time curl-free DG-like scheme that uses SSK-RK2 timestepping, b) a third order in space and time curl-free DG-like scheme that uses SSK-RK3 timestepping, and c) a fourth order in space and time curl-free DG-like scheme that uses SSP-RK(5,4) timestepping. The CFL numbers in the x- and y-directions are denoted by $C_x$ and $C_y$ and the color coding shows the absolute value of the largest eigenvalue of the amplification matrix. The white polygon identifies the full domain of stability and the white circle identifies the largest circle that can be fit within the domain of stability. The radius of the white circle, therefore, gives us the maximal CFL number.*

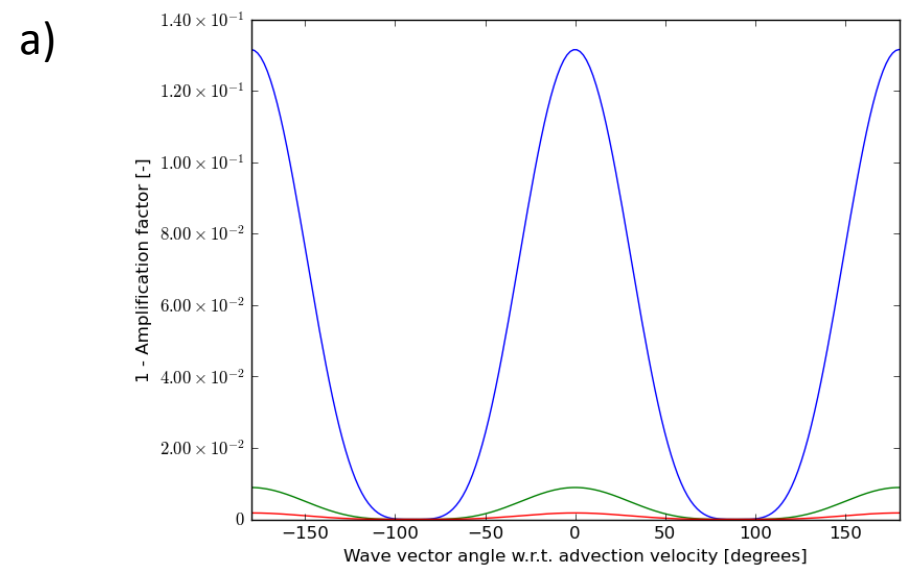
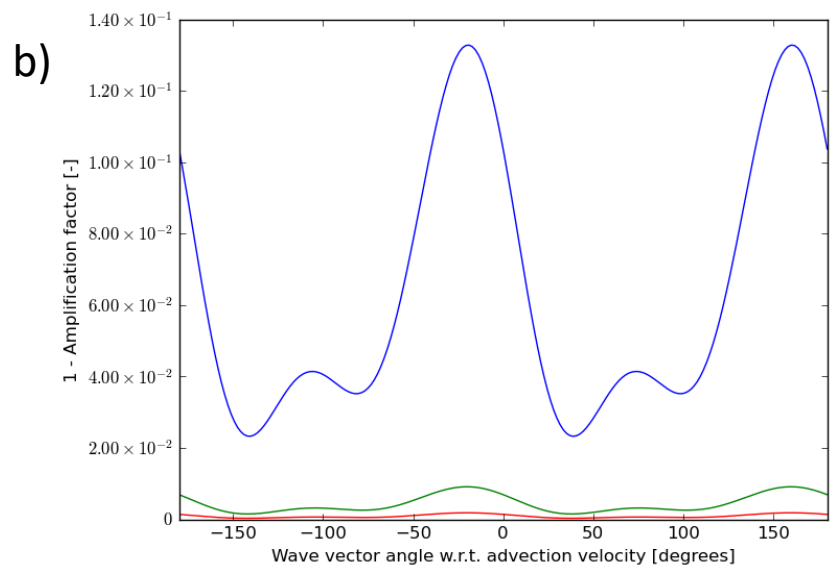
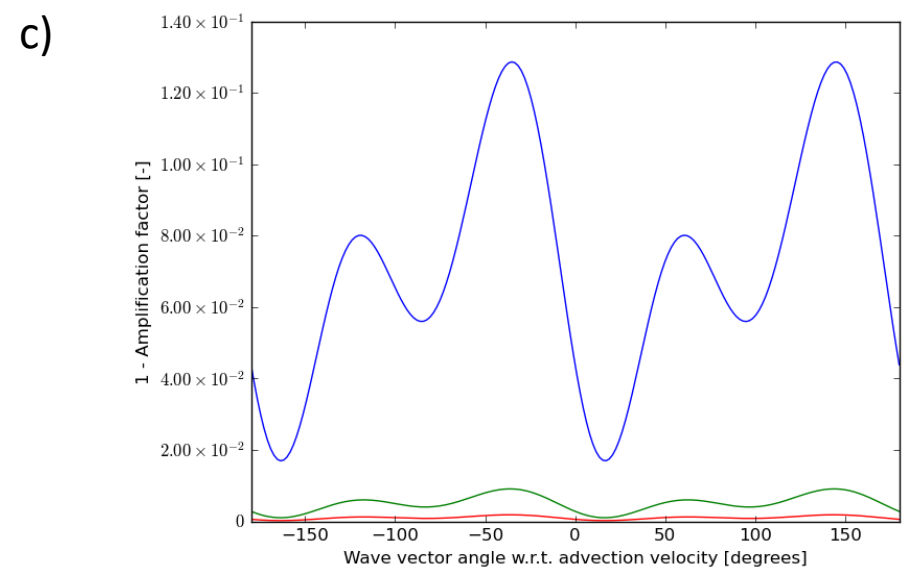
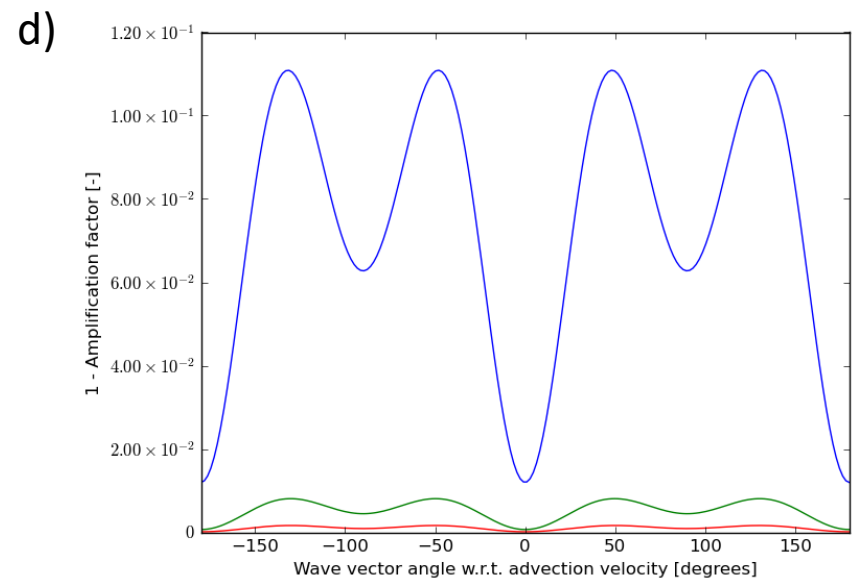

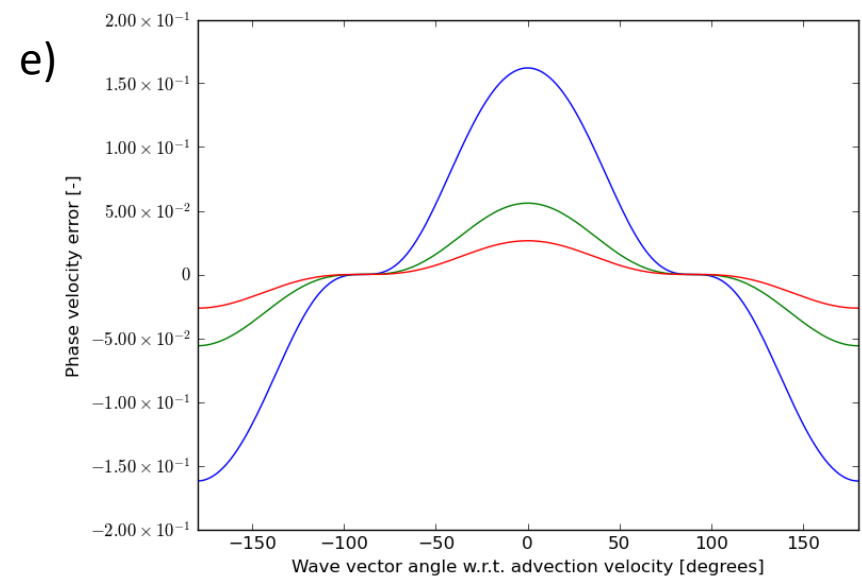
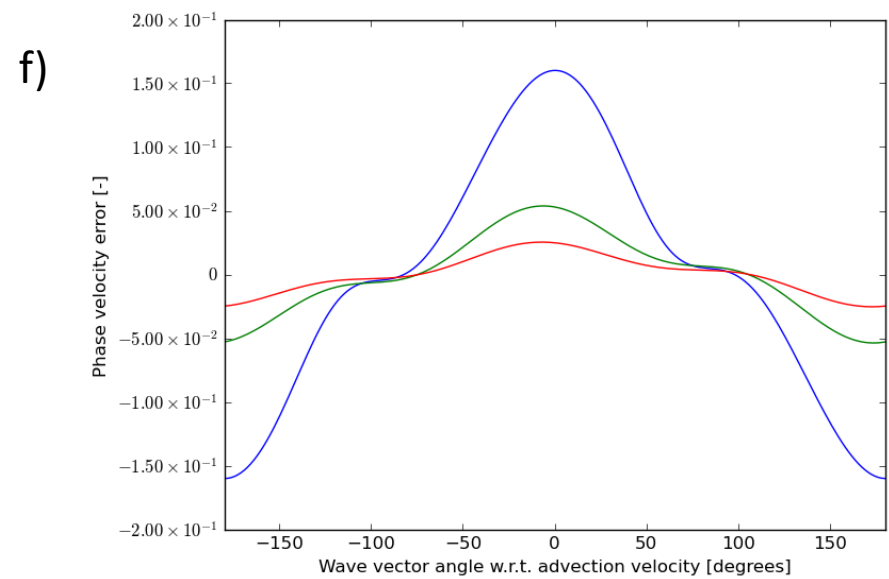
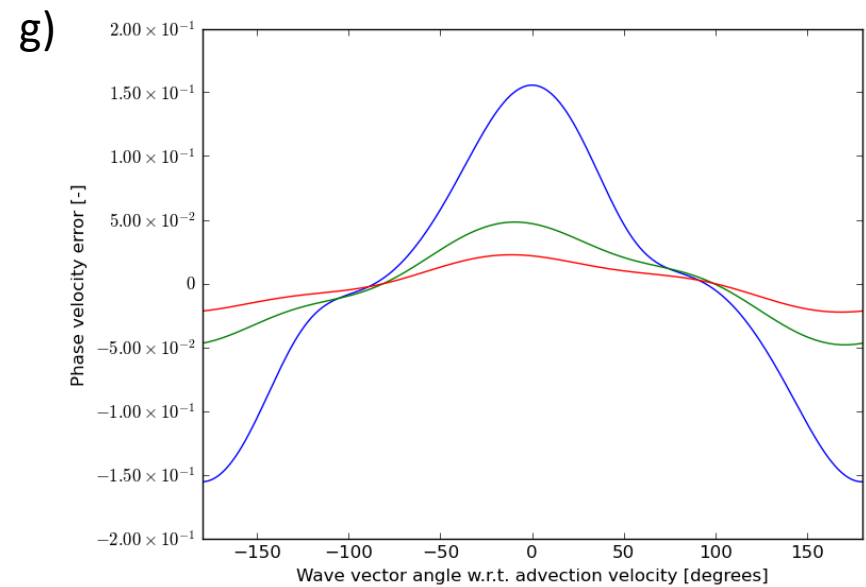
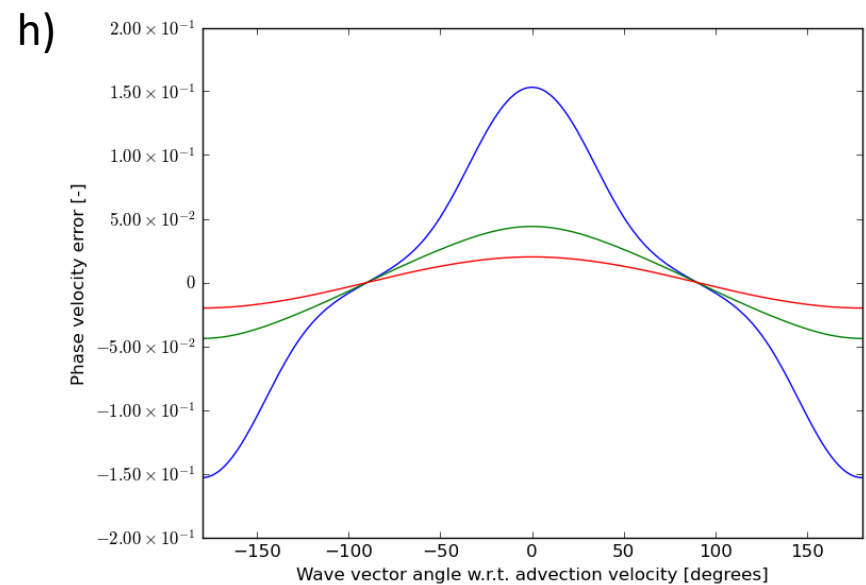

*Fig. 5 shows the wave propagation characteristics for curl-preserving second order WENO-like schemes. Figs. 5a to 5d show one minus the absolute value of the amplification factor when the velocity vector makes angles of 0º , 15º , 30º and 45º relative to the x-direction of the 2D mesh. Figs. 5e to 5h show the phase error, again for the same angles. The 2D wave vector can make any angle relative to the 2D direction of velocity propagation, therefore, the amplitude and phase information are shown w.r.t. the angle made between the velocity direction and the direction of the wave vector. In each plot, the blue curve refers to waves that span 5 cells per wavelength; the green curve refers to waves that span 10 cells per wavelength; the red curve refers to waves that span 15 waves per wavelength.*

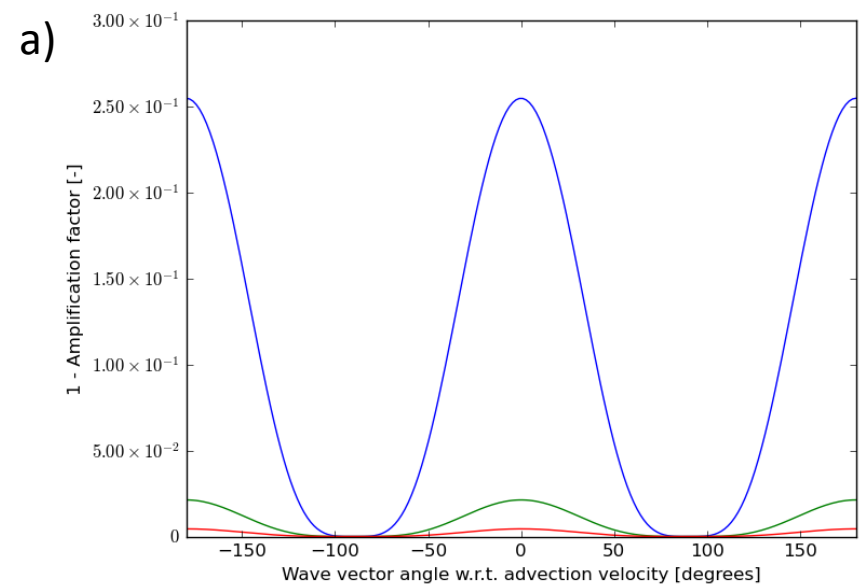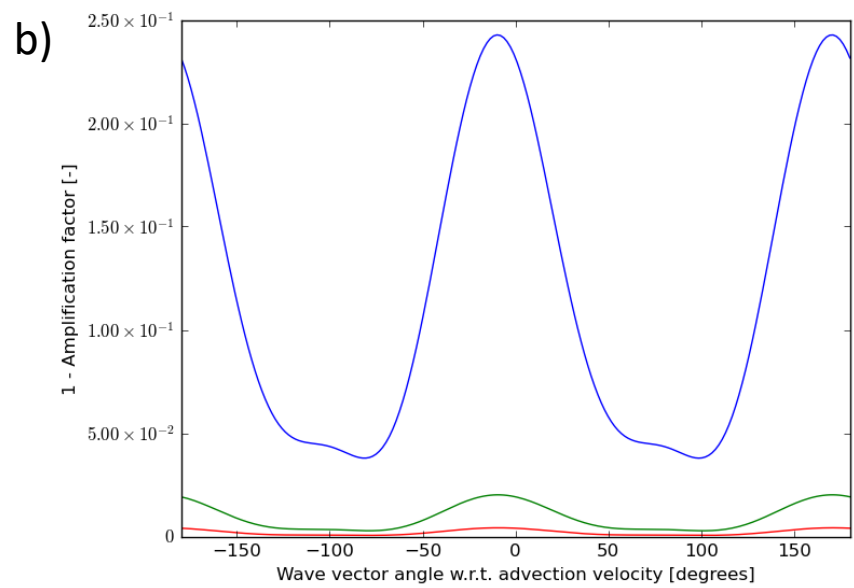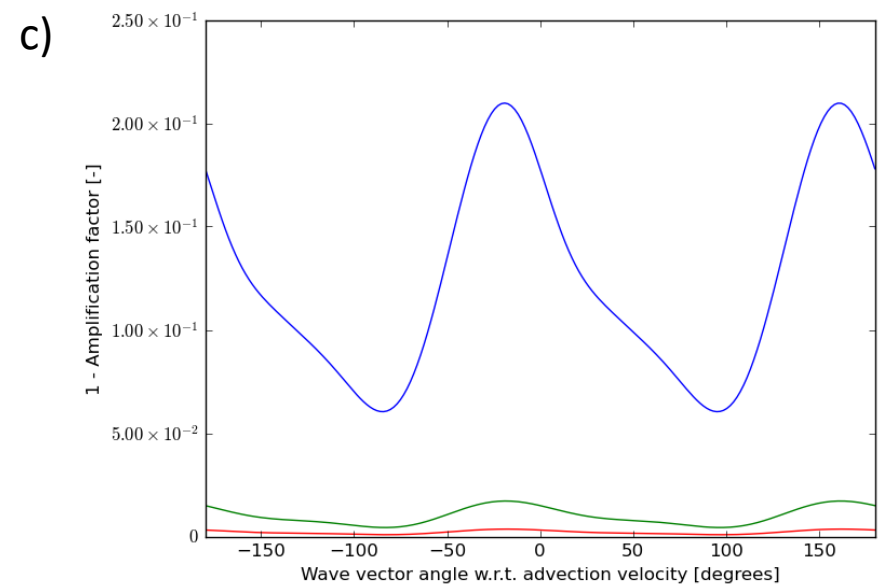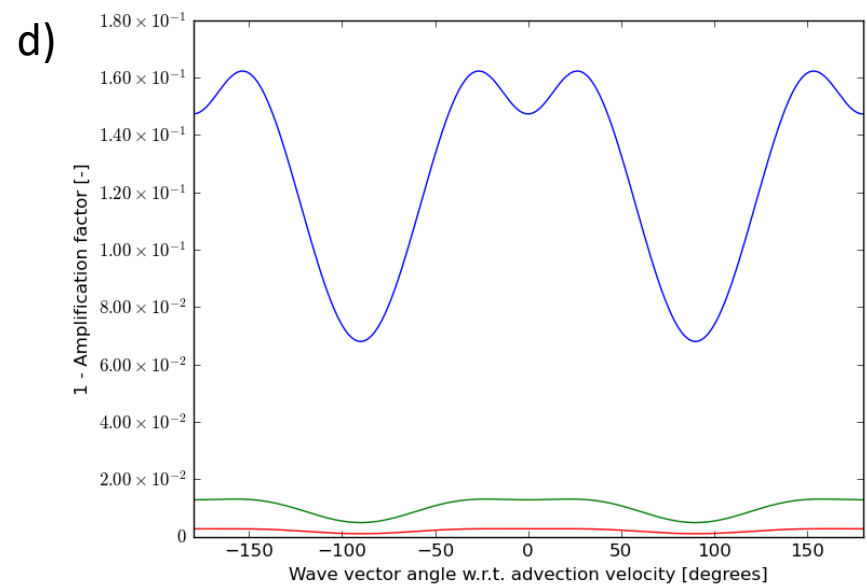

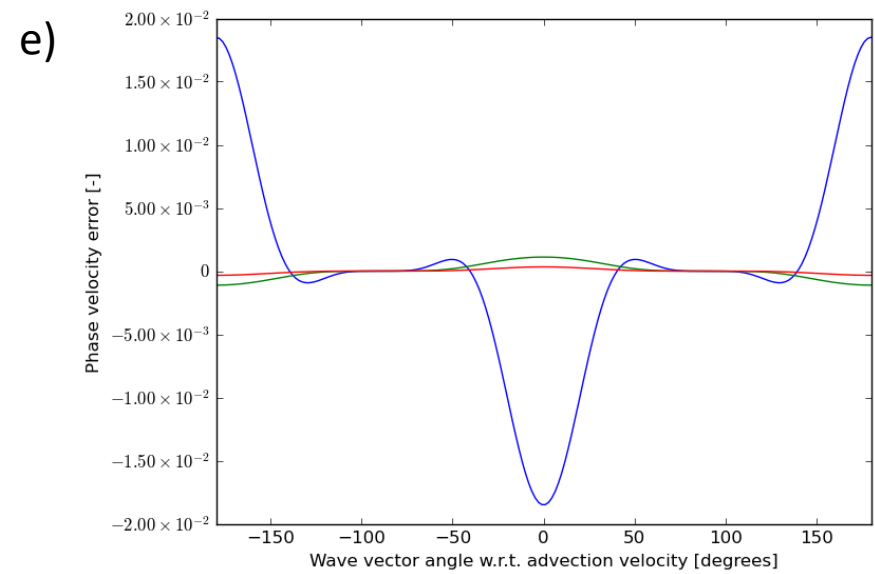
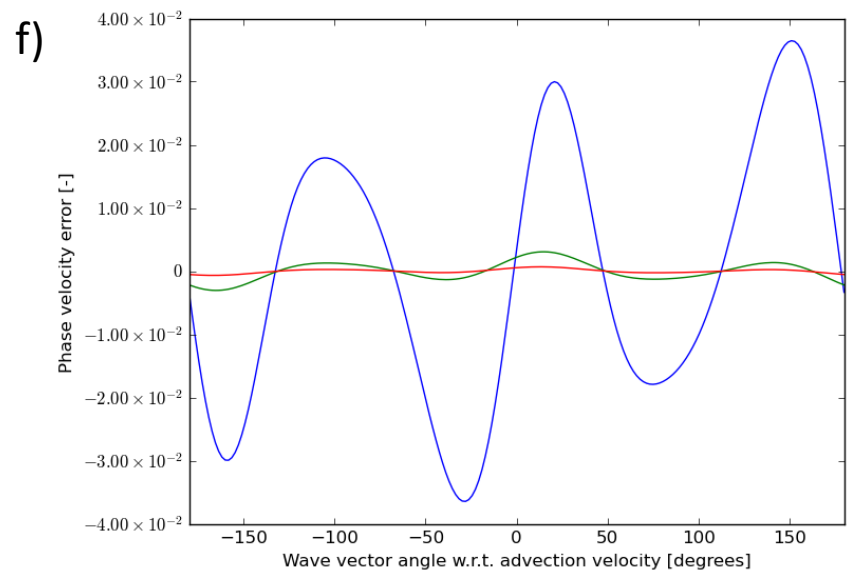
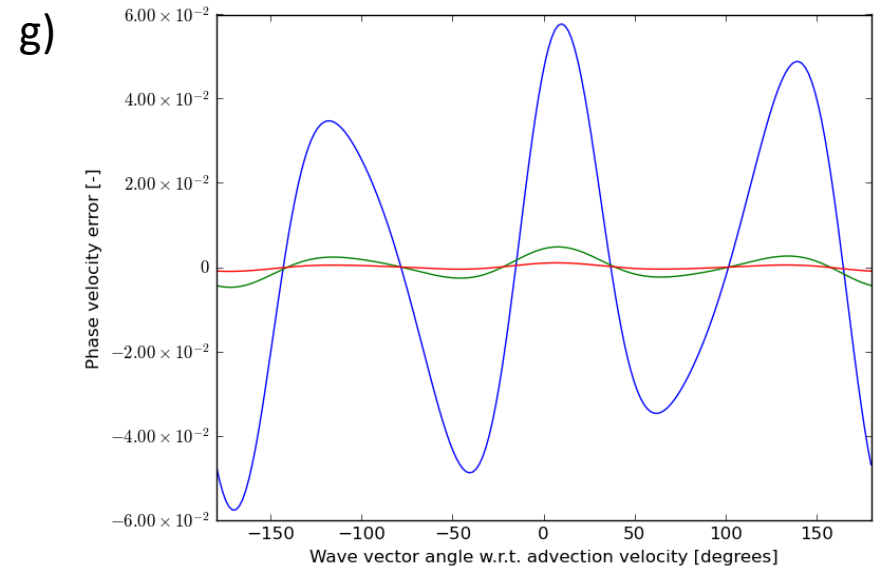
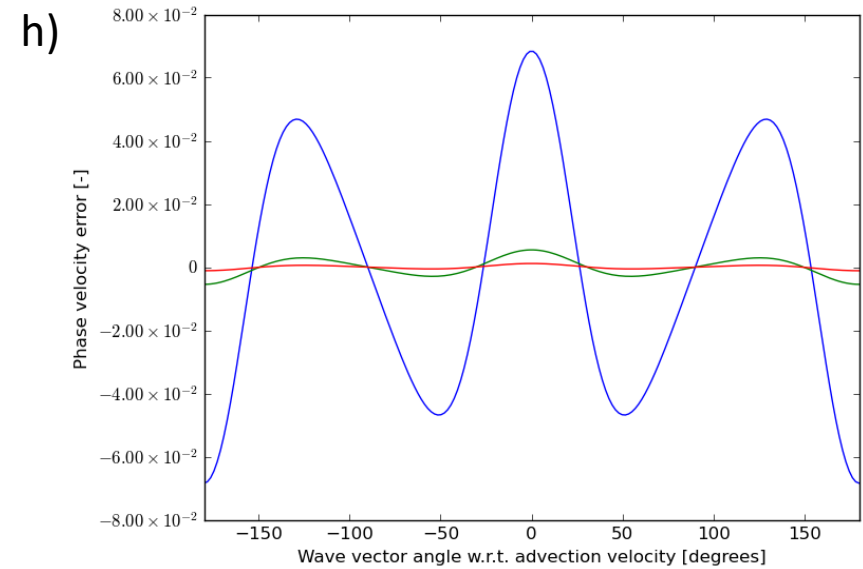

*Fig. 6 shows the wave propagation characteristics for curl-preserving third order WENO-like schemes. Figs. 6a to 6d show one minus the absolute value of the amplification factor when the velocity vector makes angles of 0º, 15º, 30º and 45º relative to the x-direction of the 2D mesh. Figs. 6e to 6h show the phase error, again for the same angles. The 2D wave vector can make any angle relative to the 2D direction of velocity propagation, therefore, the amplitude and phase information are shown w.r.t. the angle made between the velocity direction and the direction of the wave vector. In each plot, the blue curve refers to waves that span 5 cells per wavelength; the green curve refers to waves that span 10 cells per wavelength; the red curve refers to waves that span 15 waves per wavelength.*

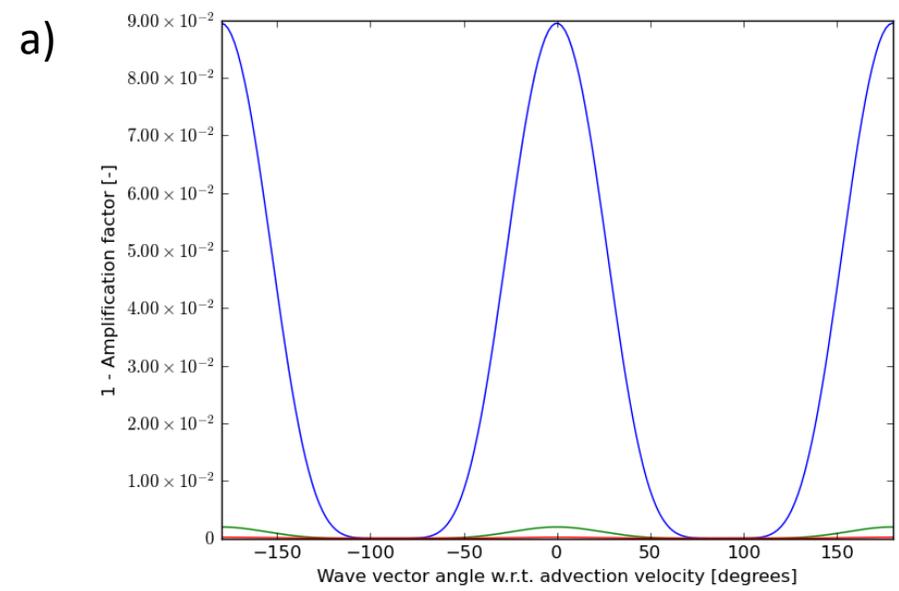
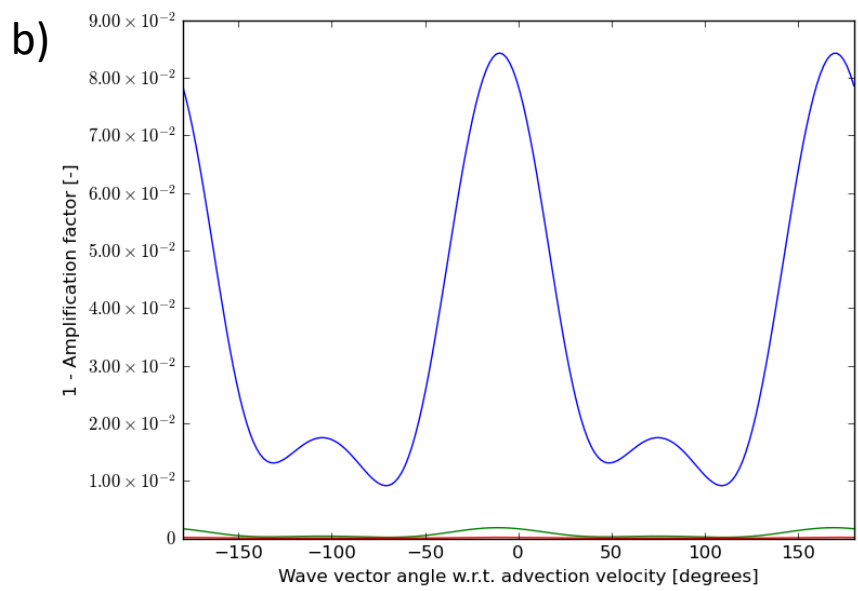
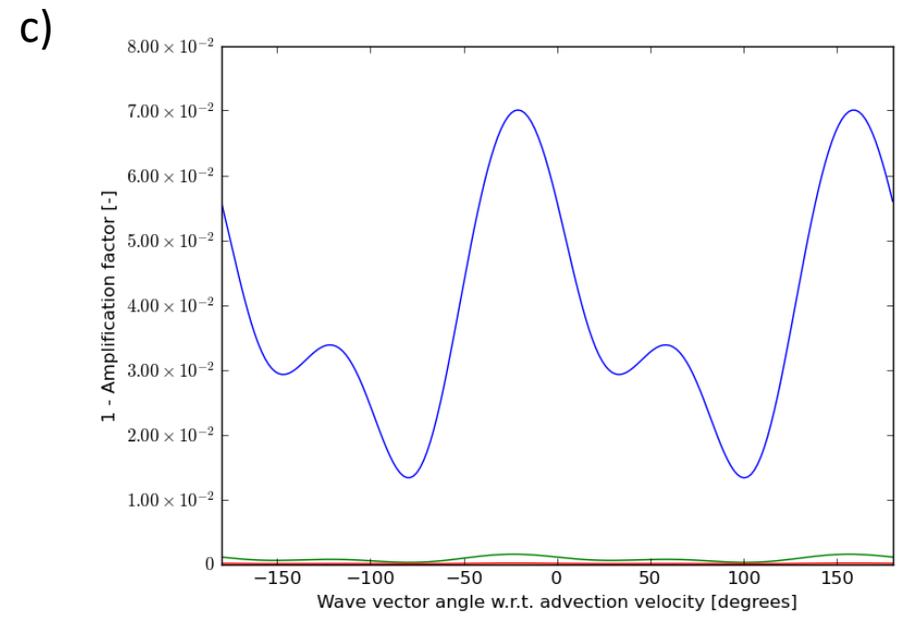
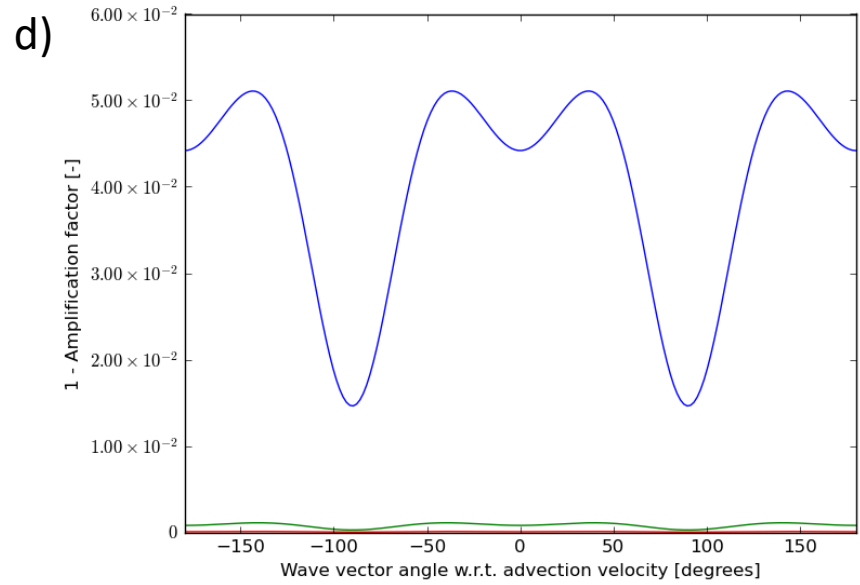

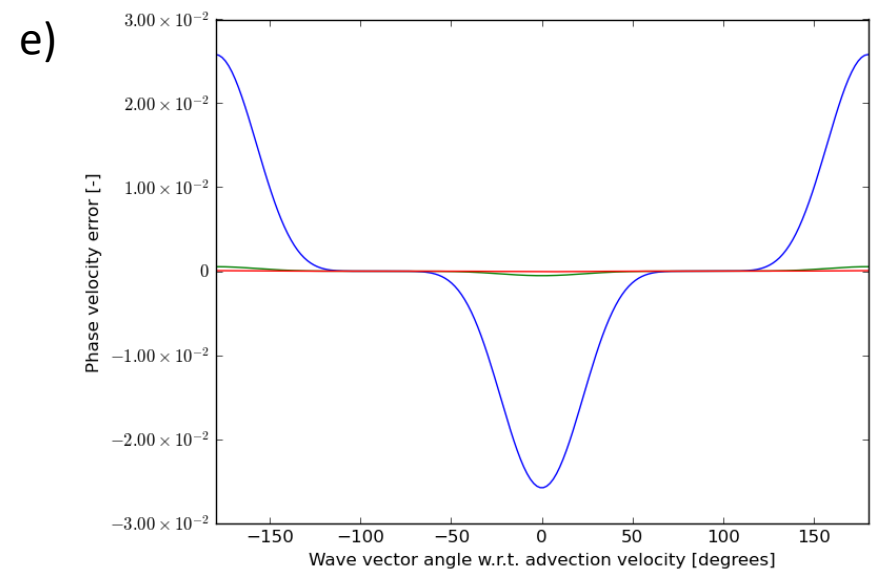
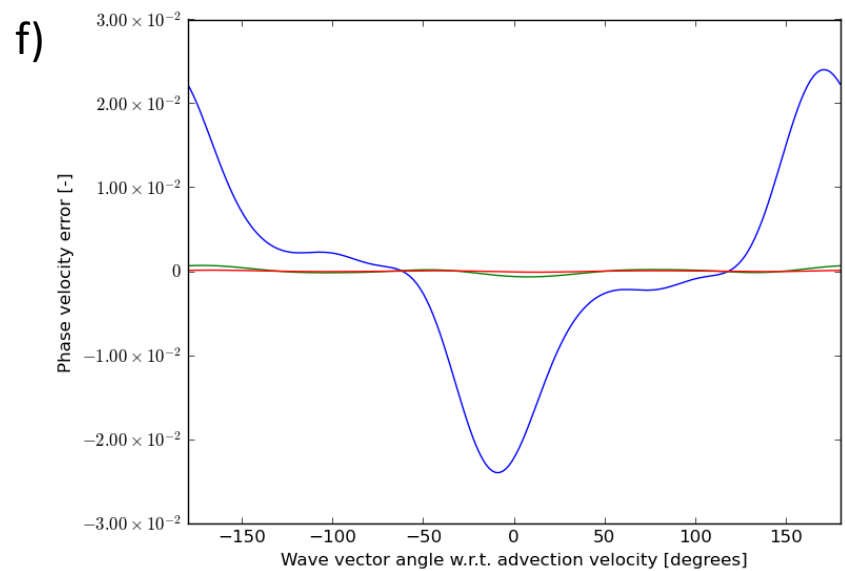
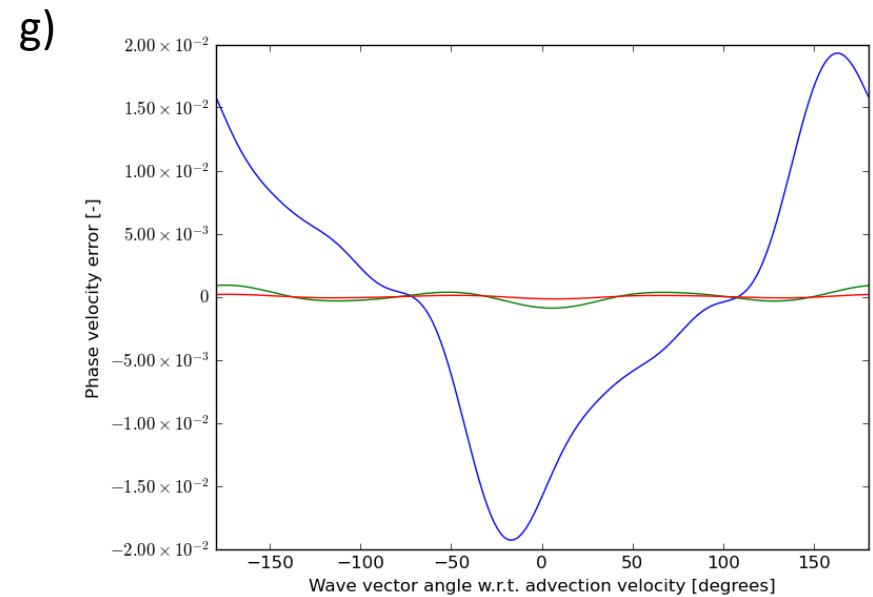
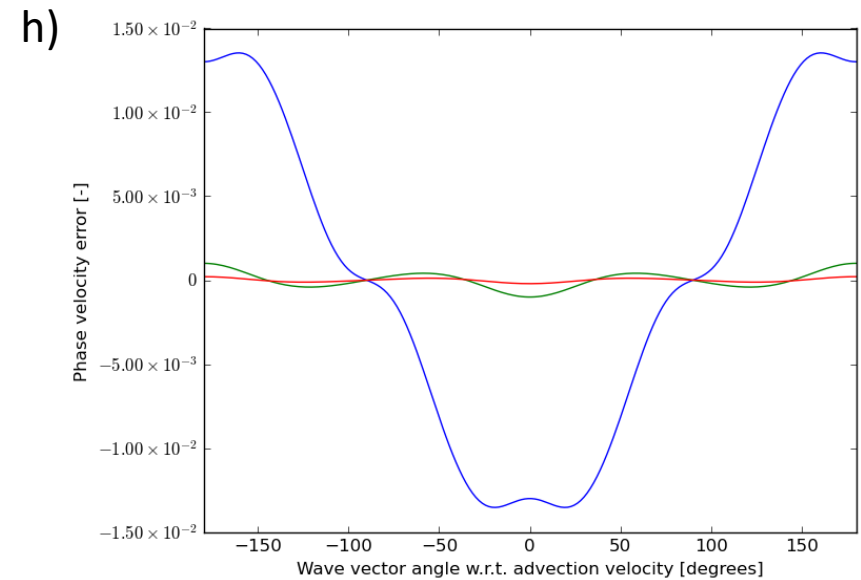

Fig. 7 shows the wave propagation characteristics for curl-preserving fourth order WENO-like schemes. Figs. 7a to 7d show one minus the absolute value of the amplification factor when the velocity vector makes angles of 0º, 15º, 30º and 45º relative to the x-direction of the 2D mesh. Figs. 7e to 7h show the phase error, again for the same angles. The 2D wave vector can make any angle relative to the 2D direction of velocity propagation, therefore, the amplitude and phase information are shown w.r.t. the angle made between the velocity direction and the direction of the wave vector. In each plot, the blue curve refers to waves that span 5 cells per wavelength; the green curve refers to waves that span 10 cells per wavelength; the red curve refers to waves that span 15 waves per wavelength.

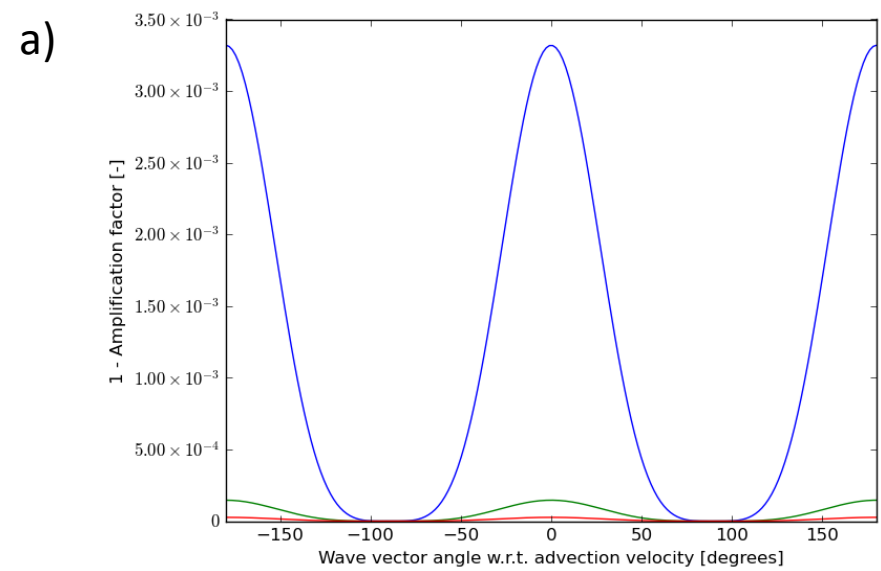
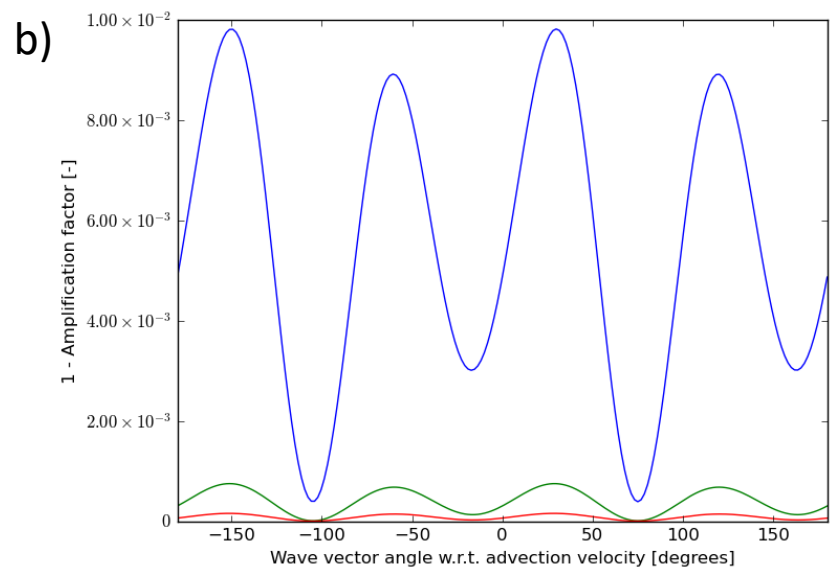
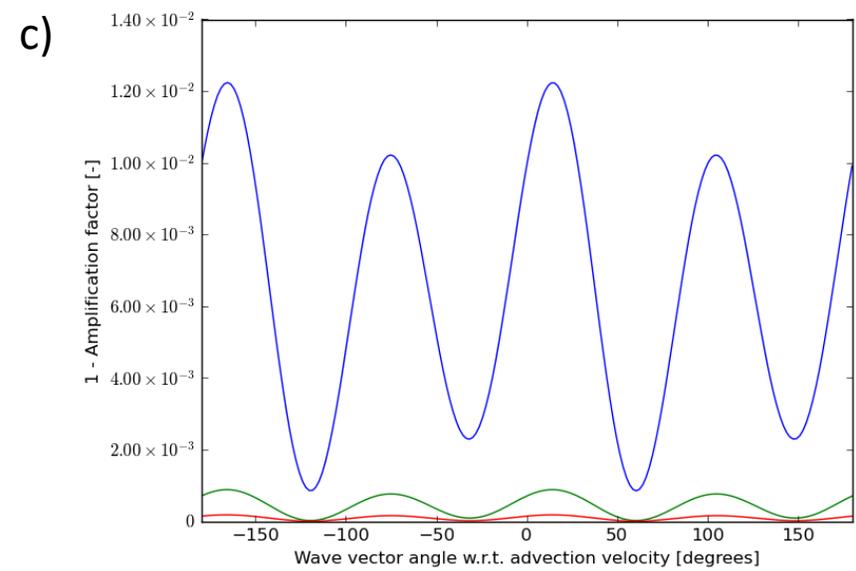
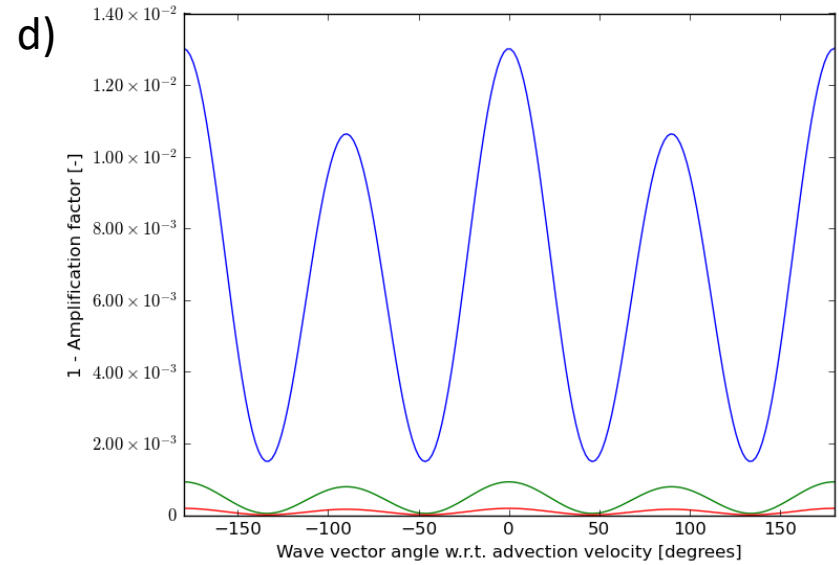

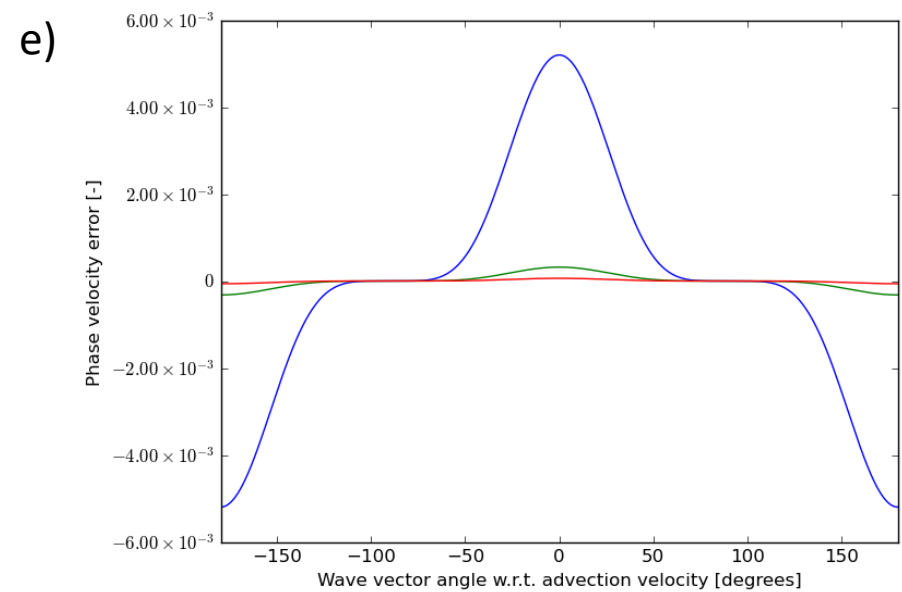
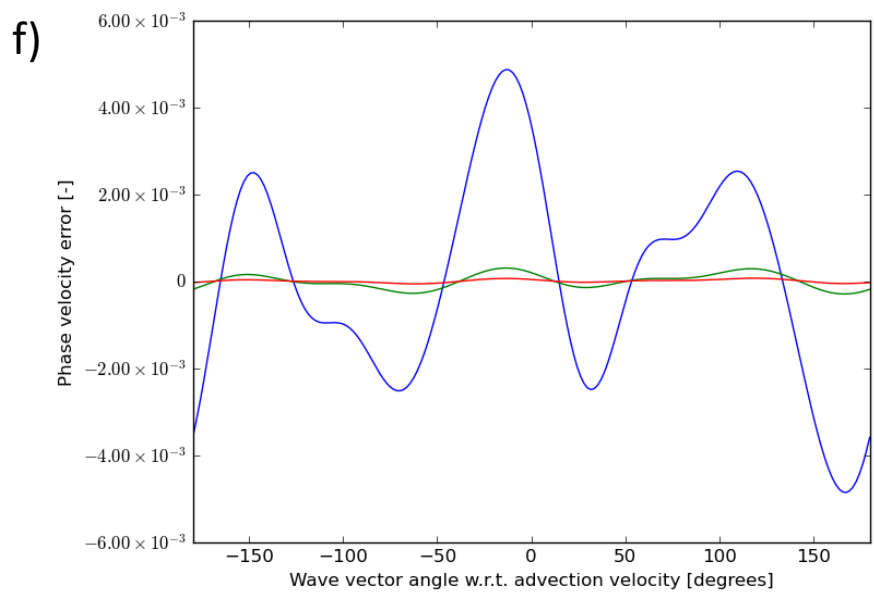
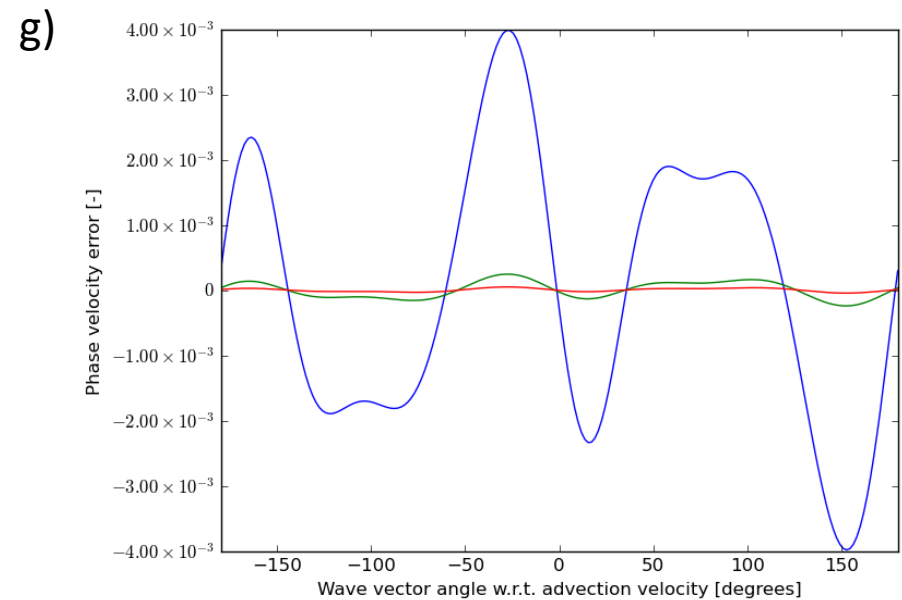
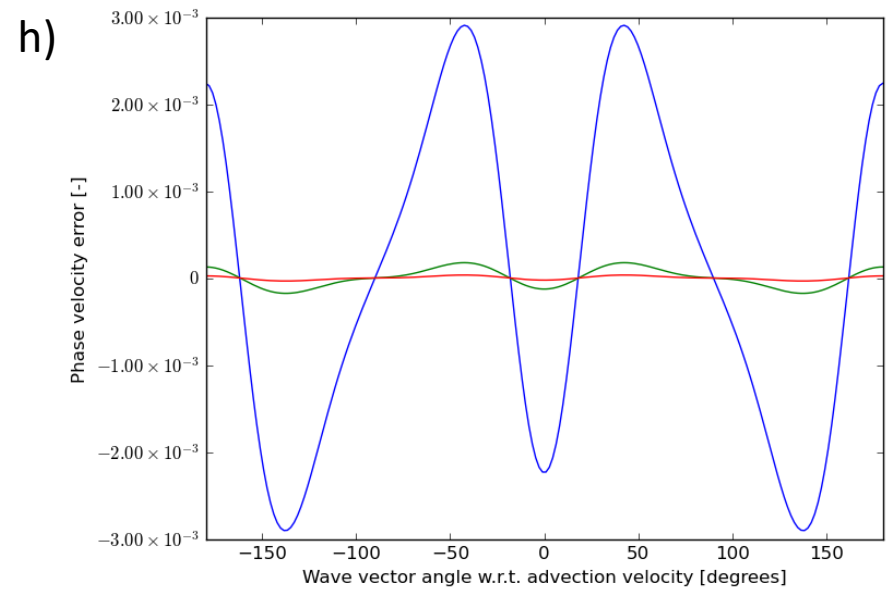

*Fig. 8 shows the wave propagation characteristics for curl-preserving third order P1P2-like schemes. Figs. 8a to 8d show one minus the absolute value of the amplification factor when the velocity vector makes angles of 0º, 15º, 30º and 45º relative to the x-direction of the 2D mesh. Figs. 8e to 8h show the phase error, again for the same angles. The 2D wave vector can make any angle relative to the 2D direction of velocity propagation, therefore, the amplitude and phase information are shown w.r.t. the angle made between the velocity direction and the direction of the wave vector. In each plot, the blue curve refers to waves that span 5 cells per wavelength; the green curve refers to waves that span 10 cells per wavelength; the red curve refers to waves that span 15 waves per wavelength.*

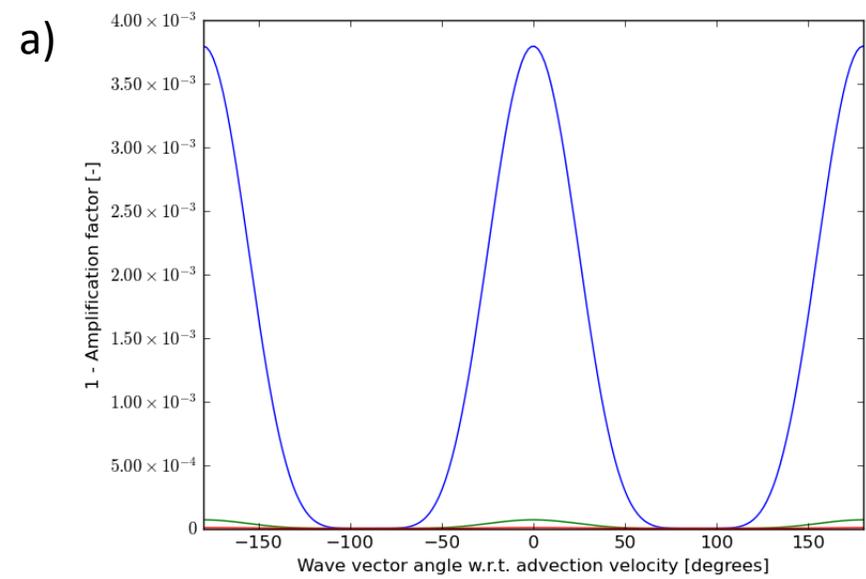 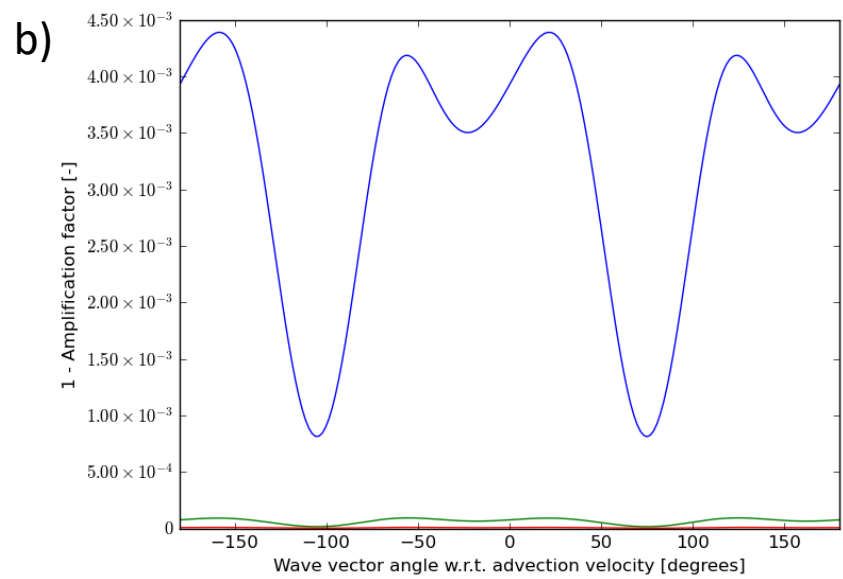

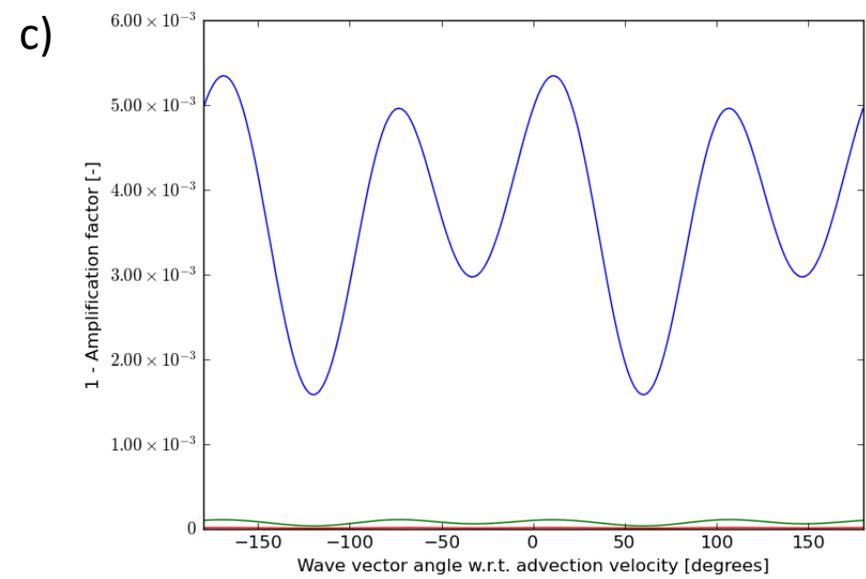 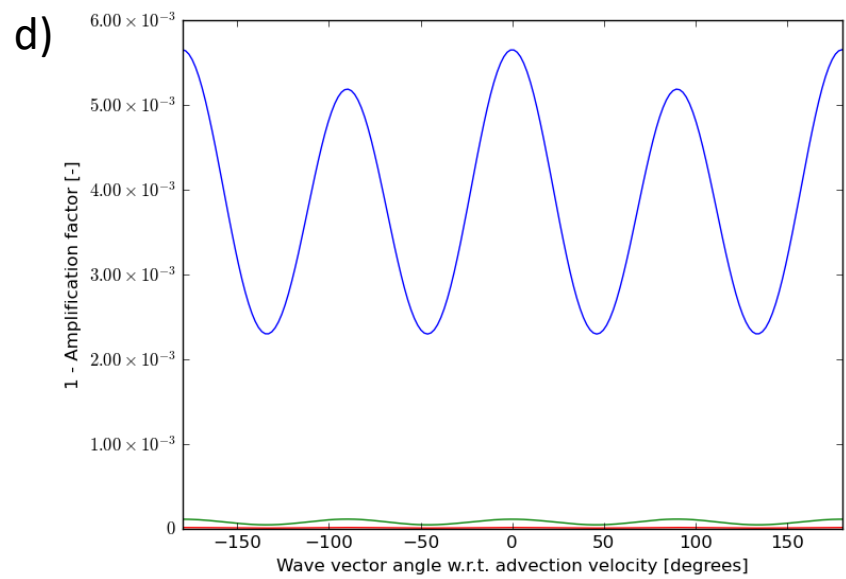

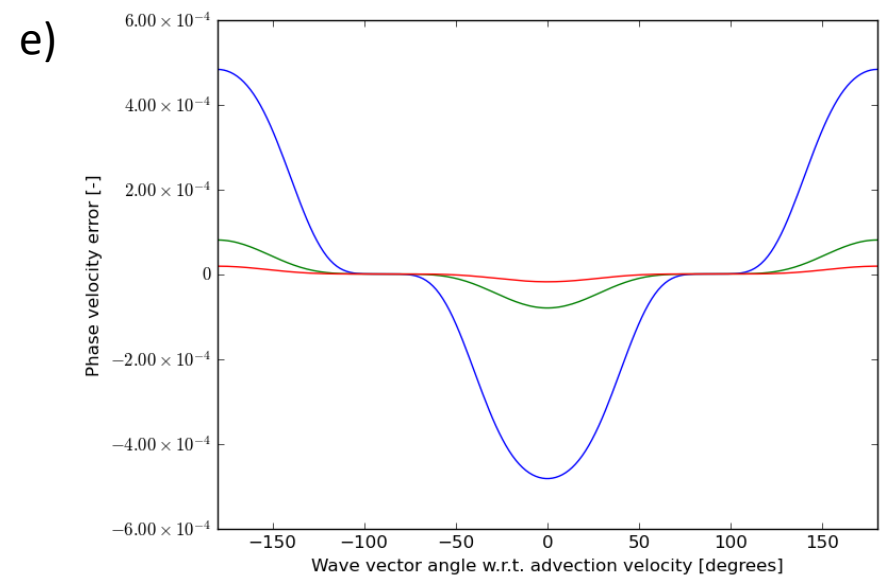
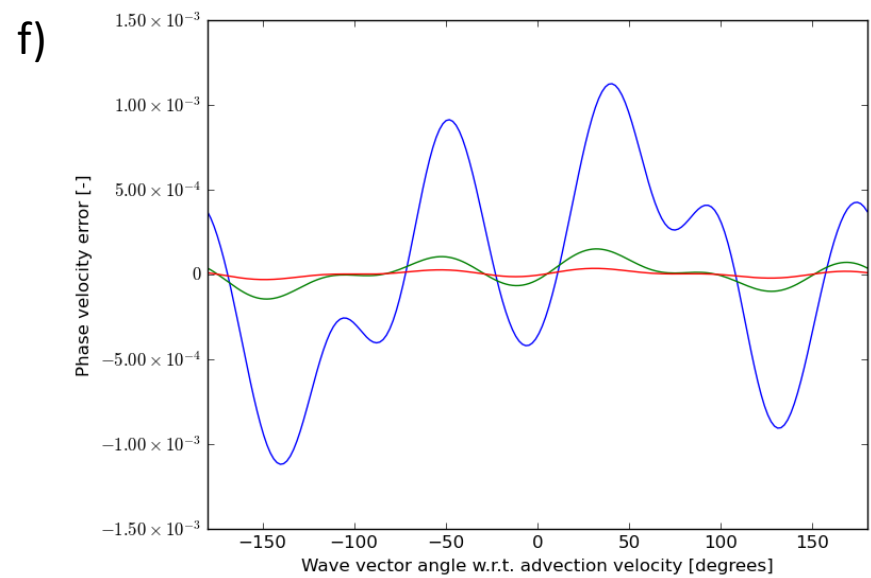
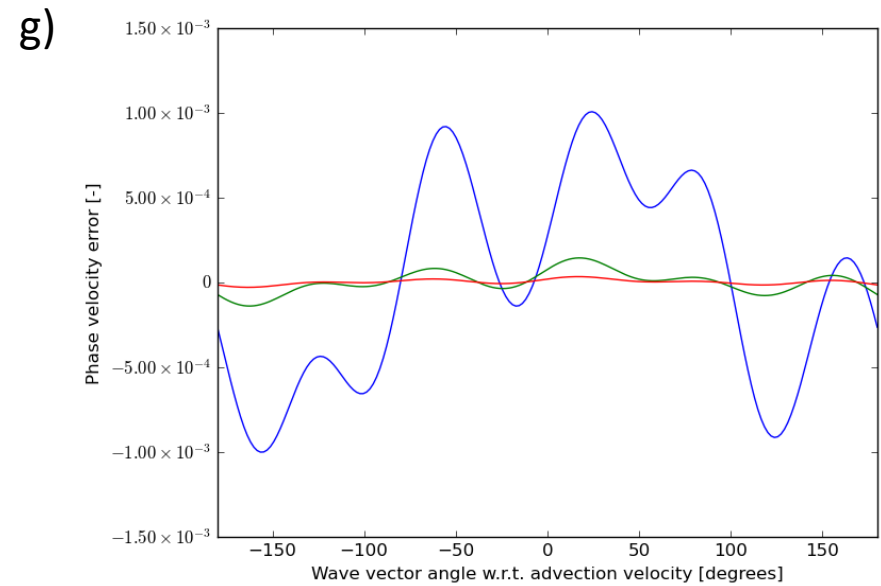
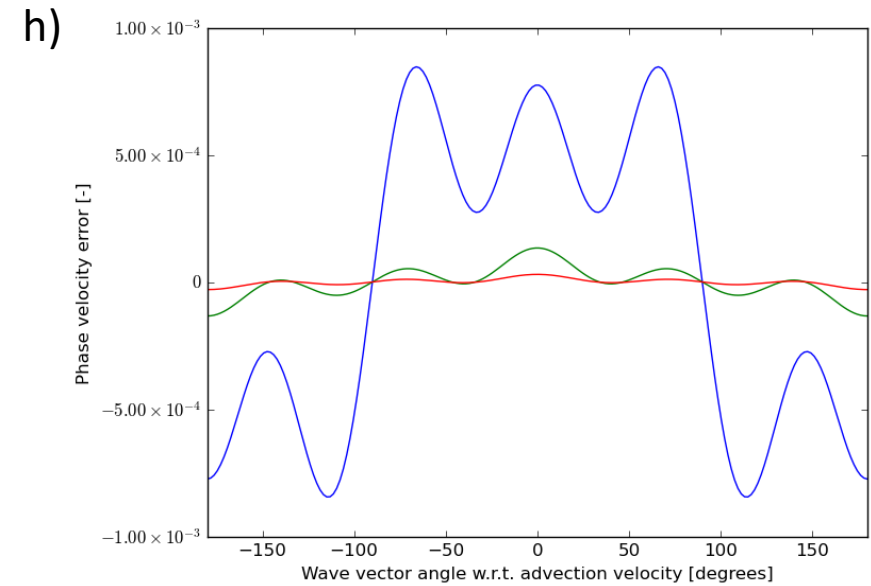

*Fig. 9 shows the wave propagation characteristics for curl-preserving fourth order P1P3-like schemes. Figs. 9a to 9d show one minus the absolute value of the amplification factor when the velocity vector makes angles of 0º , 15º , 30º and 45º relative to the x-direction of the 2D mesh. Figs. 9e to 9h show the phase error, again for the same angles. The 2D wave vector can make any angle relative to the 2D direction of velocity propagation, therefore, the amplitude and phase information are shown w.r.t. the angle made between the velocity direction and the direction of the wave vector. In each plot, the blue curve refers to waves that span 5 cells per wavelength; the green curve refers to waves that span 10 cells per wavelength; the red curve refers to waves that span 15 waves per wavelength.*

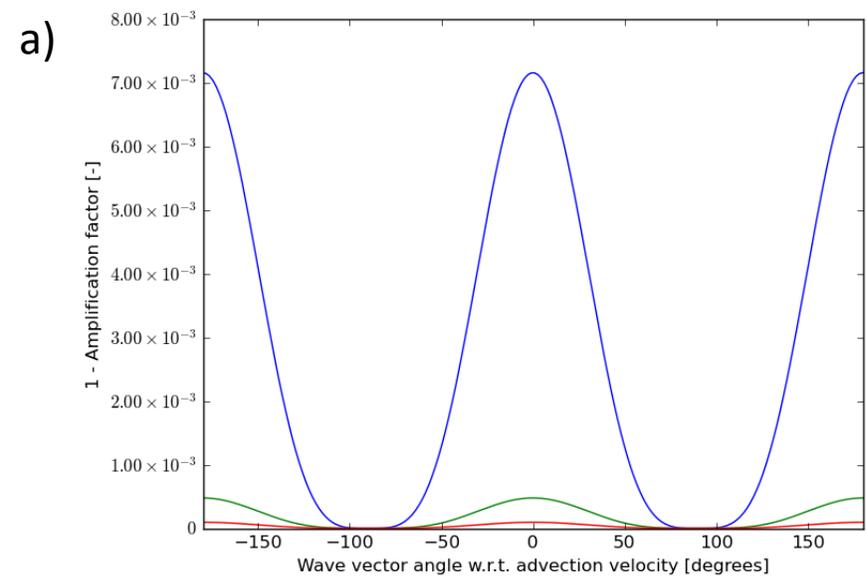 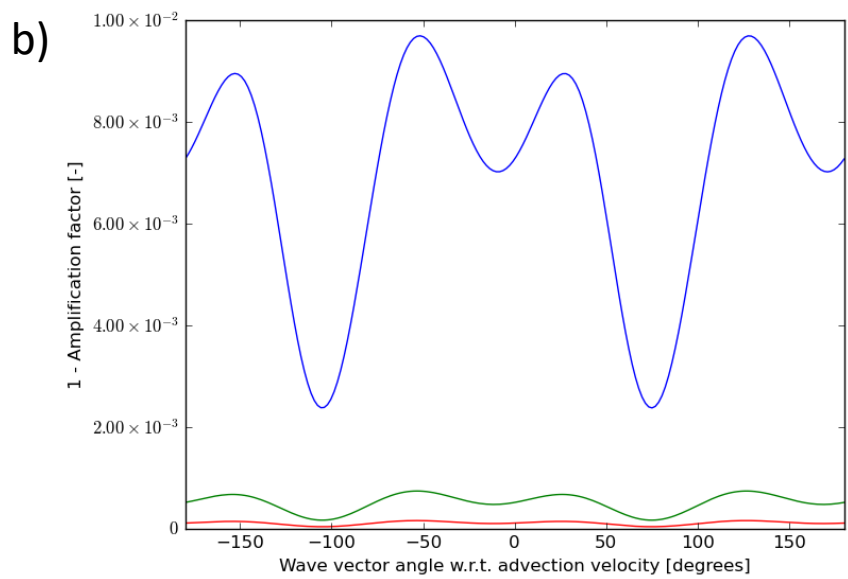

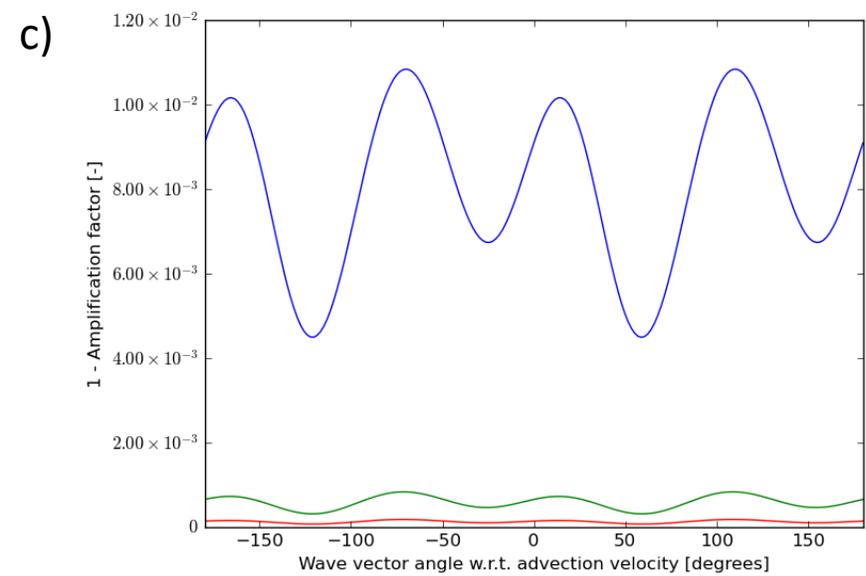 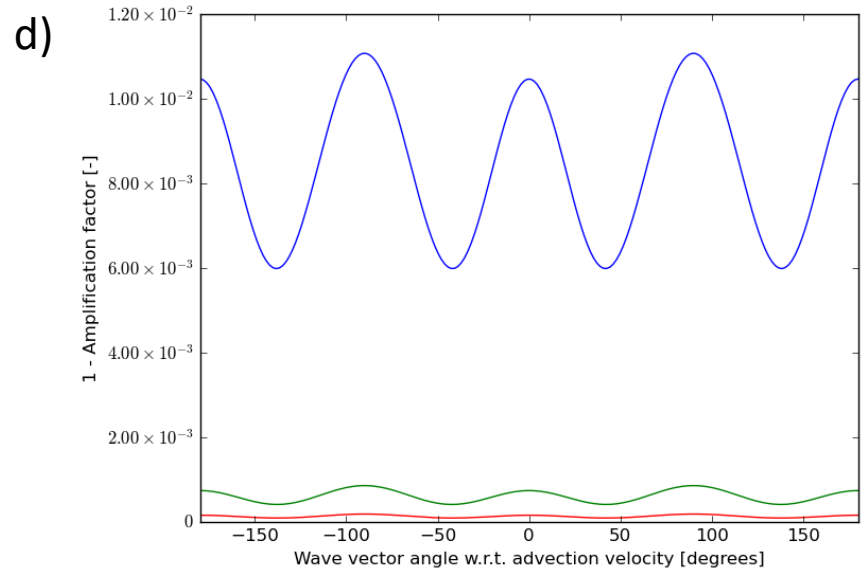

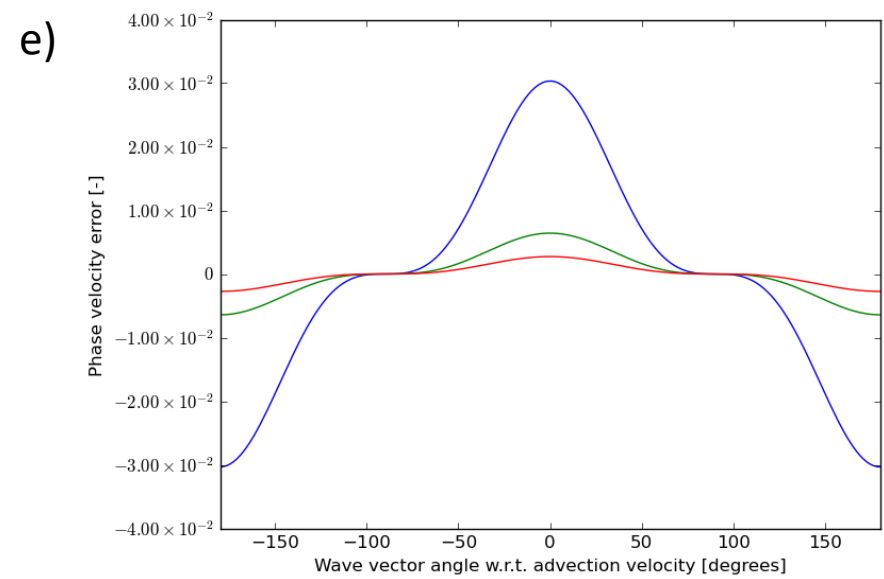
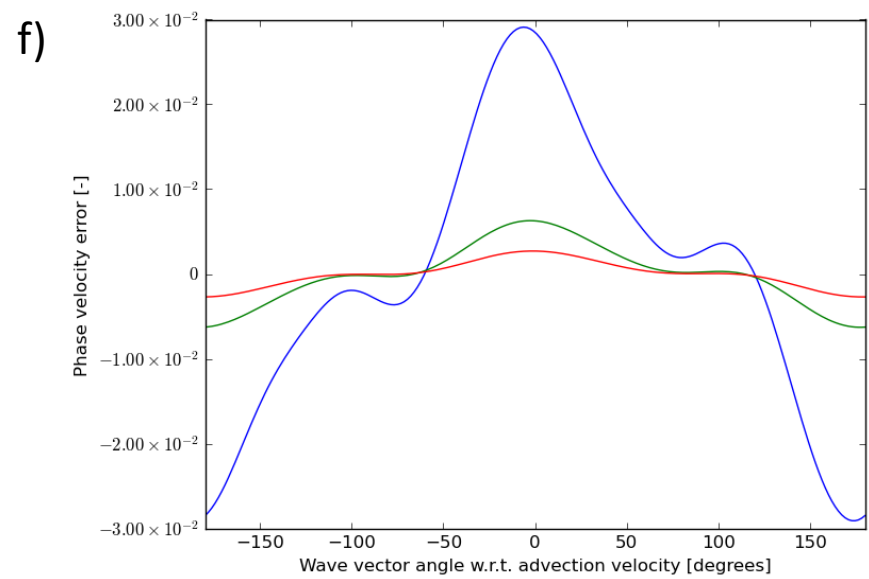
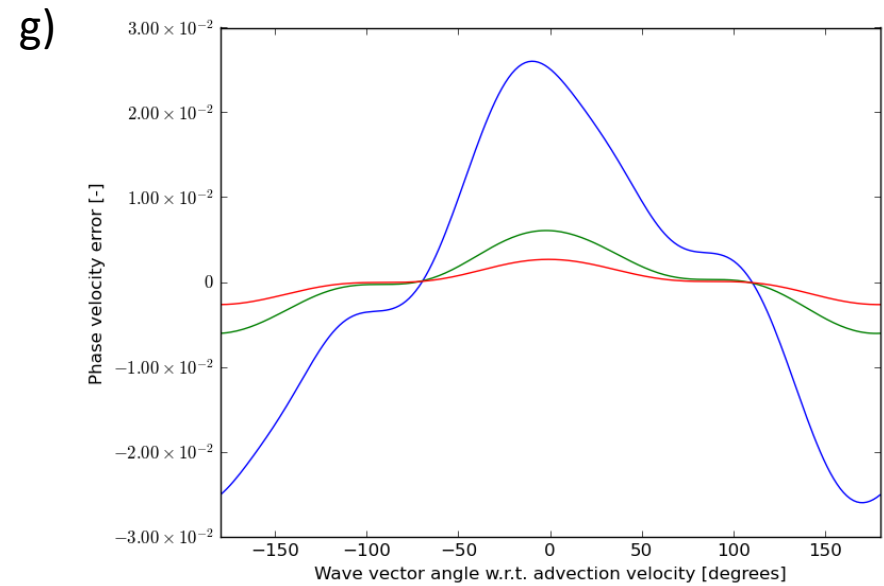
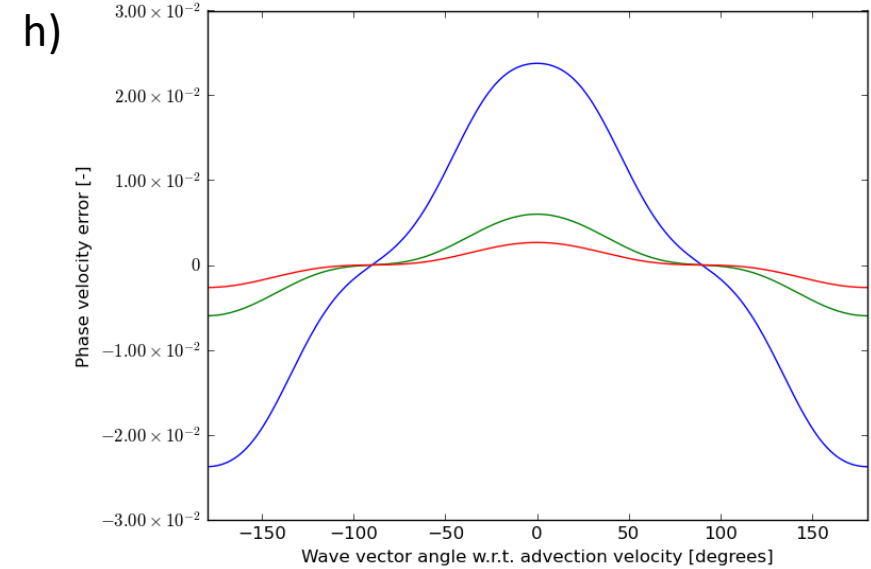

*Fig. 10 shows the wave propagation characteristics for curl-preserving second order DG-like schemes. Figs. 10a to 10d show one minus the absolute value of the amplification factor when the velocity vector makes angles of $0^o$, $15^o$, $30^o$ and $45^o$ relative to the x-direction of the 2D mesh. Figs. 10e to 10h show the phase error, again for the same angles. The 2D wave vector can make any angle relative to the 2D direction of velocity propagation, therefore, the amplitude and phase information are shown w.r.t. the angle made between the velocity direction and the direction of the wave vector. In each plot, the blue curve refers to waves that span 5 cells per wavelength; the green curve refers to waves that span 10 cells per wavelength; the red curve refers to waves that span 15 waves per wavelength.*

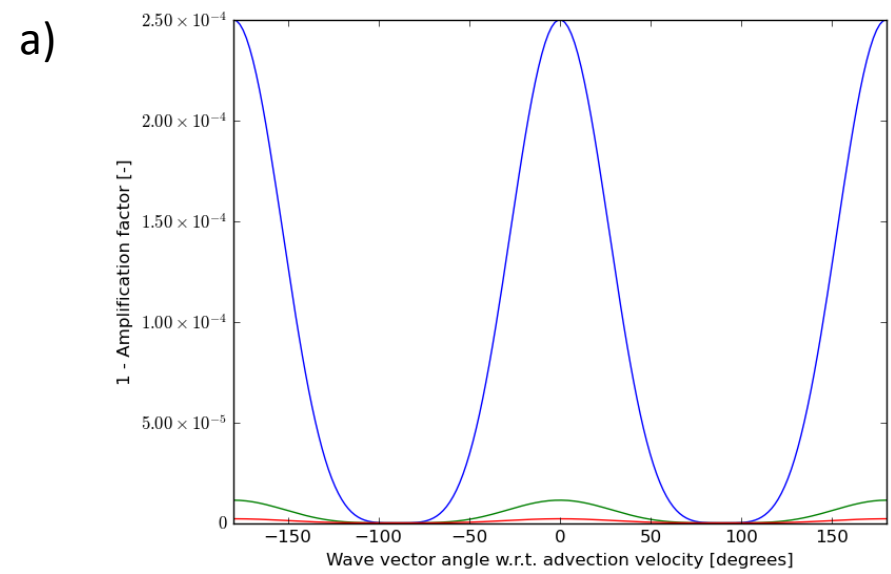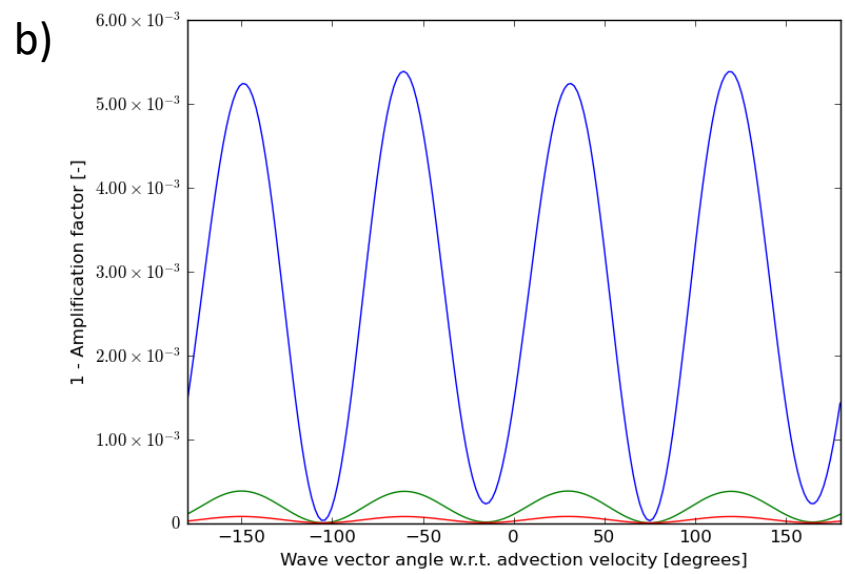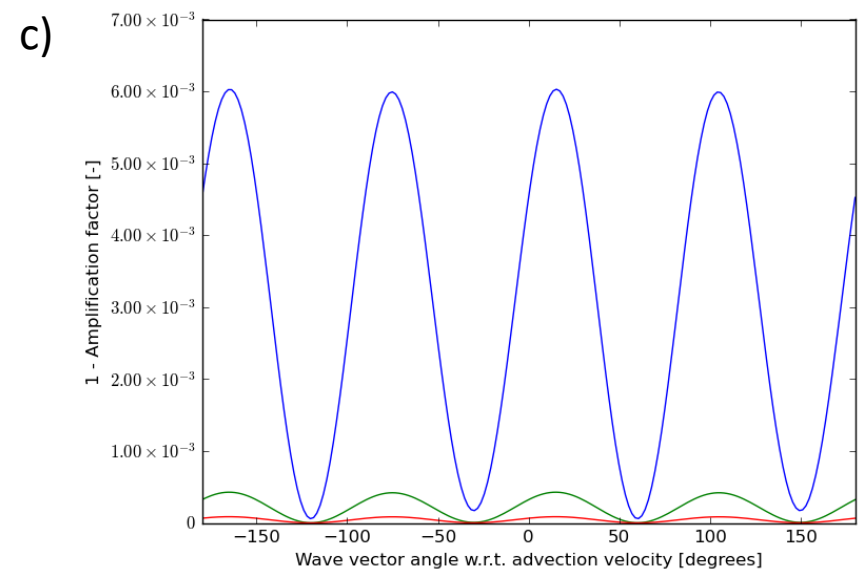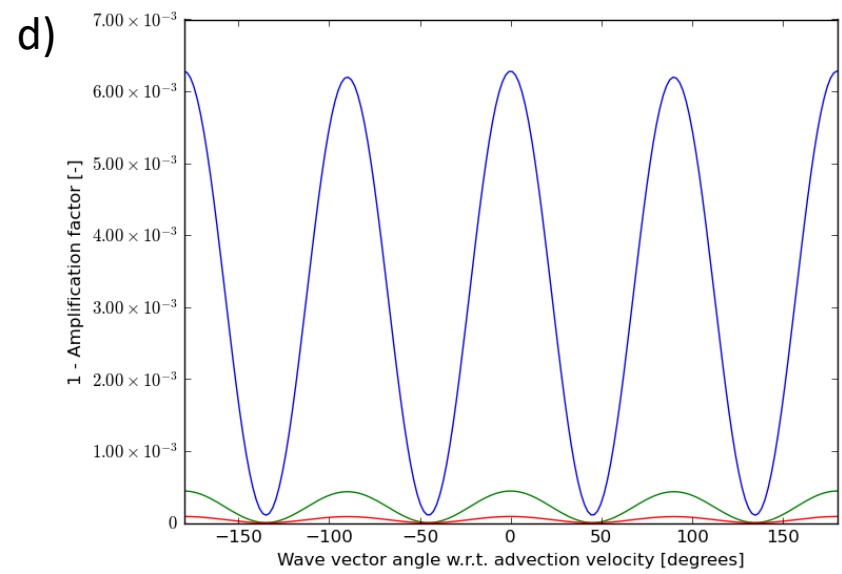

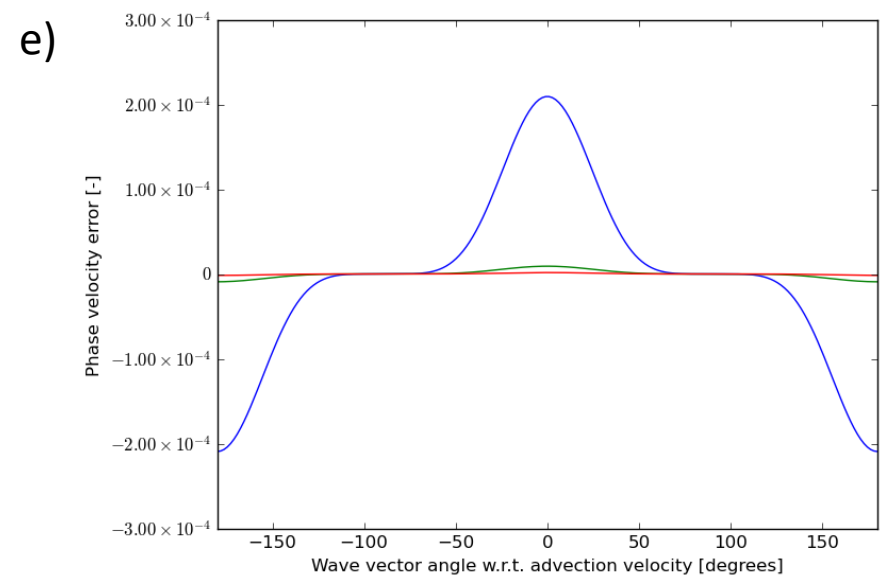
e)

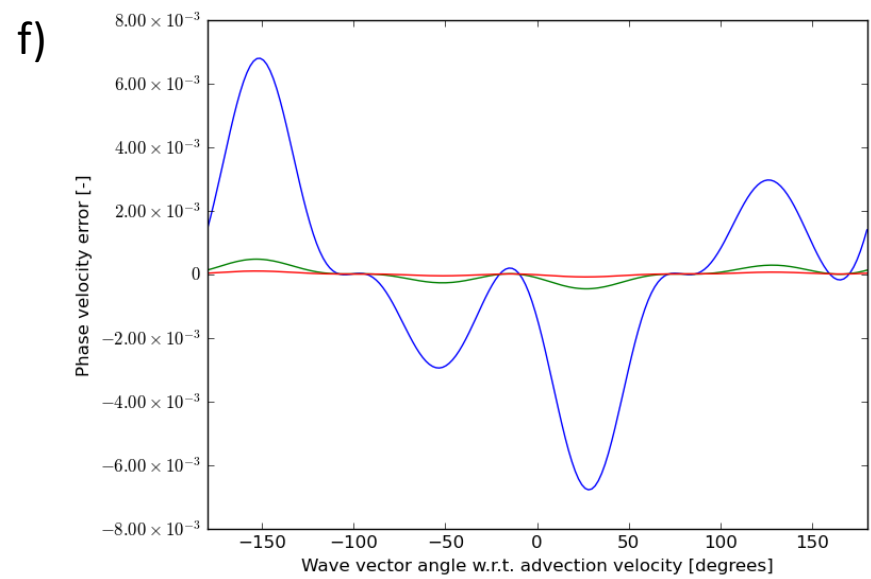
f)

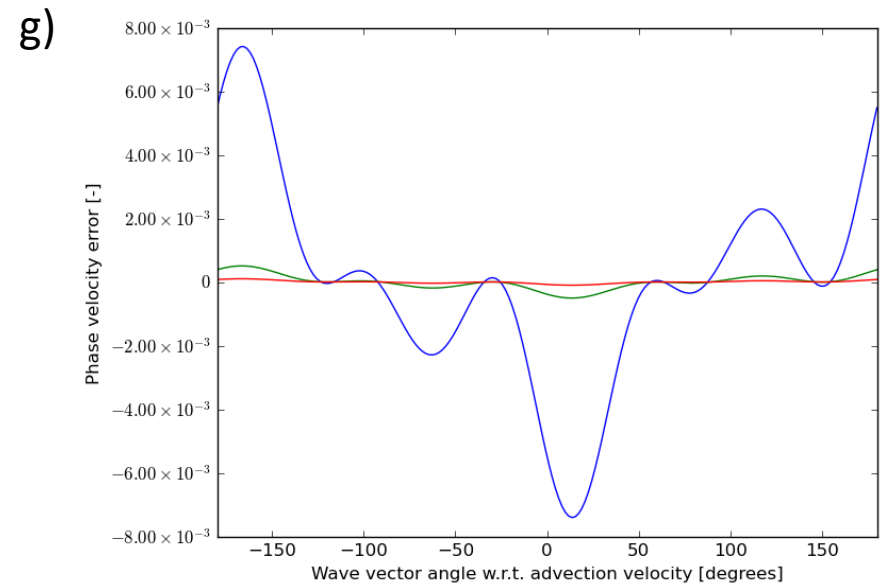
g)

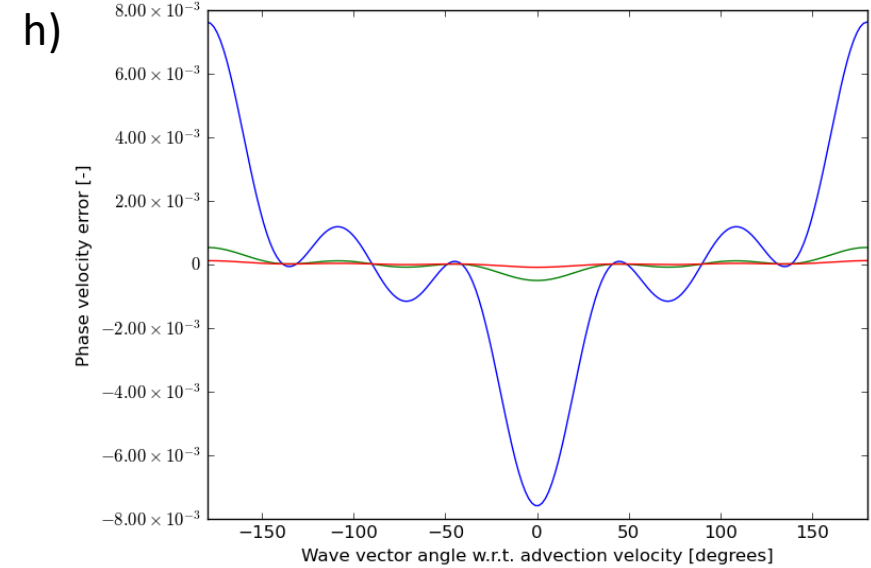
h)

*Fig. 11 shows the wave propagation characteristics for curl-preserving third order DG-like schemes. Figs. 11a to 11d show one minus the absolute value of the amplification factor when the velocity vector makes angles of 0º, 15º, 30º and 45º relative to the x-direction of the 2D mesh. Figs. 11e to 11h show the phase error, again for the same angles. The 2D wave vector can make any angle relative to the 2D direction of velocity propagation, therefore, the amplitude and phase information are shown w.r.t. the angle made between the velocity direction and the direction of the wave vector. In each plot, the blue curve refers to waves that span 5 cells per wavelength; the green curve refers to waves that span 10 cells per wavelength; the red curve refers to waves that span 15 waves per wavelength.*

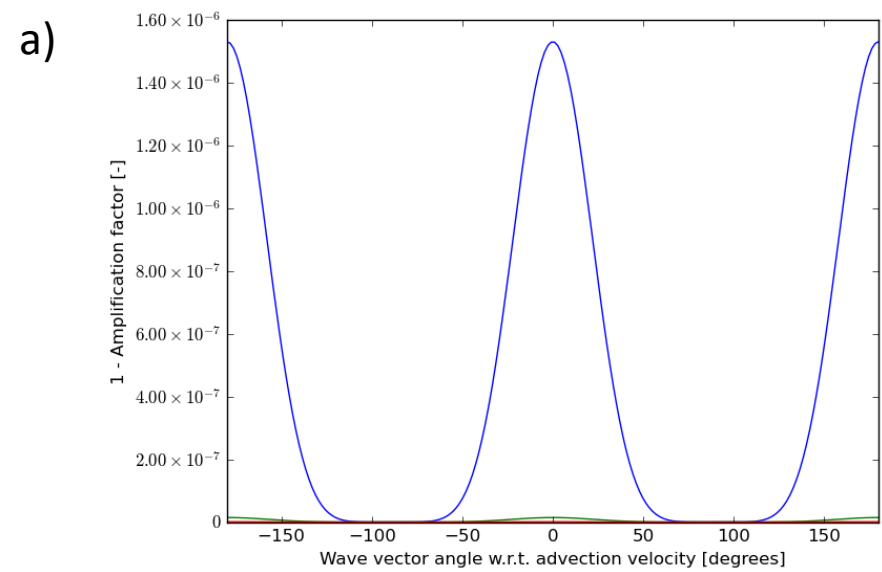 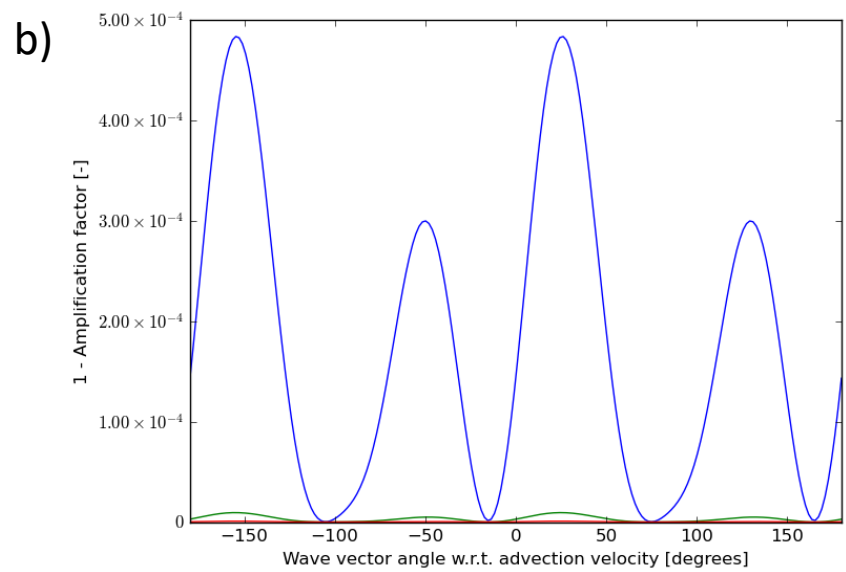

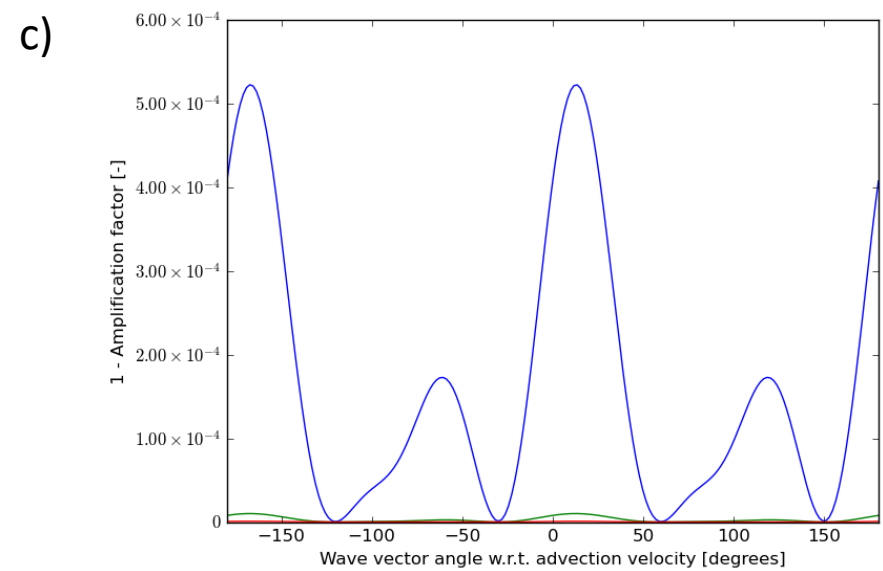 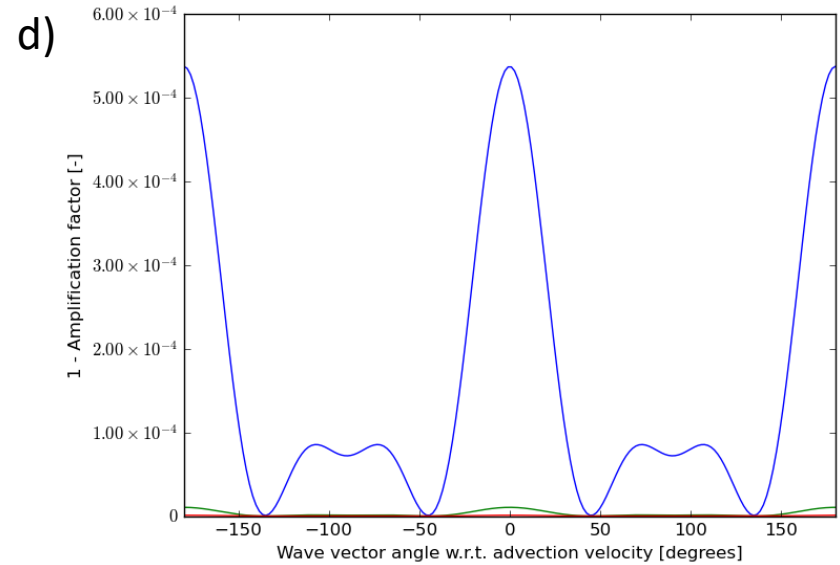

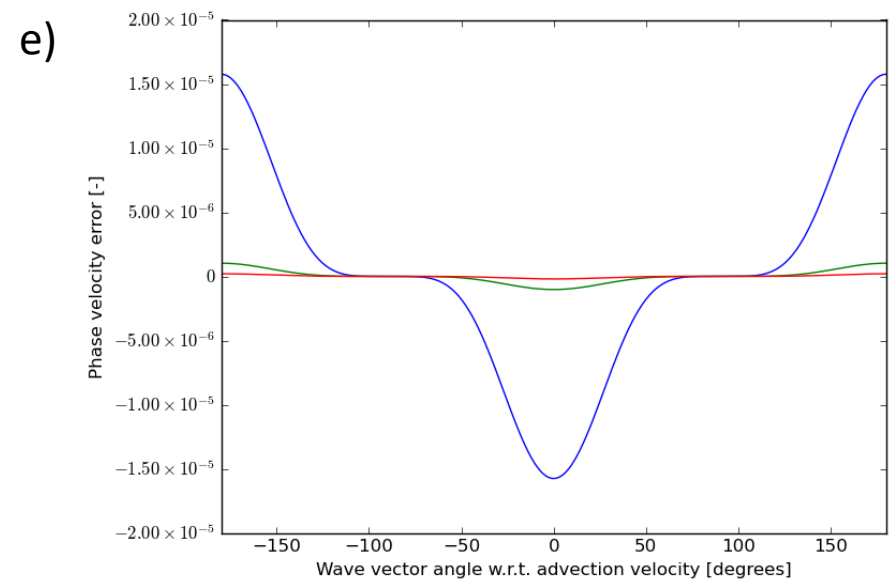
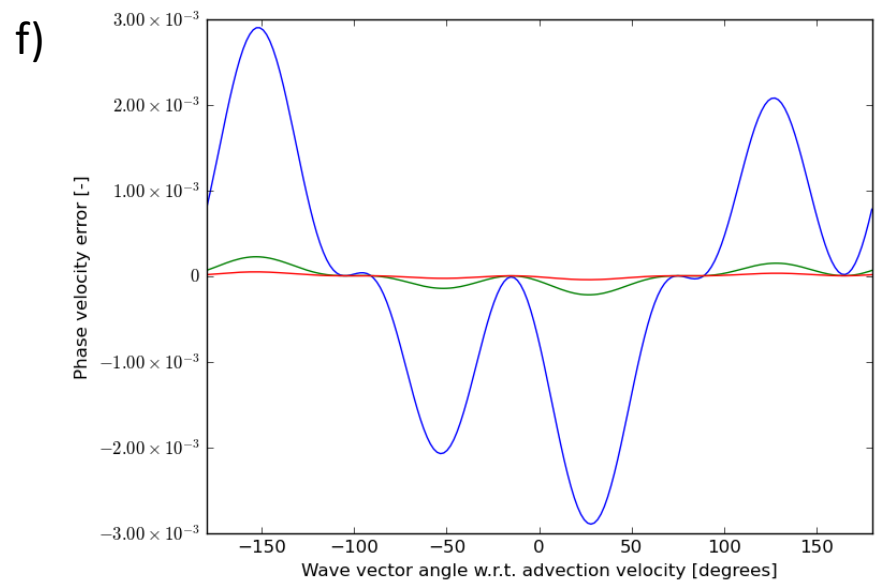
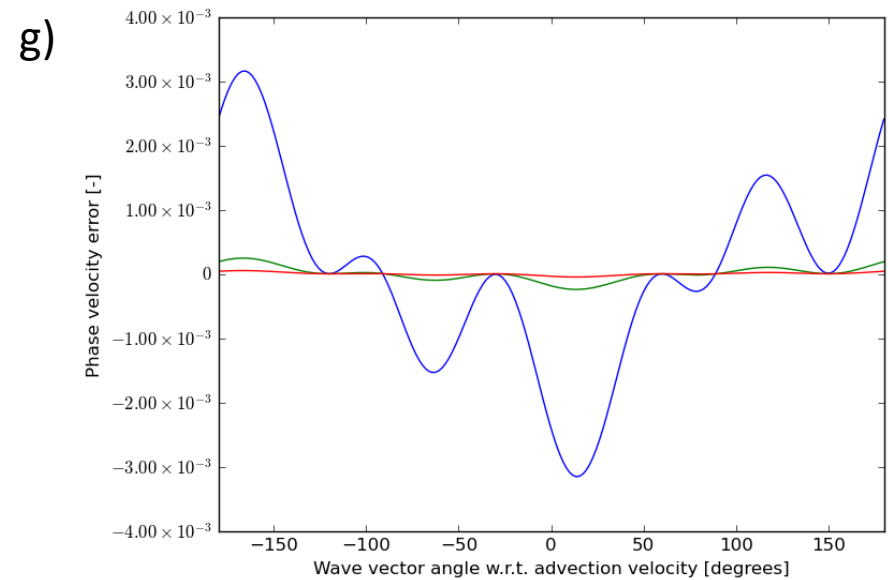
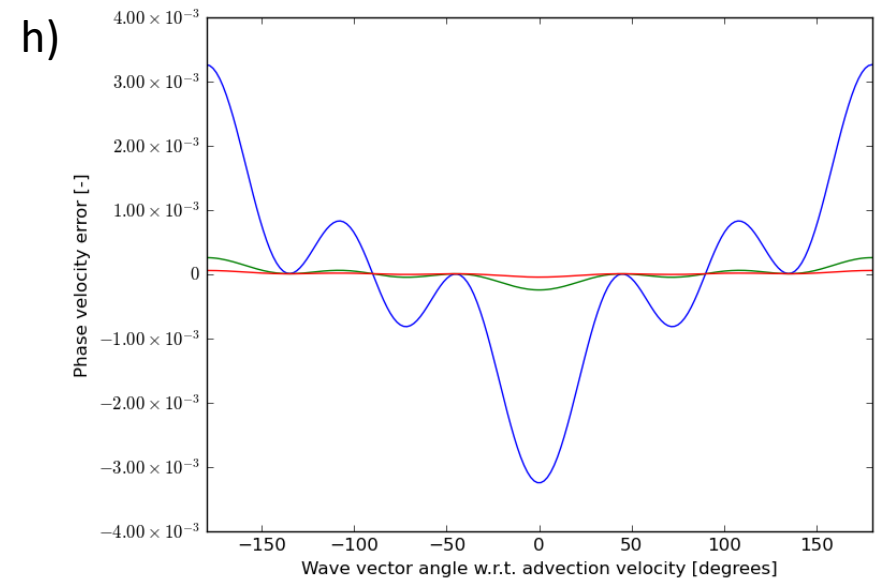

Fig. 12 shows the wave propagation characteristics for curl-preserving fourth order DG-like schemes. Figs. 12a to 12d show one minus the absolute value of the amplification factor when the velocity vector makes angles of 0º , 15 º , 30 º and 45 º relative to the x-direction of the 2D mesh. Figs. 12e to 12h show the phase error, again for the same angles. The 2D wave vector can make any angle relative to the 2D direction of velocity propagation, therefore, the amplitude and phase information are shown w.r.t. the angle made between the velocity direction and the direction of the wave vector. In each plot, the blue curve refers to waves that span 5 cells per wavelength; the green curve refers to waves that span 10 cells per wavelength; the red curve refers to waves that span 15 waves per wavelength.

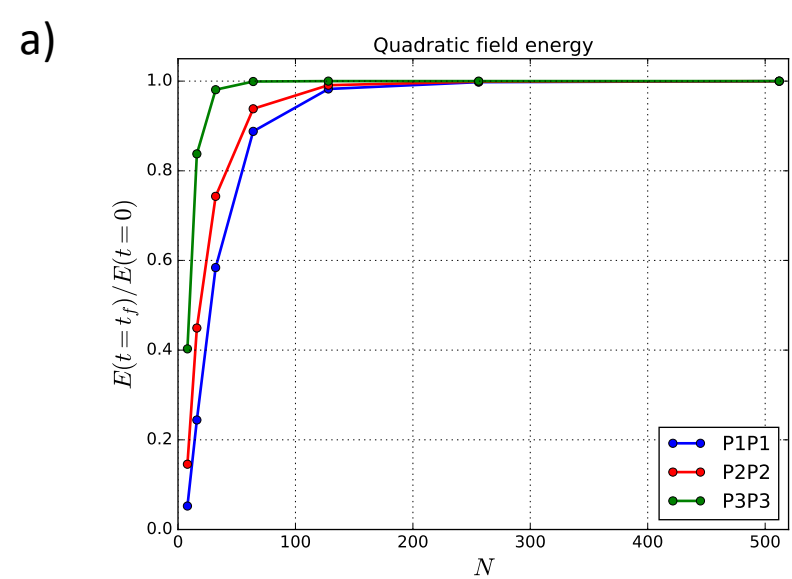
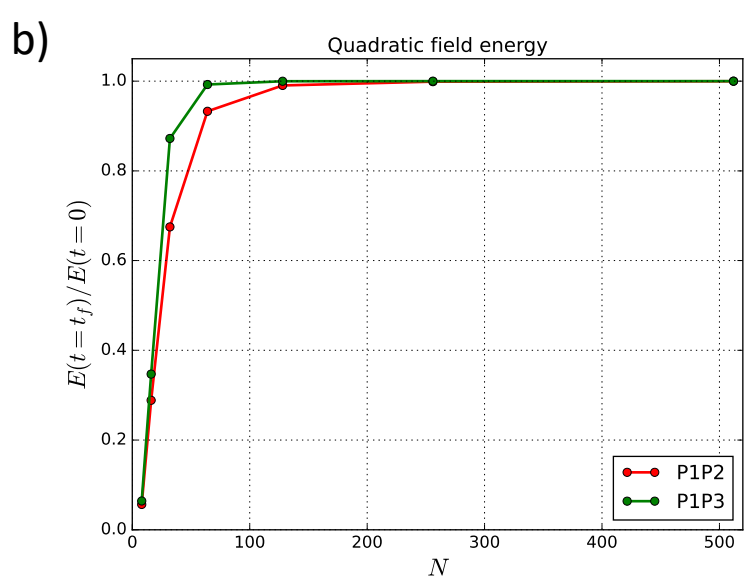

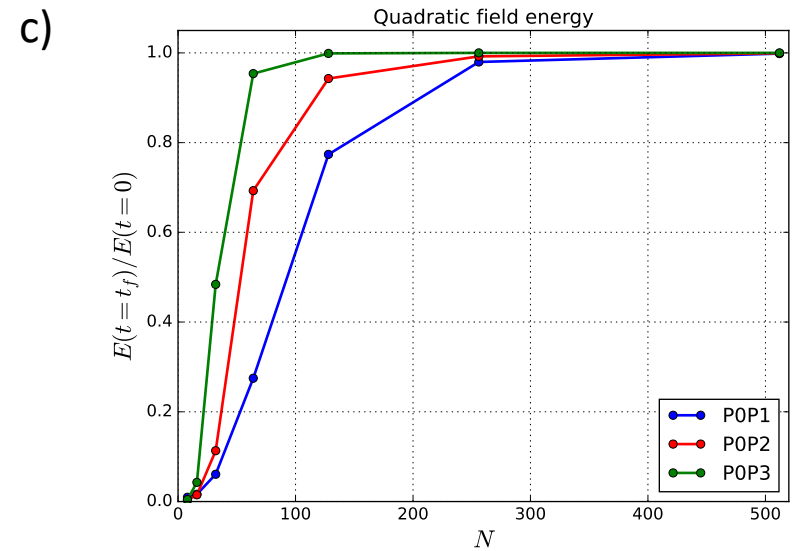

*Fig. 13 shows the quadratic field energy from the vortex problem that is preserved on the mesh at the final time in the simulation as a function of mesh size. Panel a) displays the curl-free DG-like schemes, panel b) displays the curl-free P1PN-like schemes and panel c) displays the curl-free WENO-like schemes.*

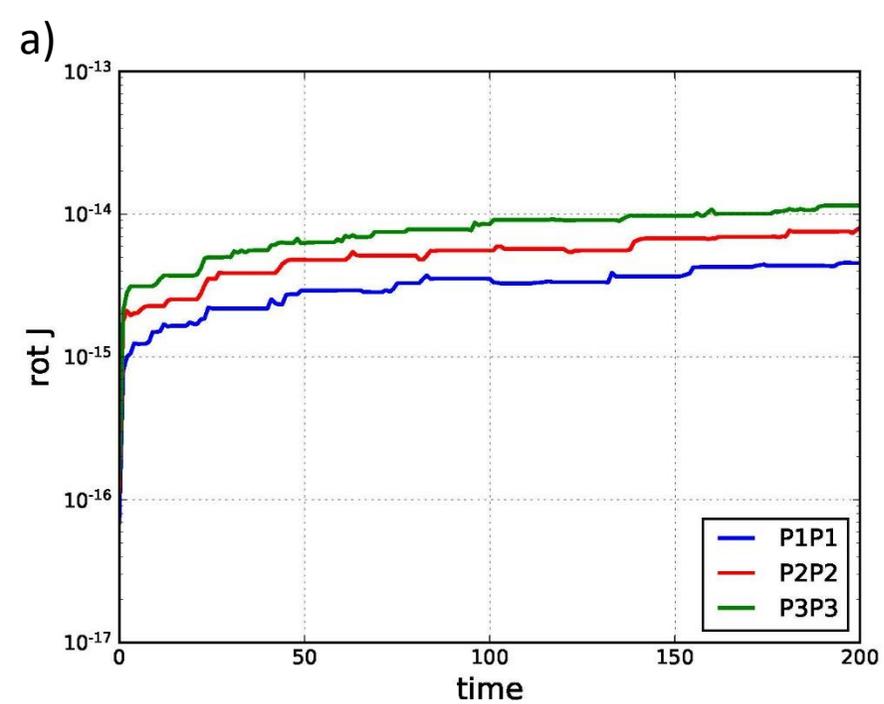
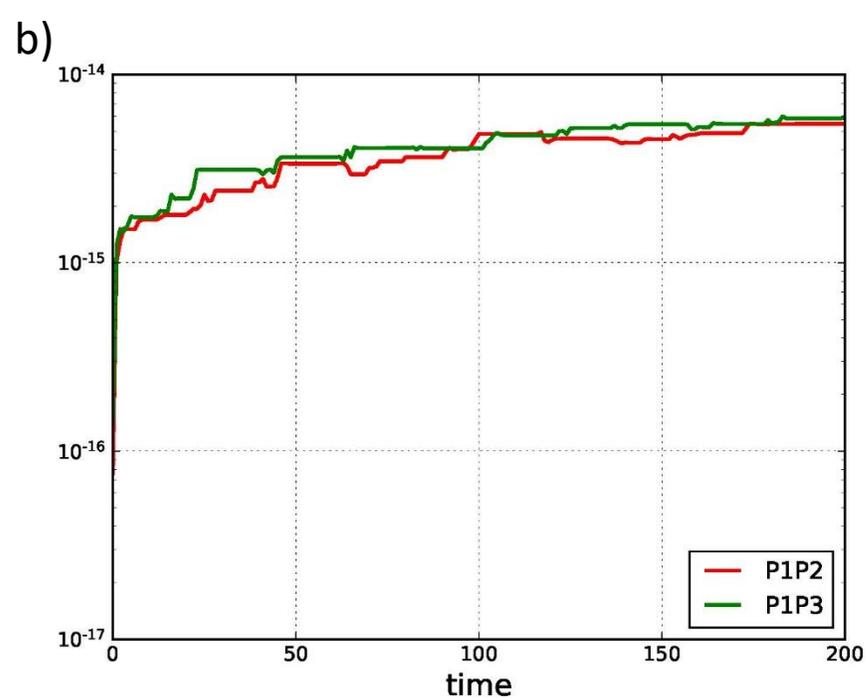
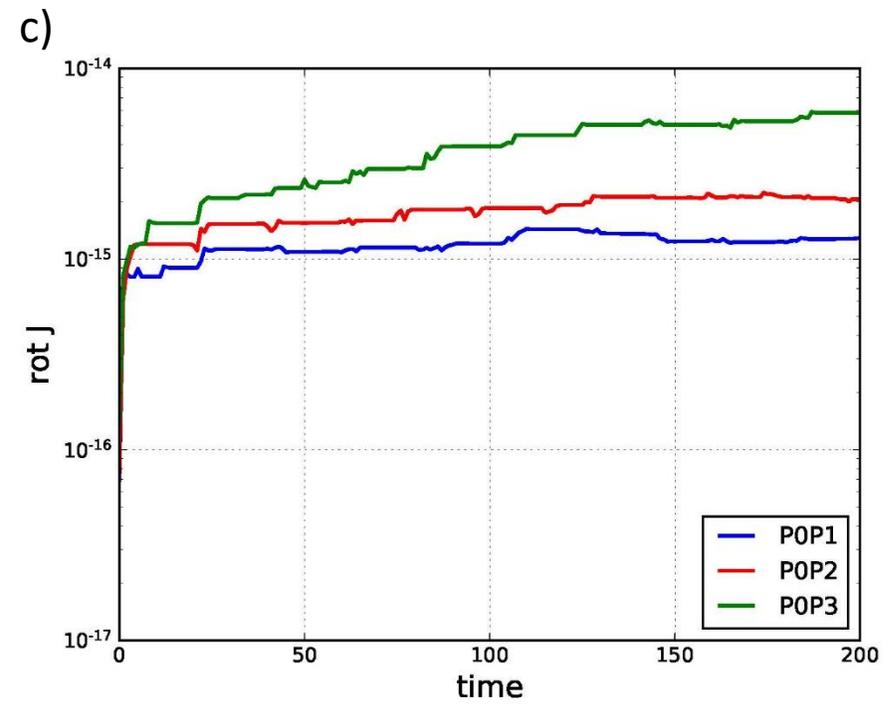

Fig. 14 shows the maximum pointwise error of the curl of **J** as a function of time for a 64×64 zone run of the vortex problem. Fig. 14a shows the evolution of the maximum pointwise curl as a function of time for the 2$^{nd}$, 3$^{rd}$ and 4$^{th}$ order curl-free DG-like schemes. Fig. 14b shows the evolution of the maximum pointwise curl as a function of time for the 3$^{rd}$ and 4$^{th}$ order curl-free P1PN-like schemes. Fig. 14c shows the evolution of the maximum pointwise curl as a function of time for the 2$^{nd}$, 3$^{rd}$ and 4$^{th}$ order curl-free WENO-like schemes. The figure shows that all our curl-preserving schemes can preserve the curl constraint up to machine accuracy.